\documentclass{article}
\newcommand{\pre}{{\noindent\noindent{} \bf Preuve. }}
\newcommand{\di}{\displaystyle}
\def\sk{\vskip 0.2cm}
\def\Z{\mathbf Z}
\def\N{{\mathbf N}}
\def\squareforqed{\hbox{\rlap{$\sqcap$}$\sqcup$}}
\def\qed{\ifmmode\else\unskip\quad\fi\squareforqed}
\def\smartqed{\def\qed{\ifmmode\squareforqed\else{\unskip\nobreak\hfil
\penalty50\hskip1em\null\nobreak\hfil\squareforqed
\parfillskip=0pt\finalhyphendemerits=0\endgraf}\fi}}
\def\smartqed{\def\qed{\ifmmode\squareforqed\else{\unskip\nobreak\hfil
\penalty50\hskip1em\null\nobreak\hfil\squareforqed
\parfillskip=0pt\finalhyphendemerits=0\endgraf}\fi}}
\newcommand{\cqfd}{\qed}

\newcommand{\Ga}{\Gamma}
\newcommand{\De}{\Delta}

\newcommand{\al}{\alpha}

\newcommand{\La}{\Lambda}
\newcommand{\be}{\beta}
\newcommand{\ep}{\varepsilon}

\newcommand{\Om}{\Omega}

\newcommand{\te}{\theta}
\newcommand{\ze}{\zeta}
\newcommand{\f}{\mathbf{\mathcal{F}}}
\newtheorem{thm}{Th\'eor\`eme}[section]
\newtheorem{cor}[thm]{Corollaire}
\newtheorem{lem}[thm]{Lemme}
\newtheorem{pro}[thm]{Proposition}

\newtheorem{exe}[thm]{Exemple}

\newtheorem{df}{D\'efinition}[section]
\newtheorem{rem}{Remarque}[section]

\title{Extensions absolument $lq$-finies\\}
\author{El Hassane Fliouet}
\date{}
\begin{document}
\maketitle
\markboth{E. Fliouet}{Extensions absolument $lq$-finies}
\renewcommand{\sectionmark}[1]{}
\begin{abstract}

Let $K/k$ be purely inseparable extension of characteristic $p>0$. By invariants, we characterize the measure of the size of $K/k$. In particular, we give a necessary and sufficient condition that $K/k$ is of bounded size. Furthermore, in this note, we continue to be interested in the relationship that connects  the restricted distribution of  finitude at the local level of intermediate fields of a purely inseparable extension $K/k$ to the absolute or global finitude of $K/k$. Part of this problem was treated successively by J.K Devney, and in my work with M. Chellali. The other part is the subject of this paper, it is a question of describing the absolutely $lq$-finite extensions.
Among others, any absolutely $lq$-finite extension decomposes into  $w_0$-genera\-t\-ed  extensions. However, we construct an example of extension of infinite size such that for any intermediate field $L$ of $K/k$, $L$ is of finite size over $k$. In addition, $K/k$ does not respect the  distribution of horizontal finitude.
\end{abstract}

{\bf Mathematics Subject Classification MSC2010:} Primary 12F15
\sk

\section{Introduction}

Soit $K/k$ une extension purement ins\'{e}parable de caract\'{e}ristique $p>0$. Une partie $B$ de $K$ est dite $r$-base de $K/k$ si $K=k(K^p)(B)$ et pour tout $x\in B$, $x\not\in k(K^p)(B\setminus \{x\})$. En vertu du (\cite{N.B}, III,  p. 49, corollaire {3}) et de la propri\'{e}t\'{e} d'\'{e}change des $r$-bases,  on en d\'{e}duit que toute extension admet une $r$-base et que le cardinal d'une $r$-base  est invariant. Si de plus $K/k$ est d'exposant fini, on v\'{e}rifie aussit\^{o}t que $B$ est une $r$-base de $K/k$ si et seulement si $B$ est un g\'{e}n\'{e}rateur minimal de $K/k$. Par suite, nous pouvons
contr\^{o}ler  la taille de toute extension purement ins\'{e}parable d'exposant born\'{e} $K/k$ et la longueur de tout corps $k$ par l'interm\'{e}diaire du degr\'{e} d'irrationalit\'{e} de $K/k$ et le degr\'{e} d'imperfection de $k$  d\'{e}finis respectivement par $di(K/k)=|G|$ o\`{u} $G$ est un g\'{e}n\'{e}rateur minimal de $K/k$, et $di(k)=|B|$, o\`{u} $B$ est une $r$-base de $k/k^p$.  Notamment, on montre que la mesure de la taille d'une extension est une fonction croissante par rapport \`{a} l'inclusion. Plus pr\'{e}cis\'{e}ment, pour toute chaine d'extensions d'exposant born\'{e} $k\subseteq L\subseteq K$, on a $di(L/k)\leq di(K/k)$. Par ailleurs, cette propri\'et\'e s'est av\'{e}r\'{e}e fort utile pour \'etendre la meseure de la taille  \`a une extension purement ins\'{e}parable $K/k$ quelconque par prolongement vertical du degr\'{e} d'irrationalit\'{e} des sous-extensions interm\'{e}diaires d'exposant fini de $K/k$. Ainsi,  on pose $di(K/k)=\di\sup_{n\in\N}(di(k^{p^{-n}}\cap K /k))$, ici le $\sup$ est utilis\'e au sens du (\cite{N.B}, III, p. 25, proposition {2}). Dans ce contexte, on montre que la mesure de la taille d'une extension est compatible avec l'inclusion et la lin\'{e}arit\'{e} disjointe. En d'autres termes,  on a :
\sk

\begin{itemize}{\it
\item[\rm{(i)}] Pour toute chaine d'extensions purement ins\'{e}parables $k\subseteq L\subseteq L'\subseteq K$, on a $di(L'/L)\leq di(K/k)\leq di(k)$
\item[\rm{(ii)}] Pour tous corps interm\'{e}diaires $K_1$ et $K_2$ de $K/k$, $k$-lin\'{e}airement disjoints, on a $di(K_1(K_2)/k)=di(K_1/k)+di(K_2/k)$ et $di(K_1(K_2)/K_2)=di(K_1/k)$.}
\end{itemize}
\sk

En outre, $di(K/k)=\sup(di(L/k))$, o\`{u} $k\subseteq L\subseteq K$. Autrement dit, la mesure de la taille de $K/k$ est vue comme limite inductive du degr\'e d'irrationalit\'e de ces sous-extensions interm\'{e}diaires.
Dans cette note, nous continuons \`{a} s'int\'{e}resser \`{a} la relation qui relit la r\'{e}partition restreinte de la finitude au niveau des corps interm\'{e}diaires d'une extension purement ins\'{e}parable $K/k$ \`{a} la finitude absolue ou globale de $K/k$.
\sk

Dans \cite{Dev2}, J. K. Deveney construit un exemple d'extension modulaire $K/k$ tel que toute sous-extension propre de $K/k$ est finie, et pour tout $n\in \N$, on a $[k^{p^{-n}}\cap K: k]= p^{2n}$. Par cons\'{e}quent, $K/k$ ne conserve pas la distribution de la finitude au niveau local des corps interm\'{e}diaires de $K/k$. Cependant, il est facile de v\'{e}rifier que
$di(k^{p^{-n}}\cap K/k))=2$, donc $K/k$ est relativement parfaite d'exposant non born\'{e} dont la mesure de la taille vaut 2, et par suite $K/k$ conserve la finitude horizontale (la taille) des corps interm\'{e}diaires. En tenant compte de cet exemple, dans \cite{Che-Fli2} on construit pour tout entier $j$ une extension purement ins\'{e}parable $K/k$ d'exposant non born\'{e}  v\'{e}rifiant :
\sk

\begin{itemize}{\it
\item[\rm{(i)}] Toute sous-extension propre de $K/k$ est finie.
\item[\rm{(ii)}] Pour tout $n\in \N$, $[k^{p^{-n}}\cap K: k]= p^{jn}$.}
\end{itemize}
\sk

Am\'{e}liorant ainsi le contre-exemple de J. K. Deveney, une telle extension est relativement parfaite d'exposant non born\'{e}, et pour tout $n\in \N$, $di(k^{p^{-n}}\cap K/ k)=j$. Il est \'{e}galement clair que $K/k$ pr\'{e}serve la mesure de la taille sans  conservation de la  finitude restreinte, donc on peut conjecturer que les extensions conservent la mesure de la taille de leurs corps interm\'{e}diaires.
La r\'esolution de ce probl\`eme nous am\`{e}ne \`{a} \'{e}tudier de pr\`{e}s la $lq$-finitude absolue.
\sk

Une extension $K/k$ est dite $lq$-finie si pour tout $n\in \N$, $k^{p^{-n}}\cap K/k$ est finie, et si de plus $di(K/k)$ est fini $K/k$ sera appel\'{e}e extension $q$-finie.  Il est trivialement \'evident que la finitude entraine la $q$-finitude entraine la $lq$-finitude, et lorsque $di(k)$ est fini, toute extension $K/k$ est $q$-finie. Pour des raisons de la non-contradiction, on donne des exemples qui  diff\'{e}rencient ces notions.
Toutefois, notre objectif pricipal consiste \`a d\'ecrire les extensions absolument $lq$-finies, \`a cet \`egard on montre que les propri\'{e}t\'{e}s suivantes sont \'{e}quivalentes :
\sk

\begin{itemize}{\it
\item[\rm{(1)}] Pour toute sous-extension $L/k$ de $K/k$, $K/L$ est $lq$-finie.
\item[\rm{(2)}] Toute sous-extension $L/k$ de $K/k$ satisfait $L/k(L^p)$ est finie.
\item[\rm{(3)}] Pour toute sous-extension $L/k$ de $K/k$, on a $L/rp(L/k)$ est finie,  o\`u $rp(L/k)$ est la cl\^oture relativement parfaite de $L/k$.
\item[\rm{(4)}] Toute suite d\'{e}croissante de sous-extensions de $K/k$ est stationnaire.}
\end{itemize}
\sk

Une extension qui v\'erifie l'une des conditions ci-dessus s'appelle extension absolument $lq$-finie. Il s'agit de la distribution absolue de la $lq$-finitude. On v\'{e}rifie aussit\^{o}t que la $q$-finitude entraine la $lq$-finitude absolue entraine la $lq$-finitude au sens large.
En particulier, la $lq$-finitude absolue est transitive, et toute extension qui est \`a la fois modulaire et  $lq$-finie est $q$-finie.
Plus g\'en\'erale\-m\-ent,  on montre que $K/k$ est $q$-finie si et seulement si  $K/k$ est absolument $lq$-finie et pour tout corps interm\'ediaire $L$ de $K/k$, la plus petite sous-extension $m$ de $L/k$ telle que $L/m$ est modulaire est non triviale ($m\not=L$).  Pour ne mentionner que cela, on montre que toute extension absolument $lq$-finie est compos\'{e}e d'\'{e}xtensions $w_0$-g\'{e}n\'{e}ratrices.
Enfin, nous approuvons cette notion par la construction  d'un exemple d'extension de taille infinie telle que pour tout corps  interm\'{e}diaire $L$ de $K/k$, on a $L/k$ est $q$-finie. En outre, $K/k$ ne respecte pas la r\'{e}partition de la finitude horizontale au niveau de leurs sous-extensions. D'autre part, cet exemple permet la distinction entre la $lq$-finitude simple et absolue.

\section{G\'{e}n\'{e}ralit\'{e}}

D'abord, nous commencerons par donner une  liste pr\'eliminaire  des notations le plus souvent utilis\'ees  tout le long de ce travail :
\sk

\begin{itemize}{\it
\item $k$ d\'esigne toujours un corps commutatif de caract\'eristique $p>0$, et $\Om$ une cl\^oture alg\'ebrique de $k$.
\item $k^{p^{-\infty}}$ indique la cl\^oture purement ins\'eparable de $\Om/k$.
\item Pour tout $a\in \Om$, pour tout $n\in \N^*$, on symbolise la racine du polyn\^ome $X^{p^n}-a$ dans $\Om$  par $a^{p^{-n}}$. En outre, on pose $k(a^{p^{-\infty}})=k(a^{p^{-1}},\dots, a^{p^{-n}},$ $\dots )=\di\bigcup_{n\in \N^*} k(a^{p^{-n}})$ et
$k^{p^{-n}}=\{a\in \Om\,|\ , a^{p^{n}}\in k\}$.
\item Pour toute famille $B=(a_i)_{i\in I}$  d'\'el\'ements de $\Om$, on note $k(B^{p^{-\infty}})= k({({a_i}^{p^{-\infty}})}_{i\in I})$.
\item Enfin, |.| sera employ\'e au lieu du terme cardinal.}
\end{itemize}
\sk

 Il est \`a signaler aussi que toutes les extensions qui interviennent dans ce papier sont des sous-extensions purement ins\'eparables de $\Om$, et il est commode de noter $[k,K]$ l'ensemble des corps interm\'{e}diaires d'une extension $K/k$.
\subsection{$r$-base, $r$-g\'en\'erateur}
\begin{df}Soit $K/k$ une extension. Une partie $G$ de $K$
est dite $r$-g\'en\'era\-t\-eur de $K/k$, si $K=k(G)$ ; et si de plus pour
tout $x\in G$, $x\not \in k(G\backslash x)$, $G$ sera appel\'{e}e
$r$-g\'en\'erateur minimal de $K/k$.
\end{df}
\begin{df}Etant donn\'ees  une extension $K/k$  de caract\'eristique $p>0$  et
une partie $B$ de $K$. On dit que $B$ est une  $r$-base de $K/k$, si $B$
est un $r$-g\'en\'erateur minimal de $K/k(K^p)$. Dans le m\^eme ordre d'id\'ees, on dit que $B$ est
$r$-libre sur $k$, si $B$ est une $r$-base de $k(B)/k$~; dans le cas
contraire $B$ est dite $r$-li\'ee sur $k$.
\end{df}

Voici quelques cas particuliers :
\sk

\begin{itemize}{\it
\item Toute $r$-base de $k/k^p$ s'appelle $p$-base de $k$.
\item Egalement, toute partie d'\'el\'ements de $k$, $r$-libre sur $k^p$ sera appel\'ee $p$-ind\'epend\-a\-n\-te (ou $p$-libre) sur $k^p$.}
\end{itemize}
\sk

Ici $B$ d\'esigne une partie d'un corps commutatif $k$ de caract\'eristique $p>0$. Comme cons\'equeces imm\'ediates on a :
\sk

\begin{itemize}{\it
\item[\rm{(1)}] $B$  est $p$-base de $k$ si et seulement si pour tout $n\in \Z$, $B^{p^n}$ l'est \'egalement  de $k^{p^n}$.
\item[\rm{(2)}] $B$ est $r$-libre sur $k^p$ si et seulement si pour tout $n\in \Z$, $B^{p^n}$ l'est auusi sur  $k^{p^{n+1}}$.
\item[\rm{(3)}] $B$  est $p$-base de $k$ si et seulement si $B$ est un $r$-g\'en\'erateur minimal de $k/k^p$.
\item[\rm{(4)}] $B$  est $p$-base de $k$ si et seulement si pour tout $n\in \N^*$, $k^{p^{-n}}=\otimes_k (\otimes_k $ $k($ $a^{p^{-n}}$ $))_{a\in B}$
et  pour tout $a\in B$, $a\not\in k^p$. En particulier, $B$  est $p$-base de $k$ si et seulement si $k^{p^{-\infty}}=\otimes_k (\otimes_k k(a^{p^{-\infty}}))_{a\in B}$ et pour tout $a\in B$, $a\not\in k^p$.}
\end{itemize}
\sk

Il est \`a noter que le produit tensoriel  est  utilis\'e conform\'ement  \`a la d\'efinition {5} (cf. \cite{N.B2}, III, p. 42). Il est vu comme limite inductive du produit tensoriel d'une famille finie de $k$-alg\`ebre.
Toutefois, la proposition ci-dessous permet de ramener  l'\'{e}tude des propri\'{e}t\'{e}s des syst\`{e}mes  $r$-libres des extensions  de haureur $\leq 1$,  ($K^p\subseteq k$) au cas fini. Plus pr\'{e}cis\'{e}ment, on a :
\begin{pro} {\label{pr1}} Soit $K/k$ une extension de caract\'{e}ristique $p>0$. Une partie $B$ de $K$ est $r$-libre sur $k(K^p)$ si et  seulement s'il en est de m\^{e}me pour toute sous-partie finie de $B$.
\end{pro}
\pre Imm\'{e}diat. \cqfd
\begin{pro} {\label{pr2}} Soit $K/k$ une extension de caract\'{e}ristique $p>0$. Toute partie finie $B$ de $K$ satisfait  $[k(K^p)(B) :k(K^p)]\leq p^{|B|}$, et il y'a \'{e}galit\'{e} si et seulement si $B$ est $r$-libre sur $k(K^p)$.
\end{pro}
\pre  Notons $B=\{x_1,\dots,x_n\}$, comme pour tout $i\in \{1,\dots,n\}$, on a $x_{i}^p \in k(K^p)\subseteq k(K^p)(x_1,\dots,x_{i-1})$, alors $[k(K^p)(x_1,\dots,x_{i}) : k(K^p)(x_1,\dots,x_{i-1})]$ $\leq p$, et il y'a  \'{e}galit\'{e} si et seulement si $x_i\not\in k(K^p)(x_1,\dots,x_{i-1})$. Compte tenu de la transitivit\'e de la finitude,
$[k(K^p)(x_1,\dots,x_n) :k(K^p)]=\di\prod_{i=1}^{n}[k(K^p)(x_1,\dots,x_{i}) :k(K^p)(x_1,\dots,x_{i-1})]\leq p^n$, et il y'a \'{e}galit\'{e} si et seulement si $B$ est $r$-libre sur $k(K^p)$. \cqfd
\begin{cor} {\label{cor1}} Soit $K/k$ une extension de caract\'{e}ristique $p>0$. Une partie  $B$ de $K$ est $r$-libre sur $k(K^p)$ si et seulement si pour toute sous-partie finie $B'$ de $B$, on a $[k(K^p)(B') :k(K^p)]=p^{|B'|}$.
\end{cor}

Comme application, le r\'{e}sultat ci-dessous montre que la $r$-ind\'{e}pendance est transitive dans le cas des extensions de hauteurs $1$. Autrement dit :
\begin{pro} {\label{pr3}} Etant donn\'ee  une extension $K/k$ de caract\'eristique $p>0$. Deux parties
$B_1$ et $B_2$ de $K$ sont respectivement  $r$-libres sur $k(K^p)$ et $k(B_1)(K^p)$
si et seulement si  $B_1\cup B_2$ l'est
sur $k(K^p)$. En particulier, si $B_1$ est une $r$-base de
$k(K^p)/K^p$, et $B_2$ est une $r$-base de $K/k(B_1)(K^p)$, alors $B_1\cup
B_2$ est  $p$-base de $K$.
\end{pro}
\pre La condition suffisante r\'esulte aussit\^ot de la d\'efinition des $r$-bases. Par ailleurs, d'apr\`{e}s la proposition {\ref{pr1}}, on se ram\`{e}ne   au cas o\`{u} $B_1$ et $B_2$ sont finies.  En vertu de la proposition {\ref{pr2}},  on a $[k(K^p)(B_1\cup B_2) :k(K^p)]= [k(K^p)(B_1\cup B_2) :k(K^p)(B_1)]. [k(K^p)(B_1) :k(K^p)]=p^{|B_2|}.p^{|B_1|}=p^{|B_2|+|B_1|}=p^{|B_1\cup B_2|}$, et par suite $B_1\cup B_2$ est $r$-libre sur $k(K^p)$.  \cqfd \sk

Comme cons\'{e}quences imm\'ediates :
\begin{cor} {\label{adcor}} Soient $k\subseteq L\subseteq K$ des extensions purement ins\'eparables et $B_1$, $B_2$ deux parties respectivement de $K$ et $L$. Si $B_1$ est une $r$-base de $K/L$ et $B_2$ une $r$-base de $L(K^p)/k(K^p)$, alors $B_1\cup B_2$ est une $r$-base de $K/k$.
\end{cor}
\begin{cor} {\label{cor2}}  Soient $K/k$ une extension de caract\'eristique $p>0$, $x$ un \'el\'em\-ent de $K$, et $B$
une partie $r$-libre sur $k(K^p)$. Pour que $B\cup \{x\}$ soit  $r$-libre sur $k(K^p)$ il faut et il suffit que $x\not\in k(K^p)(B)$.
\end{cor}
\pre imm\'{e}diat. \cqfd
\begin{thm} {\label{thm1}} [th\'eor\`eme de la $r$-base incompl\`ete] Etant donn\'{e}es  une extension $K/k$ de caract\'eristique $p>0$, et une partie $B$ de $K$,  $r$-libre sur $k(K^p)$. Pour tout
$r$-g\'en\'erateur $G$ de $K/k(K^p)$, il existe un sous-ensemble $G_1$ de $G$
tel que $B\cup G_1$ est une $r$-base de $K/k$.
\end{thm}
\pre  Le cas o\`u $k(K^p)(B)=K$ est trivialement \'{e}vident.  Si $k(K^p)(B)\not= K$, il existe $x\in G$ tel que
$x\not\in k(K^p)(B)$. En effet, si pour tout $x\in G$, $x\in
k(K^p)(B)$,  comme $G$ est un $r$-g\'en\'erateur de $K/k(K^p)$, on aura $k(K^p)(G)=K\subseteq
k(K^p)(B)$, absurde. D'apr\`es le lemme pr\'ec\'edent, $B\cup \{x\}$
est une partie $r$-libre sur $k(K^p)$. Posons ensuite $H=\{L\subset G$ tel
que $B\cup L$ est $r$-libre sur $k(K^p)\}$. IL est clair que $H$ est inductif, et donc d'apr\`es
le lemme de Zorn, $H$ admet un \'el\'ement maximal que l'on note $M$. Soit  $B_1=M\cup B$, n\'ecessairement $K=k(K^p)(B_1)$, si
$K\not= k(K^p)(B_1)$, il existe \'{e}galement un \'el\'ement $y$ de $G$ tel que $y\not \in
k(K^p)(B_1)$, et donc $B_1\cup \{y\}$ serait
 $r$-libre sur $k(K^p)$~;  c'est une contradiction avec
le fait que $M$ est maximal.\cqfd \sk

Voici quelques cons\'equences imm\'ediates :
\sk

\begin{itemize}{\it
\item[\rm{(1)}]  De tout $r$-g\'en\'erateur de $K/k(K^p)$ on
peut en extraire  une $r$-base de $K/k$.
\item[\rm{(2)}]  Toute partie $r$-libre sur $k(K^p)$ peut
\^{e}tre compl\'et\'{e}e en une $r$-base de $K/k$. En particulier, toute partie $p$-ind\'ependante sur $k^p$ peut \^{e}tre \'{e}tendue en une $p$-base de $k$.
\item[\rm{(3)}]  Toute extension $K/k$ admet une $r$-base. En outre, tout corps commutatif de caract\'eristique $p>0$ admet une $p$-base.}
\end{itemize}
\sk

Par ailleurs, toutes les $r$-bases d'une m\^{e}me extension ont m\^{e}me cardinal comme le pr\'ecise le  r\'{e}sultat suivant.
\begin{thm} {\label{thm2}} Soit $K/k$ une extension de caract\'eristique $p>0$. Si $B_1$
et $B_2$  sont deux $r$-bases de $K/k$, alors $|B_1|=|B_2|$.
\end{thm}

Pour la preuve de ce th\'eor\`eme on se s\'ervira des r\'esultats
suivants.
\begin{lem} {\label{lem1}} {\bf [Lemme d'\'echange]} Sous les conditions du th\'eor\`eme
pr\'ec\'edent, pour tout $x\in B_2$, il existe $x_1\in B_1$ tel
que $(B_1\backslash \{x_1\})\cup \{x\}$ est une $r$-base de $K/k$.
\end{lem}
\pre  Choisissons un \'{e}l\'{e}ment arbitraire $x$ de $B_2$, comme $B_2$ est une $r$-base de $K/k$, il en r\'esulte que
 $\{x\}$ est $r$-libre sur $k(K^p)$. Compte tenu du
th\'eor\`eme {\ref{thm1}}, il existe $B'_1\subset B_1$ tel que $B'_1\cup
\{x\}$ est une $r$-base de $K/k$.  D'o\`{u}, $p=[k(K^p)(B'_1)(\{x\}) :k(K^p)(B'_1)] =[K :k(K^p)(B'_1)]=[k(K^p)(B'_1)( B_1\backslash B'_1) :k(K^p)(B'_1)]$, et comme
 $B_1\backslash B'_1$ est $r$-libre sur  $k(K^p)(B'_1)$,   on en d\'{e}duit que $|B_1\backslash B'_1|=1$, c'est-\`{a}-dire $ B_1\backslash B'_1$ est r\'eduit \`a un singleton. \cqfd
\begin{pro} {\label{pr4}}Soit $K/k$ une extension de caract\'eristique $p>0$. Si $K/k$ admet au moins  une $r$-base finie, alors  toutes les $r$-bases de $K/k$ sont finies et ont m\^eme cardinal.
\end{pro}
\pre Imm\'ediat, il suffit d'appliquer la proposition {\ref{pr2}}. \cqfd \sk

\noindent{} {\bf Preuve du th\'eor\`eme {\ref{thm2}}.} D'apr\`es la prroposition {\ref{pr1}}, on se ram\`ene  au cas  o\`u  $|B_1|$ et $|B_2|$ sont infinies.
Comme $B_1$ est une $r$-base de $K/k$, pour tout $x\in B_2$,  il existe une partie finie  $D(x)$ de $B_1$ telle
que $x\in k(K^p)(D(x))$, et par suite $K=k(K^p)(B_2)\subseteq k(K^p)(\di\bigcup_{x\in B_2}(D(x)))$.  Il en r\'esulte que $\di\bigcup_{x\in
B_2}(D(x))=B_1$, et  en vertu du (\cite{N.B}, III, p. 49, cor 3), on obtient $|B_1|\leq
|B_2|.|\N|=|B_2|$. De la m\^eme fa\c{c}on on montre que
$|B_2|\leq |B_1|$~;  d'o\`u
$|B_1|=|B_2|$. \cqfd \sk

Comme cons\'equence, on a :
\begin{cor} {\label{cor3}}
Pour toute  partie $B_1$ de $K$, $r$-libre sur $k(K^p)$, et tout
$r$-g\'en\'erateur $G$ de $K/k$, on a $|B_1|\leq |G|$.
\end{cor}
\pre Imm\'ediat, puisque tout $r$-g\'en\'erateur peut se r\'{e}duire (respectivement toute famille $r$-libre peut se compl\'eter) en une $r$-base. \cqfd \sk

Dans le cas o\`u $K/k(K^p)$ est fini, compte tenu du th\'eor\`eme de la $r$-base incompl\`ete, un $r$-g\'en\'erateur $G$ de $K/k(K^p)$ est une $r$-base de $K/k$ si et seulement si $|G|=Log_p([K:k(K^p)])$. En particulier, si $B$ est une $r$-base de $K/k$ et $G$ un
$r$-g\'en\'erateur de $K/k(K^p)$ tels que $|B|=|G|<+\infty$, alors $G$ est une $r$-base de $K/k$.
\sk

Soit $K /k$ une extension purement ins\'{e}parable de caract\'{e}ristique $p>0$. On rappelle que $K$ est dit d'exposant fini sur $k$,
s'il existe $e\in \N$ tel que $K^{p^e}\subseteq k$, et le plus petit entier qui satisfait cette relation sera appel\'e exposant (ou hauteur) de $K/k$. Certes,
la proposition  suivante  permet de ramener l'\'{e}tude des propri\'{e}t\'{e}s  des $r$-g\'en\'erateurs minimals des extensions
d'exposant fini au cas des extensions  de hauteur $1$,  lesquelles sont  plus riches.
\begin{pro} {\label{pr5}} Soit $K/k$ une extension purement ins\'eparable d'exposant fini. Pour qu'une  partie de $K$ soit une $r$-base de $K/k$ il faut et il suffit que elle soit $r$-g\'en\'erateur minimal de $K/k$.
\end{pro}
\pre
Soit $G$ une $r$-base de $K/k$, donc $K=k(K^p)(G)=\dots= k(K^{p^e})(G)=k(G)$, et s'il existe $x\in G$ tel que $x\in
k(G\backslash \{x\})$, on aura $x\in k(K^p)(G\backslash \{x\})$, c'est une contradiction avec le fait que
$G$ est une $r$-base de $K/k$. Inversement, pour tout $r$-g\'en\'erateur minimal $G$  de $K/k$, on a
$K=k(G)=k(K^p)(G)$, et s'il existe $x\in G$ tel que $x\in K=
k(K^p)(G\backslash x)=\dots =k(K^{p^e}) (G\backslash
\{x\}=k(G\backslash \{x\})$, on aura une contradiction avec le fait que $G$ est
un $r$-g\'en\'erateur minimal de $K/k$. \cqfd
\begin{thm} {\label{thm3}} Soit $L/k$ une sous extension d'une extension purement ins\'epa\-r\-a\-ble d'exposant fini $K/k$. Pour toutes $r$-bases $B_L$ et $B_K$ respectivement de
$L/k$ et $K/k$, on a $|B_L|\leq |B_K|$.
\end{thm}
\pre On distingue deux cas :
\sk

1-ier cas. Si $K/k$ est d'exposant 1, c'est-\`{a}-dire $K^p\subseteq k$, donc $L^p\subseteq k$. D'apr\`{e}s le th\'eor\`eme {\ref{thm1}}, il existe $B_1\subseteq B_K$
tel que $B_L\cup B_1$ est une $r$-base de $K/k$, et par suite $|B_L|\leq |B_L\cup B_1|=|B_K|$.
\sk

2-i\`{e}me cas. Etant donn\'e  un entier naturel $e$ distinct de $0$ et $1$. Raisonnons par r\'ecurrence en supposant  que le th\'{e}or\`{e}me est v\'{e}rifi\'e  pour toute extension d'exposant $<e$, et soit $K/k$ une extension purement ins\'eprable d'exposant $e$.  Il est clair que  $k(K^p)\subseteq L(K^p)\subseteq K$, et donc il existe  $B_1\subseteq B_L$ et $B_2\subseteq B_K$ telles que $B_1$ et $B_2$ sont deux $r$-bases respectivement de $L(K^p)/k(K^p)$ et $K/L(K^p)$. D'apr\`es la transitivit\'e de la $r$-ind\'ependance, $B_1\cup B_2$ est une $r$-base de $K/k(K^p)$. Posons ensuite $k_1=k(B_1)$ et  $B'_L=B_L\setminus B_1$ ; on v\'erifie aussit\^ot que  $L\subseteq k_1(K^p)=k_1({B_2}^p)$, et $k_1(K^p)/k_1$ est d'exposant $<e$. Par application de la propri\'et\'e de r\'ecurrence et du corollaire {\ref{cor3}},
on obtient $|B'_L|\leq  |{B_2}^p|=|B_2|$. Comme $B_1\cap B'_L=\emptyset$ et $B_1\cap B_2=\emptyset$, alors $|B_1\cup B'_L|\leq  |B_1\cup B_2|$, et par suite $|B_L|\leq |B_K|$. \cqfd
\section{Degr\'{e} d'irrationalit\'{e}}

 Soit $K/k$ une extension purement ins\'{e}parable. D\'esormais, et sauf mention expresse du contraire, pour tout $n\in \N^{*}$, on note $k_n=k^{p^{-n}}\cap K$, on obtient  ainsi  $k\subseteq k_1\subseteq \dots \subseteq k_n\subseteq \dots \subseteq K$, et $k_n/k$ est d'exposant fini. Soit $B_n$ une $r$-base de $k_n/k$, d'apr\`{e}s le th\'{e}or\`{e}me {\ref{thm3}}, $|B_n|\leq |B_{n+1}|$.
 Ensuite, on pose  $di(K/k)= \di \sup_{n\in \N^*}(|B_n|)$, on rappelle que le $\sup$ est utilis\'e ici au sens du (\cite{N.B}, III, p. 25, proposition {2}).
\begin{df}
 L'invariant $di(K/k)$ d\'{e}fini ci-dessus s'appelle le degr\'{e} d'irratio\-n\-al\-i\-t\'{e} de $K/k$.
\end{df}

En particulier, et pour des raisons de diff\'erenciation,  le degr\'{e} d'irrationalit\'{e} de $k/k^p$ sera appel\'e degr\'e d'imperfection de $k$ et sera not\'e $di(k)$.
\begin{rem}
$di(K/k)$ permet de mesurer la taille de l'extension $K /k$, et $di(k)$ la longueur  de $k$.
\end{rem}

Toutefois, on v\'erifie aussit\^ot que :
\sk

\begin{itemize}{\it
\item $di(K/K)=0$.
\item Pour tout $n\in \Z$, $di(k)=di(k^{p^n})=di(k^{p^{-\infty}}/k)$.
\item Compte tenu du corollaire {\ref{adcor}}, pour toute sous-extension $L/k$ de $K/k$, on a $di(K/k(K^p))=di(K/L(K^p))+di(L(K^p)/k(K^p))$. Plus g\'en\'eralement, si $K/k$ est d'exposant $1$, on a $di(K/k)=di(K/L)+di(L/k)$.
\item En vertu de la proposition {\ref{pr2}}, pour toute extension purement ins\'eparable d'exposant finie $K/k$, on a $di(K/k)=di(K/k(K^p))$.}
\end{itemize}
\begin{thm} {\label{thm4}} Soient $k\subseteq L\subseteq K$ des extensions purement ins\'{e}parables, on a
$di(L/k)\leq di(K/k)$. En outre, $di(K/k)=\di\sup(di(L/k))_{L\in [k, K]}$.
\end{thm}
\pre D'apr\`{e}s le th\'eor\`eme {\ref{thm3}}, il suffit de remarquer  que pour tout $n\geq 1$,  on a $di(k^{p^{-n}}\cap L/k)\leq di(k_n/k)$, et donc $\di\sup(di(k^{p^{-n}}\cap L/k))_{n\geq 1}\leq \di\sup(di(k_n/k))_{n\geq 1}$ ; ou encore $di(L/k)\leq di(K /k)$. \cqfd \sk

Une cons\'{e}quence type est le r\'{e}sultat suivant :
\begin{thm} {\label{thm5}}
Pour toute extension purement ins\'{e}parable $K /k$, on a $di(K /k)$ $\leq di(k)$.
\end{thm}
\pre Il  suffit de remarquer qu'une partie $B$ de $k$ est une $p$-base de $k$ si et seulement si  $B^{p^{-n}}$ est  une $r$-base de $k^{p^{-n}}/k$ pour tout $n\geq 1$. Comme $k^{p^{-\infty}}=\di\bigcup_{n\geq 1} k^{p^{-n}}$, on a  pour tout $n\geq 1$, $k^ {p^{-n}}\cap K\subseteq k^{p^{-\infty}},$  et par suite $di(K/k)\leq di(k^{p^{-\infty}}/k)=di(k)$. \cqfd
\begin{pro} {\label{pr6}} Soit $(K_n/k)_{n\in \N}$ une famille croissante de sous-extensions purement ins\'{e}parables d'une extension $\Om/k$. On a : $$di(\di\bigcup_{n\in\N}(K_n)/k)= \di\sup_{n\in\N}(di(K_n/k)).$$
\end{pro}
\pre  Notons $K=\di\bigcup_{n\in\N}K_n$, et soit $j$ un entier naturel non nul. Il est imm\'ediat que $k_j=k^{p^{-j}}\cap K=\di\bigcup_{n\in \N} (k^{p^{-j}}\cap K_n)$. Dans la suite on distingue deux cas :
\sk

 1-ier cas : si $di(k_j/k)$ est fini, ou encore $k_j/k$ est finie. Comme pour tout $n\in N$, on a $k^{p^{-j}}\cap K_n\subseteq k^{p^{-j}}\cap K_{n+1}\subseteq k^{p^{-j}}\cap K$, alors la suite d'entiers $([k^{p^{-j}}\cap K_n:k])_{n\in \N}$ est croissante et born\'ee, donc stationnaire \`a partir d'un rang $n_0$ ; et par cons\'equent pour tout $n\geq n_0$, $k^{p^{-j}}\cap K_n=k^{p^{-j}}\cap K_{n+1}$. En outre, $di(k^{p^{-j}}\cap K/k)=di(k^{p^{-j}}\cap K_{n_0}/k)=\di\sup_{n\in N} (di(k^{p^{-j}}\cap K_n/k))$.
\sk

 2-i\`eme cas : si $di(k^{p^{-j}}\cap K/k)$ est infini, ou encore $\di\sup_{n\in N} (di(k^{p^{-j}}\cap K_n/k))$ n'est pas fini. Comme $k^{p^{-j}}\cap K=\di\bigcup_{n\in \N} (k^{p^{-j}}\cap K_n)$, donc si $B^{j}_{n}$ est une $r$-base de $k^{p^{-j}}\cap K_n/k$, alors $\di\bigcup_{n\in\N} B^{j}_{n}$ est un $r$-g\'{e}n\'{e}rateur de $k^{p^{-j}}\cap K/k$. En vertu du corollaire {\ref{cor3}},  $di(k^{p^{-j}}\cap K/k)\leq |\di\bigcup_{n\in\N} B^{j}_{n}|$, et d'apr\`es (\cite{N.B}, III,  p.49, corollaire {3}),  $|\di\bigcup_{n\in\N} B^{j}_{n}|\leq \di\sup_{n\in\N}(|B^{j}_{n}|)=\di\sup_{n\in\N}(di(k^{p^{-j}}\cap K_n/k))$.
 \sk

 Compte tenu de ces deux cas, on en d\'eduit que  $di(K/k)\leq \di\sup_{n\in\N}(di(K_n/k))$. Mais comme $K_n\subseteq K$ pour tout $n\geq 1$, d'apr\`{e}s le th\'{e}or\'{e}me {\ref{thm4}} on obtient $\di\sup_{n\in\N}(di(K_n/k))\leq di(K/k)$, et par suite $di(K/k)=\di\sup_{n\in\N}(di(K_n/k))$.\cqfd \sk

Le r\'{e}sultat  suivant qui est une cons\'{e}quence bien connue de la lin\'{e}arit\'{e} disjointe   intervient souvent dans le reste de ce papier.
\begin{pro} {\label{pr7}} Soient $K_1/k$ et $K_2/k$ deux sous-extensions d'une m\^eme extension $K/k$, $k$-lin\'{e}airement disjointes.  Pour touts corps interm\'{e}diaires  $L_1$ et $L_2$  respectivement de $K_1$ et $K_2$, on a $L_2(K_1)$ et $L_1(K_1)$ sont $k(L_1, L_2)$-lin\'{e}airement-disjointes. En particulier, $L_2(K_1)\cap L_1(K_2)=k(L_1,L_2)$.
\end{pro}

Une famille $(F_i/k)_{i\in J}$ d'extensions est dites $k$-lin\'{e}airement disjointes,  si pour toute partie  $G$ d'\'{e}l\'{e}ments finis de $J$, $(F_n/k)_{n\in G}$ sont $k$-lin\'{e}airement disjointes (cf. \cite{F-J}, p. 36). Il est trivialement \'{e}vident que $k((F_i)_{i\in J})=\di \prod_{i\in J} F_i\simeq \otimes_k (\otimes_k F_i)_{i\in J}$ si et seulement si $(F_i/k)_{i\in J}$ sont $k$-lin\'{e}airement disjointes. De plus, les propri\'{e}t\'{e}s de la lin\'{e}arit\'{e} disjointe du cas fini se prolonge naturellement \`{a} une famille quelconques d'extensions  $k$-lin\'{e}airement disjointes. En particulier,  pour tout $i\in J$, soit $L_i$ un sous-corps de $F_i$, si $(F_i/k)_{i\in J}$ sont $k$-lin\'{e}airement disjointes,  compte tenu de la transitivit\'{e} de la lin\'{e}arit\'{e} disjointe,  $(L_i/k)_{i\in J}$ (resp. $((\di \prod_{n\in J}L_n)F_i/k)_{i\in J}$) sont $k$-lin\'{e}airement (resp. $\di \prod_{n\in J}L_n$-lin\'{e}airement) disjointes.
\sk

Consid\'erons maintenant  deux sous-extensions $K_1/k$ et $K_2/k$ d'exposant fini d'une m\^eme extension purement ins\'eparable  $K/k$. On v\'erifie aussit\^ot que si $B_1$ et $B_2$ sont deux $r$-bases respectivement de $K_1/k$ et $K_2/k$, alors $B_1$ et $B_1\cup B_2$ sont  deux $r$-g\'en\'erateurs respectivement de $K_1(K_2)/K_2$ et $K_1(K_2)/k$. En outre, $di(K_1(K_2)/K_2)\leq di(K_1/k)$ et $di(K_1(K_2)/k)\leq di(K_1/k)+di(K_2/k)$. D'une fa\c{c}on  plus pr\'ecise, on a :
\begin{pro} {\label{pr8}}Sous les conditions ci-dessus, et si de plus $K_1/k$ et $K_2/k$ sont
$k$-lin\'{e}airement disjointes, on a :
\sk

\begin{itemize}{\it
\item[\rm{(i)}] $B_1\cup B_2$ est une $r$-base de $K_1(K_2)/k$.
\item[\rm{(ii)}] $B_1$ est une $r$-base de $K_1(K_2)/K_2$.}
\end{itemize}
\end{pro}
\pre Ici, on se contente de pr\'{e}senter uniquement la preuve du premier item, puisque  les deux assertions  utilisent les m\^{e}mes techniques de raisonnement.
Il est clair que $K_1(K_2)=k(B_1\cup B_2)$, il suffit donc de montrer que $B_1\cup B_2$ est  minimal.  Pour cela, on suppose par exemple l'existence d'un  \'el\'ement $x$ dans $B_1$ tel que $x\in k((B_1\setminus\{x\})\cup B_2)=K$.  Comme $K_1/k$ et $K_2/k$ sont $k$-lin\'{e}airement disjointes,  par transitivit\'{e}, on a $k(B_1)=K_1$ et $K_2(B_1\setminus \{x\})=K$ sont $k(B_1\setminus \{x\})$-lin\'{e}airement disjoints, et
donc $K_1=K\cap K_1=k(B_1\setminus \{x\})$, c'est une contradiction avec le fait que $B_1$ est une $r$-base de $K_1/k$. \cqfd \sk

Comme cons\'{e}quence imm\'{e}diate, on a
\begin{cor} {\label{cor4}}  Soient $K_1$ et $K_2$ deux corps interm\'{e}diaires d'une m\^{e}me extension purement ins\'{e}parable $\Om/k$. Alors :
\sk

\begin{itemize}{\it
\item[\rm{(i)}] $di(K_1(K_2)/k)\leq di(K_1/k)+di(K_2 /k)$, et il y'a \'egalit\'e si $K_1$ et $K_2$ sont $k$-lin\'{e}airement disjoints.
\item[\rm{(ii)}] $di(K_1(K_2)/K_2)\leq di(K_1/k)$, et il y'a \'egalit\'e si $K_1$ et $K_2$ sont $k$-lin\'{e}air\-e\-ment disjoints.}
\end{itemize}
\end{cor}
 \pre Il suffit de remarquer que $K_1(K_2)=\di\bigcup_{j\in \N} (k^{p^{-j}}\cap K_1)(k^{p^{-j}}\cap K_2)=\di\bigcup_{j\in \N} K_2(k^{p^{-j}}\cap K_1)$, et si $K_1$ et $K_2$ sont $k$-lin\'{e}airement disjoints, d'apr\`{e}s la transitivit\'{e} de la lin\'{e}arit\'{e} disjoint, $k^{p^{-j}}\cap K_1$ et $k^{p^{-j}}\cap K_2$ sont aussi $k$-lin\'{e}airement disjoints pour tout $j\geq 1$.  On se ram\`{e}ne ainsi  au cas o\`{u} $K_1/k$ et $K_2/k$  sont d'exposant fini auquel cas le r\'{e}sultat d\'{e}coule imm\'{e}diatement de la proposition pr\'ec\'edente. \cqfd \sk

 Comme cons\'equence imm\'ediate, on a :
 \begin{cor} {\label{cor5}} Pour toute sous-extension $L/k$ d'une extension purement ins\'e\-p\-a\-rable $K/k$, on a $di(L(K^p)$ $ /k(K^p))\leq di(L/k(L^p))$, et il y'a \'egalit\'e si  $k(K^p)$ et $L$ sont $k(L^p)$-lin\'eairement disjointes.
 \end{cor}
\pre Due au corollaire {\ref{cor4}}. \cqfd \sk

Le r\'esultat suivant am\'eliore naturelement les conditions du th\'eor\`eme {\ref{thm4}}
\begin{thm} {\label{thm6}} Pour toute famille  d'extensions purement ins\'eparables $k\subseteq L\subseteq L'\subseteq K$, on a $di(L/L')\leq di(K/k)$.
\end{thm}
\pre Il est clair que $K=\di\bigcup_{j\in \N} L(k_j)$, et d'apr\`es la proposition {\ref{pr6}}, et le th\'eor\`eme {\ref{thm4}}, on a $di(L'/L)\leq di (K/L)=\di\sup_{j\in \N}(di(L(k_j)/k))\leq \di\sup_{j\in \N}(di(k_j/k))=di(K/k)$. \cqfd \sk

Comme cons\'equence imm\'ediate, on a :
\begin{cor} {\label{cor6}} Pour toute extension purement ins\'eparable $K/k$, on a $di(K)\leq di(k)$.
\end{cor}
\pre Il suffit de remarque que $K\subseteq k^{p^{-\infty}}$, et $di(K)=di(K/K^p)\leq di(k^{p^{-\infty}}/k^p)=di(k)$.\cqfd
 \subsection{Extensions relativement parfaites}

Au cours de cette section, on reprend, en les am\'eliorant,  quelques notions et r\'esultats de {\cite{Che-Fli4}}, puisqu'ils sont utilis\'es fr\'equemment ici.
\sk

Un corps $k$ de caract\'eristique $p$ est dit parfait si $k^{p}=k$ ; dans le m\^{e}me ordre d'id\'{e}es,
on dit que $K/k$ est relativement parfaite  si $k(K^{p})=K$. On v\'erifie ais\'ement que :
\sk

\begin{itemize}{\it
\item La relation "\^etre relativement parfaite" est transitive, c'est-\`a-dire si $K/L$ et $L/k$ sont relativement parfaites, alors $K/k$ l'est aussi.
\item Si $K/k$ est relativement parfaite, il en est de m\^eme de $L(K)/k(L)$.
\item La propri\'et\'e "\^etre relativement parfaite" est stable par un produit quelconque portant sur $k$. Autrement dit, pour toute famille $(K_i/k)_{i\in I}$ d'ext\-e\-nsions relativement parfaites, on a alors $\displaystyle \di \prod_{i}^{}K_{i}/k$ est aussi relativement parfaite.}
\end{itemize}
\sk

Par suite, il existe une plus grande sous-extension relativement parfaite de $K/k$ appel\'ee cl\^oture relativement parfaite de $K/k$, et se note $rp(K/k)$.
On a les relations
d'associativit\'{e}-transitivit\'{e} suivantes.
\begin{pro} Soit $L$ un corps interm\'{e}diaire de
$K/k$. Alors
$$
rp(rp(K/L)/k)=rp(K/k) \quad \mbox{ et } \quad
rp(K/rp(L/k))=rp(K/k).
$$
\end{pro}
\pre Cf. {\cite{Che-Fli4}}, p. {50},  proposition {5.2}. \cqfd
\begin{cor}Pour tout $L\in
[k,K]$, on a $K/L \hbox{ finie} \Longrightarrow rp(K/k) \subset  L.$
\end{cor}

En particulier, si $K/k$ est relativement parfaite, on a $K/L \hbox{
$finie$ } \Longrightarrow   L=K.$
Sch\'{e}matiquement on a un $trou$
\sk

\[\begin{array}{rcl}
k\longrightarrow &&K;\\
&\uparrow&\\
&\hbox{$trou$ }&
\end{array}\]

\noindent et ce $trou$ caract\'{e}rise le fait que $K/k$ est
relativement parfaite. En effet,  supposons que
$K/k$ v\'{e}rifie le $trou$ et soit $B$ une $r$-base de $K/k$.
Supposons $B\neq \emptyset$;
soit $x\in B$ et $L=k(K^{p})(B\setminus\left \{x\right \})$; on a $K/L$
finie, donc $K=L$ ce qui est absurde.

\begin{pro} {\label{arpaa1}} Soit $K/k$ une extension purement ins\'eparable telle que $[K:k(K^p)]$ est fini. Alors on a :
\sk

\begin{itemize}{\it
\item[{\rm (i)}] $K$ est relativement parfaite sur une extension finie de $k$.
\item[{\rm (ii)}] La suite d\'{e}croissante  $(k(K^{p^{n}}))_{n \in
\N}$ est stationnaire sur $k(K^{p^{n_{0}}})=rp(K/$ $k)$.}
\end{itemize}
\end{pro}
\pre Cf. {\cite{Che-Fli4}}, p.  {51},
 lemme {2.1}. \cqfd \sk

Comme  cons\'equence de la proposition pr\'ec\'edente, on a :
\begin{pro} Soit $K/k$ une extension purement ins\'eparable telle que $[K:k(K^p)]$ est fini. Pour tout $L\in [k,K]$, on a $rp(K/L)=L(rp(K/k)).$
\end{pro}
\pre Cf. {\cite{Che-Fli4}}, p. {51},
 proposition {6.2}. \cqfd \sk

 En utilisant le lemme 1.16 qui se trouve dans (\cite{Mor-Vin}, p. 10), on peut affirmer que la condition de finitude de $[K:k(K^p]$ est n\'ec\'essaire, et par suite, le r\'esultat pr\'ec\'edent peut tomber en d\'efaut si $K/k(K^p)$ n'est pas finie.
Par ailleurs, on v\'erifie aussit\^ot que $k(K^p)=rp(K/k)(K^p)$, et donc pour qu'une partie $G$ de $K$ soit $r$-base de $K/k$ il faut et il suffit qu'elle en soit de m\^eme de $K/rp(K/k)$. De plus, comme $2$-i\`eme cons\'equence de la proposition {\ref{arpaa1}}, le r\'esultat suivant exprime une condition n\'ecessaire et suffisant pour que $K/rp(K/k)$ soit finie. Plus pr\'ecis\'ement, on a :
\begin{pro} {\label{arp1}} Soit $K/k$ une extension purement ins\'eparable, alors $K/rp($ $K/$ $k)$ est finie si est seulement il en est de m\^eme de $K/k(K^p)$.
\end{pro}
\pre R\'esulte de la proposition {\ref{arpaa1}}. \cqfd
\section{Extensions $q$-finies}
\begin{df}  Toute extension de degr\'e d'irrationalit\'{e}  fini s'appelle extension $q$-finie.
\end{df}

En d'autres sens, la $q$-finitude est synonyme de la finitude horizontale. Toutefois, la finitude se traduit par la finitude horizontale et verticale, il s'agit de la finitude au point de vue taille et hauteur. Autrement dit, $K/k$ est finie si et seulement si $K/k$ est $q$-finie d'exposant born\'e. Par ailleurs, on v\'{e}rifie que {\it
le degr\'e d'irrationalit\'e d'une extension $K/k$ vaut 1 si est seulement si l'ensemble de corps interm\'ediaires de $K/k$ est totalement ordonn\'e.}
Ensuite, on appelle extension $q$-simple toute extension qui satisfait l'affirmation pr\'ec\'edente.
\begin{rem} On rappelle que lorsque $di(k)$ est fini, et apr\`es avoir montr\'e dans \cite{Che-Fli2} que $K/k(K^p)$ est finie et $di(K)\leq di(k)$,  le degr\'e d'irrationalit\'e d'une extension purement ins\'eparable $K/k$ a \'et\'e d\'efini  par l'entier  $di(K/k)=di(k)-di(K)+di(K/k(K^p))$. En outre, toute extension est $q$-finie si $di(k)$ est fini.  Avec quelques  modifications l\'eg\`eres, on peut toujours prolonger cette d\'efinition au cas o\`u $di(k)$ est non born\'e.  Commen\c{c}ons par le choix d'une extension $K/k$ relativement parfaite et $q$-finie.
Etant donn\'ee une $p$-base $B$ de $k$, donc $k=k^p(B)$, et par suite $k(K^p)=K^p(B)$. Comme $K/k$ est relativement parfaite, alors $K=k(K^p)=K^p(B)$. D'apr\`es le th\'eor\`eme {\ref{thm1}}, il existe $B_1\subseteq B$ telle que $B_1$ est une $p$-base de $K$. Ainsi,  on aura $k^{p^{-\infty}}=k(B^{p^{-\infty}})=k({B_1}^{p^{-\infty}})\otimes_k  k((B\setminus B_1)^{p^{-\infty}})\simeq K^{p^{-\infty}}\simeq  K \otimes_k k({B_1}^{p^{-\infty}} )$. En particulier, d'apr\`es le corollaire {\ref{cor4}}, $di(K/k)=di(K \otimes_k k({B_1}^{p^{-\infty}} )/ k({B_1}^{p^{-\infty}} ))=di(k^{p^{-\infty}}/ k({B_1}^{p^{-\infty}} ))=di(k((B\setminus B_1)^{p^{-\infty}})/k)=| B\setminus B_1|$. Si on interpr\`ete (par abus de langage) $| B\setminus B_1|$ comme diff\'erence de degr\'e d'imperfection de $k$ et $K$ en \'ecrivant $| B\setminus B_1|=di(k)-di(K)$, on obtiendra  $di(K/k)=di(k)-di(K)$.    Dans le cas g\'en\'eral, supposons que $K/k$ est $q$-finie quelconque, donc $K/rp(K/k)$ est finie, d'o\`u $di(K)=di(rp(K/k))$; et par suite  $di(K/k)=di(rp(K/k)/k)+di(K/k(K^p))=di(k)-di(K)+di(K/k(K^p))$ (cf. proposition {\ref{pr11}} ci-dessous).
\sk

Il est \`a signaler en tenant compte de cette consid\'eration que tous les r\'esultats des articles \cite{Che-Fli1}, \cite{Che-Fli2}, \cite{Che-Fli3} se g\'en\'eralisent naturellement par translation  \`a une extension $q$-finie quelconque.
\end{rem}

Soient $L/k$ une sous-extension d'une extension $q$-finie $K/k$, pour tout $n\in \N$, on note toujours $k_n=k^{p^{-n}}\cap K$. On v\'{e}rifie aussit\^{o}t que :
\sk

\begin{itemize}{\it
\item[\rm{(i)}] La $q$-finitude est transitive, en particulier, pour tout $n\in \N$, $K/k(K^{p^n})$ et $k_n/k$ sont finies.
\item[\rm{(ii)}] Il existe $n_0\in \N$, pour tout $n\geq n_0$, $di(k_n/k)=di(K/k)$.}
\end{itemize}
\sk

Par ailleurs, voici quelques applications imm\'ediates des propositions {\ref{arpaa1}} et {\ref{arp1}}.
\begin{pro} {\label{pr9}} Soit $K/k$ une extension $q$-finie. La suite $(k(K^{p^n}))_{n\in\N}$ s'arr\-\^{e}\-te sur $rp(K/k)$ \`{a} partir d'un $n_0$. En particulier,  $K/rp(K/k)$ est finie.
\end{pro}

Comme cons\'{e}quence, on a :
\begin{cor} {\label{cor7}} La cl\^{o}ture relativement parfaite d'une extension $q$-finie $K/k$ n'est pas triviale. Plus pr\'{e}cis\'{e}ment,  $rp(K/k)/k$ est d'exposant non born\'{e} si  $K/k$ l'est.
\end{cor}
\pre Imm\'ediat. \cqfd
\begin{pro} {\label{pr10}} Pour toute  extension $q$-finie $K/k$, il existe $n\in \N$ tel que $K/k_n$ est relativement parfaite. En outre,  $k_n(rp(K/k))=K$.
\end{pro}
\begin{pro} {\label{pr11}} Le degr\'{e} d'irrationalit\'{e} d'une extension $q$-finie  $K/k$ v\'{e}rifie l'\'{e}galit\'{e} suivante : $di(K/k)$ $=di(rp(K/k)/k)+di(K/k(K^p))=di(K/rp(K/k))$ $+di(rp(K/k)/k)$.
\end{pro}
\pre Soient $G$ une $r$-base de $K/k$ et $K_r=rp(K/k)$, donc $k(G)/k$ admet un exposant  finie not\'e $m$ et, $K=K_r(G)$. En paticulier,  pour tout $n\geq m$,  $k(G)\subseteq k_n$. Compte tenu de la $r$-ind\'{e}pendance de $G$ sur $k(K^p)$ et vu que $k({k_n}^p)$ est un sous-ensemble de $k(K^p)$, on en d\'{e}duit que $G$ est $r$-libre sur $k({k_n}^p)$ pour tout $n\geq m$. Compl\'{e}tons $G$ en une $r$-base de $k_n/k$ par une partie $G_n$ de $k_n$. Dans ces conditions, pour $n$ suffisamment grand, on aura  $|G|+|G_n|=\di \sup_{j\geq m}(|G|+|G_j|)=di(K/k)=di(K_r(G)/k) \leq di(K_r/k)+di(k(G)/k)=di(K_r/k)+|G|$, et donc $|G_n|\leq di(K_r/k)$. Toutefois, comme $\di\bigcup_{n\geq m} k({k_n}^{p^m})=\di\bigcup_{n\geq m} k({G_n}^{p^m},G^{p^m})=\di\bigcup_{n\geq m} k({G_n}^{p^m})=k(K^{p^m})=K_r(K^{p^m})$, d'apr\`{e}s le th\'{e}or\`{e}me {\ref{thm4}}, pour $n$ suffisamment grand, on aura \'egalement $ di(K_r/k)\leq di(K_r(K^{p^{m}})/k)= di(k({k_n}^{p^m})/k)\leq |{G_n}^{p^m}|=|G_n|$. D'o\`u, $|G_n|=di(K_r/k)$ pour $n$ assez grand, et par suite $di(K/k)=di(K_r/k)+di(K/k(K^p))$. \cqfd
\sk

Comme cons\'equence imm\'ediate, on a :
\begin{cor} {\label{cor8}}Pour qu'une extension $q$-finie $K/k$ soit finie il faut et il suffit que $di(K/k)=di(K/k(K^p))$.
\end{cor}
\begin{thm} {\label{thm7}} Pour toutes extensions $q$-finies $k\subseteq L\subseteq K$, on a $di(K/k)\leq di(K/L)+di(L/k)$, avec l'\'{e}galit\'{e} si et seulement si $L/k(L^p)$ et $k(K^p)/k(L^p)$ sont $k(L^p)$-lin\'{e}airement disjointes.
\end{thm}
\pre Comme $K=\di\bigcup_{n\in \N} L^{p^{-n}}\cap K$ et $K/k$ est $q$-finie, d'apr\`{e}s le th\'{e}or\`{e}me {\ref{thm4}}, pour $n$ assez grand, on a $di(K/k)=di(L^{p^{-n}}\cap K/k)$ ; donc on est amen\'e au cas o\`{u} $K/L$ est finie, ou encore $rp(K/k)=rp(L/k)$.  Dans la suite,  on posera $L_r=K_r=rp(K/k)$. D'apr\`{e}s la proposition {\ref{pr11}} ci-dessus,  on aura $di(K/k)=di(K_r/k)+di(K/k(K^p))= di(L_r/k)+di(K/L(K^p))+di(L(K^p)/k(K^p))= di(L_r/k)+di(K/L)+di(L(K^p)/k(K^p))$. Compte tenu du corollaire {\ref{cor5}}, $di(L(K^p$ $)$ $/k(K^p))\leq di(L/k(L^p))$, donc $di(K/k)\leq di(L_r/k)+di(K/L) +di(L/k(L^p))=di(L/k)+di(K/L)$, et il y'a  \'{e}galit\'{e} si et seulement si $di(L/k(L^p))= di(L(K^p)/k($ $K^p))$, ou encore $[L : k(L^p)]=[L(K^p) :k(K^p)]$, c'est-\`a-dire $L/k(L^p)$ et $k(K^p)/k($ $L^p)$ sont $k(L^p)$-lin\'{e}airement disjointes.\cqfd
\begin{rem} La condition de la lin\'{e}arit\'{e} disjointe qui figure dans la proposition ci-dessus  se traduit en terme de $r$-ind\'{e}pendance par toute $r$-base de $L/k$ se compl\`{e}te en une $r$-base de $K/k$.
\end{rem}

Comme application imm\'{e}diate, on a :
\begin{cor}{\label{ccor1}} Toute sous-extension relativement parfaite $L/k$ d'une extension $q$-finies $K/k$ v\'erifie $di(K/k)=di(K/L)+di(L/k)$.
\end{cor}

D'une fa\c{c}on assez g\'en\'erale, on a :
\begin{pro} {\label{pr12}} Pour toute suite de sous-extensions relativement parfaites  $k=K_0 \subseteq K_1\subseteq \dots\subseteq K_n$ d'une extension $q$-finie $K/k$,  on a $di(K/k)=\di \sum_{i=0}^{n-1} di(K_{n+1}/K_n) +di(K/K_n)$.
\end{pro}
\pre R\'{e}sulte imm\'{e}diatement du corollaire pr\'{e}c\'{e}dent. \cqfd \sk

Dans la suite on va \'{e}tudier de plus pr\`{e}s  les propri\'{e}t\'{e}s des  exposants d'une extension $q$-finie.
\subsection{Exposants d'une extension $q$-finie}

Dans cette section nous distinguons  deux cas~:
\sk

{\bf Cas o\`u ${K/k}$ est purement ins\'{e}parable finie.}
Soit $x\in K$, posons $o( x/k ) = \inf\{$ $m \in \mathbf{N}|\; x^{p^m}\in k \}$
et $o_1(K/k) = \inf\{m\in \mathbf{N}|\; K^{p^{m}}\subset k\}$. Une
$r$-base $B=\{a_{1},a_{2},\dots, a_{n}\}$ de $K/k$ est dite
canoniquement ordonn\'{e}e si pour $j=1,2,\dots,n$, on a
$o(a_{j}/k(a_{1},a_{2},$ $\dots$ $ ,a_{j-1}))=
o_{1}(K/k(a_{1},a_{2},\dots ,a_{j-1})).$ Ainsi, l'entier
$o(a_j/k(a_1,\dots, $ $a_{j-1}))$ d\'{e}fini ci-dessus v\'{e}rifie
$o(a_j/k(a_1,\dots, a_{j-1}))=\inf\{m$ $\in \mathbf{N}|
\; di(k(K^{p^m})/k)\leq j-1\}$ (cf. {\cite{Che-Fli2}}, p.
{138}, lemme {1.3}). On en
d\'eduit aussit\^ot  le r\'{e}sultat de ({\cite{Pic}}, p.
{90}, satz {14}) qui confirme
l'ind\'{e}p\-e\-n\-d\-a\-n\-ce des entiers $o(a_i/k(a_1,$ $\dots ,
a_{i-1}) )$, $(1\leq i\leq n)$, vis-\`a-vis au choix des $r$-bases
canoniquement ordonn\'{e}es $\{a_1,\dots , $ $a_n\}$  de $K/k$. Par
suite, on pose  $o_i(K/k)=o(a_i/k(a_1,$ $\dots , a_{i-1}) )$ si
$1\leq i\leq n$, et $o_i(K/k)=0$ si $i>n$, o\`u $\{a_1,\dots , a_n\}$
est une $r$-base canoniquement ordonn\'{e}e de $K/k$. L'invariant $o_i(K/k)$
ci-dessus s'appelle le $i$-\`{e}me exposant de $K/k$. Voici les principales relations dont on aura besoin, et qui
font intervenir les exposants.
\begin{pro} {\label{pr13}}
Soient $K$ et $L$ deux corps interm\'{e}diaires d'une
extension $\Omega /k$, avec $K/k$ purement ins\'{e}parable finie.
Alors pour tout entier $j$, on a $o_{j}(K(L)$ $/k(L))$ $ \leq
o_{j}(K/k)$.
\end{pro}
\pre Cf. {\cite{Che-Fli}}, p. {373},
 proposition {5}. \cqfd
\begin{pro} {\label{pr14}}
Soit $K/k$ une extension purement ins\'{e}parable
finie. Pour toute sous-extension $L/L'$ de $K/k$, et  pour tout
$j\in\mathbf{N}$, on a  $o_j(L/L')\leq o_j(K/k)$.
\end{pro}
\pre cf. {\cite{Che-Fli}},
p. {374}, proposition {6}. \cqfd
\begin{pro} {\label{pr15}}
 Soient $\{\al_1 ,\dots,\al_n\}$
une $r$-base canoniquement ordonn\'{e}e de $K/k$, et $m_j$ le $j$-i\`{e}me exposant de $K/k$,  $1\leq j\leq n$. On a~:
\sk

\begin{itemize}{\it
\item[{\rm(1)}]  $k(K^{p^{m_j}})=k(\al_{1}^{p^{m_j}},\dots,
\al_{j-1}^{p^{m_j}})$.
\item[{\rm(2)}] Soit $\Lambda_j=\{(i_1,\dots, i_{j-1})$ tel que
$0\leq i_1<p^{m_1-m_j},\dots,0\leq i_{j-1}<p^{m_{j-1}-m_j}\}$,
alors $\{{(\al_1,\dots, \al_{j-1})}^{{p^{m_j}}\xi}$ tel que $\xi\in
\Lambda_j\}$ est une base de $k(K^{p^{m_j}})$ sur $k$.
\item[{\rm(3)}] Soient $n\in\mathbf{N}$ et $j$ le plus grand entier tel que
$m_j>n$. Alors $\{\al_{1}^{p^{n}},\dots, \al_{j}^{p^{n}}\}$ est une
$r$-base canoniquement ordonn\'{e}e de $k(K^{p^n})/k$, et sa liste
des exposants est $(m_1-n,\dots, m_j-n)$}.
\end{itemize}
\end{pro}
\pre  cf. {\cite{Che-Fli2}}, p.
{140}, proposition {5.3}. \cqfd
\begin{pro} {\label{pr16}}
Soient $K_1/k$ et $K_2/k$ deux
sous-extensions purement ins\'ep\-a\-rables de $K/k$. $K_1$ et $K_2$
sont $k$-lin\'{e}airement disjointes si et seulement si
$o_j(K_1(K_2)/K_2)=o_j(K_1/k)$ pour tout $ j\in\mathbf{N}$.
\end{pro}
\pre cf. {\cite{Che-Fli}}, p.  {374},
 proposition  {7}. \cqfd
\begin{pro} {\label{pr17}} (Algorithme de la compl\'etion des r-bases) Soient K/k une
extension purement ins\'eparable finie, $G$ un $r$-g\'en\'erateur de $K/k$, et $\{\al_1,\dots, \al_s\}$
un syst\`eme de $K$ tel que pour tout $j\in \{1,\dots, s\}$, $o(\al_j,k(\al_1,\dots,\al_{j-1}))=o_j(K/k)$.
Pour toute suite $\al_{s+1}, \al_{s+2}, \dots,$  d'\'el\'ements de $G$  v\'erifiant
$o(\al_m,k(\al_1,$ $\dots,\al_{m-1}))=\di\sup_{a\in G}(o(a,k(\al_1,\dots,\al_{m-1})))$,
la suite $(\al_i)_{i\in \N^*}$ s'arr\^ete sur un plus grand entier $n$  tel que $o(\al_n, k(\al_1,\dots,\al_{n-1}))>0$. En particulier, $\{\al_1,\dots,\al_n\}$ est une r-base canoniquement ordonn\'ee de $K/k$.
\end{pro}
\pre Cf. {\cite{Che-Fli2}}, p. {139},  proposition {1.3}.\cqfd
\sk

{\bf Cas o\`u $K/k$ est $q$-finie d'exposant non born\'e.}
Soit $K/k$ une extension $q$-finie. Rappelons que pour tout $n\in \N^*$, $k_n$ d\'esigne toujours $k^{p^{-n}}\cap K$. En vertu de la proposition {\ref{pr14}},  pour tout $j\in \N^*$, la suite des entiers naturels $(o_j(k_n/k))_{n\geq 1}$ est croissante, et donc
$(o_j(k_n/k))_{n\geq 1}$ converge vers $+\infty$, ou $(o_j(k_n/k))_{n\geq 1}$ est stationnaire \`{a} partir d'un certain rang. Lorsque $(o_j(k_n/k))_{n\geq 1}$ est born\'{e}e,  par construction,  pour tout $t\geq j$, $(o_t(k_n/k))_{n\geq 1}$ est aussi born\'{e}e (et donc stationnaire).
\begin{df} Soient $K/k$ une extension $q$-finie  et $j$ un entier naturel non nul. On appelle le $j$-i\`{e}me exposant de $K/k$ l'invariant $o_j(K/k)=\di\lim_{n\rightarrow +\infty} (o_j(k_n/k))$.
\end{df}
\begin{lem} {\label{lem2}} Soit $K/k$ une extension $q$-finie, alors  $o_s(K/k)$ est fini si et seul\-ement s'il existe un entier naturel $n$ tel que $di(k(K^{p^n})/k)<s$, et on a $o_s(K/k)=\inf\{m\in \N\,|\, di(k(K^{p^m})/k)<s\}$. En particulier, $o_s(K/k)$ est infini si et seulement si pour tout $m\in \N$, $di(k(K^{p^m})/k)\geq s$.
\end{lem}
\pre Pour simplifier l'\'{e}criture, on note $e_t=o_t(K/k)$ si $o_t(K/k)$ est fini. Compte tenu du  \cite{Che-Fli2}, p. {138}, lemme {1.3}, on v\'erifie aussit\^ot que $o_s(K/k)$ est infini si et seulement si pour tout $m\in \N$, $di(k(K^{p^m})/k)\geq s$, donc on se ram\`{e}ne au cas  o\`{u} $o_s(K/k)$ est fini. Par suite, il existe un entier $n_0$, pour tout $n\geq n_0$, $e_s=o_s(k_n/k)$. D'apr\`{e}s \cite{Che-Fli2} p. {138}, lemme {1.3},  $di(k({k_n}^{p^{e_s}})/k)<s$ et $di(k({k_n}^{p^{e_s-1}})/k)\geq s$. En vertu du th\'{e}or\`{e}me {\ref{thm4}},  $di(k(K^{p^{e_s}})/k)<s$ et  $di(k({K}^{p^{e_s-1}})/k)\geq s$. Autrement dit,  $o_s(K/k)=\inf\{m\in \N\,|\, di(k(K^{p^m})/k)<s\}$. \cqfd \sk

Le r\'{e}sultat ci-dessous permet de ramener l'\'{e}tude des propri\'{e}t\'{e}s des exposants des extensions $q$-finies aux extensions finies par le biais  des cl\^{o}tures relativement parfaites.
\begin{thm} {\label{thm8}} Soit $K_r/k$ la cl\^{o}ture relativement parfaite de degr\'{e} d'irratio\-n\-alit\'{e} $s$ d'une extension $q$-finie $K/k$, alors on a :
\sk

\begin{itemize}{\it
\item[\rm{(i)}] Pour tout $t\leq s$, $o_t(K/k)=+\infty$.
\item[\rm{(ii)}] Pour tout $t> s$, $o_t(K/k)=o_{t-s}(K/K_r)$.}
\end{itemize}
\sk

En outre, $o_t(K/k)$ est fini si et seulement si $t> s$.
\end{thm}
\pre  Pour tout $t\in {\N}^*$, notons $e_t=o_t(K/K_r)$. Comme pour tout entier $e$, on a $k(K^{p^e})=K_r(K^{pe})=\di\bigcup_{n\in \N}k({k_n}^{p^e})$, donc  $s=di(K_r/k)\leq di(k(K^{p^e})/k)=di(k({k_n}^{p^e})/k)$ pour $n$ suffisament grand. D'apr\`{e}s le lemme {\ref{lem2}}, on aura d'une part $o_t(K/k)=+\infty$ pour tout $t\leq s$, et
d'autrs part pour tout $n>s$, $di(K_r(K^{p^{e_{n-s}}})/k)=di(K_r/k)+di(K_r(K^{p^{e_{n-s}}})/K_r)<s+n-s=n$ et $di(K_r(K^{p^{e_{n-s}-1}})/k)=di(K_r/k)+di(K_r(K^{p^{e_{n-s}-1}})/K_r)\geq n$. Notamment, pour tout $n>s$, $o_n(K/k)=o_{n-s}(K/K_r)$. Toutefois, $o_n(K/k)$ est fini si et seulement si $n\leq s$.   \cqfd \sk

Voici une liste de cons\'equences imm\'ediates :
\begin{pro} {\label{pr18}} Soient $K$ et $L$ deux corps interm\'{e}diaires  d'une extension $q$-finie $M/k$. Pour tout  $j\in \N^* $, on a $o_j(L(K)/L)\leq o_j(K/k)$.
\end{pro}
\pre  Due au lemme {\ref{lem2}},  et \`{a} l'in\'{e}galit\'{e} suivante r\'{e}sultant du corllaire {\ref{cor4}}: $di(L(L^{p^{n}}, K^{p^n})/L)=di(L(K^{p^n})/L)\leq di (k(K^{p^n})/k)$ pour tout $n\in \N$. \cqfd
\begin{pro} {\label{pr19}} Etant donn\'ees des extensions $q$-finies $k\subseteq L\subseteq K$. Pour tout  $j\in \N^*$, on a $o_j(L/k)\leq o_j(K/k).$
\end{pro}
\pre Application imm\'ediate du lemme {\ref{lem2}},  et de l'in\'{e}galit\'{e} suivante r\'{e}sultant du th\'{e}or\`{e}me {\ref{thm4}} : $di(k(L^{p^{n}})/k) $ $\leq di (k(K^{p^n})/k)$ pour tout $n\in \N$. 
\sk

Par ailleurs la taille d'une extension relativement parfaite reste invariant, \`a une extension finie pr\`es comme l'indique le r\'esultat suivant.
\begin{pro} {\label{pr20}}  Etant donn\'ee  une sous-extension relativement parfaite  $K/k$ d'une extension $q$-finie $M/k$. Pour toute sous-extension finie $L/k$ de $M/k$, on a $di(L(K)/L)=di(K/k)$.
\end{pro}
\pre  En vertu du corollaire {\ref{cor4}}, il suffit de montrer que $di(L(K)/L)\geq di(K/k)$. Pour cela, on
pose d'abord $e=o_1(L/k)$ et $t=di(K/k)$. D'apr\`es le th\'eor\`eme {\ref{thm8}}, pour tout $s\in \{1,\dots, t\}$, $o_s(K/k)=+\infty$, donc pour $n$ assez grand,  on aura $o_t(k_n/k)>e+1$, en outre $L\subseteq k_n$ et $di(k_n/k)=di(K/k)$. Soit $\{\al_1,\dots, \al_t\}$ une $r$-base canoniquement ordonn\'ee de $k_n/k$, s'il existe $s\in \{1,\dots, t\}$ tel que $\al_s\in L({k_n}^p)(\al_1,\dots, \al_{s-1})$, d'apr\`es la proposition {\ref{pr14}}, on aura  $e<o_t(k_n/k)\leq o_s(k_n/k)=o(\al_s, k(\al_1,\dots,\al_{s-1}))\leq o_1(L({k_n}^p)(\al_1,\dots,\al_{s-1})/$ $k(\al_1,\dots,\al_{s-1}))\leq \di \sup(o_1(L/k), o_s(k_n/k)-1)=o_s(k_n/k)-1$, et donc $o_s(k_n/k)$ $\leq o_s(k_n/k)-1$, contradiction. D'o\`u, $\{\al_1,\dots, \al_t\}$ est une $r$-base de $L(k_n)/L$, et par suite, $t=di(K/k)=di(L(k_n)/L)\leq di(L(K)/L)$. \cqfd
\section{Extensions modulaires}

On rappelle qu'une extension
$K/k$ est dite modulaire si et seulement si pour tout
$n\in\mathbf{N}$, $K^{p^{n}}$ et $k$ sont $K^{p^{n}}\cap
k$-lin\'{e}airement disjointes. Cette notion a \'{e}t\'{e} d\'efinie
pour la premi\`{e}re fois par Swedleer dans {\cite{Swe}}, elle
caract\'{e}rise les extensions purement ins\'{e}parables, qui sont
produit tensoriel sur $k$ d'extensions simples sur $k$. Par ailleurs, toute $r$-base $B$ de $K/k$ telle que $K\simeq \otimes_k (\otimes_k k(a))_{a\in B}$ sera appel\'ee
$r$-base modulaire. En particulier, d'apr\`es le th\'eor\`eme de Swedleer,  si $K/k$ est d'exposant born\'e, il est \'equivalent de dire que :
\sk

\begin{itemize}{\it
\item[\rm{(i)}] $K/k$ admet une $r$-modulaire.
\item[\rm{(ii)}] $K/k$ est modulaire.}
\end{itemize}
\sk

Soient $m_j$ le $j$-i\`{e}me exposant d'une extension purement ins\'eparable finie $K/k$ et $\{\al_1 ,\dots,\al_n\}$
une $r$-base canoniquement ordonn\'{e}e de $K/k$,
donc d'apr\`es la proposition ${\ref{pr15}}$, pour tout $j\in \{2,\dots,  n\}$, il existe des constantes uniques
$C_{\ep}\in k$ telles que ${\al_{j}}^{p^{m_j}}=\di\sum_{\ep \in \Lambda_j}{C_{\ep}}{(\al_1,
\dots, \al_{j-1})}^{p^{m_j}\ep}$, o\`u $\Lambda_j=\{(i_1,\dots, i_{j-1})$ tel que
$0\leq i_1<p^{m_1-m_j},\dots,0\leq i_{j-1}<p^{m_{j-1}-m_j}\}$.
Ces relations s'appellent les \'{e}quations de d\'{e}finition de $K/k$.
\sk

Le crit\`ere ci-dessous permet de tester la modularit\'e d'une extension.
\begin{thm} {\label{thm9}}
{\bf [Crit\`ere de modularit\'{e}]} Sous les notations ci-dessus,
les propri\'{e}t\'{e}s suivantes sont \'{e}quivalentes~:
\sk

\begin{itemize}{\it
\item[\rm{(1)}] $K/k$ est modulaire.
\item[\rm{(2)}] Pour toute $r$-base  canoniquement ordonn\'{e}e $\{\al_1,\dots , \al_n\}$
de $K/k$, les $C_{\ep}\in k \cap K^{p^{m_j}}$ pour tout $j\in \{2,\dots,  n\}$.
\item[\rm{(3)}] Il existe une $r$-base canoniquement ordonn\'{e}e  $\{\al_1,\dots , \al_n\}$
de $K/k$ telle que les $C_{\ep}\in k \cap K^{p^{m_j}}$ pour tout $j\in \{2,\dots,  n\}$.}
\end{itemize}
\end{thm}
\pre cf. {\cite{Che-Fli2}}, p.  {142}, proposition {1.4}. \cqfd
\begin{exe} Soient $Q$ un corps parfait de
caract\'eristique $p>0$, $k=Q(X,Y,Z)$ le corps des fractions
rationnelles aux ind\'etermin\'ees $X,Y,Z$, et $K=k(\al_1,\al_2)$ avec $\al_1=X^{p^{-2}}$ et
$\al_2=X^{p^{-2}}Y^{p^{-1}}+Z^{p^{-1}}$. On v\'erifie
aussit\^ot que
\sk

\begin{itemize}{\it
\item[$\bullet$] $o_1(K/k)=2$ et $o_2(K/k)=1$,
\item[$\bullet$] $\al_2^{p}=Y\al_1^p+Z$.}
\end{itemize}
\sk

Si $K/k$ est
modulaire, d'apr\`es le crit\`ere du modularit\'e, on aura $Y\in k\cap K^p$ et $Z\in k\cap K^p$, et donc $Y^{p^{-1}}$ et $Z^{p^{-1}}\in K$. D'o\`u
$k(X^{p^{-2}},Y^{p^{-1}}, Z^{p^{-1}})$ $\subset $ $K$, et par suite,
$di(k(X^{p^{-2}},Y^{p^{-1}},Z^{p^{-1}})/k)=3<$ $di(K/k)=2$, contradiction.
\end{exe}

Le r\'{e}sultat suivant est cons\'{e}quence imm\'{e}diate de la  modularit\'{e}.
\begin{pro} {\label{apr2}} Soient $m,n\in \Z$ avec $n\geq m$. Si
$K/k$ est modulaire, alors
 $K^{p^{m}}/k^{p^{n}}$  est modulaire.
\end{pro}

La condition $n\geq m$ assure $k^{p^{n}}\subset K^{p^{m}}$.
\begin{pro} {\label{proa1}} Soit $K/k$ une extension purement ins\'{e}parable finie (respectivement, et modulaire), et soit $L/k$ une
sous-extension de $K/k$ (respectivement, et modulaire) avec
$di (L/k) = s$. Si $K^p\subseteq L$,  il existe une r-base canoniquement ordonn\'{e}e (respectivement, et modulaire)
$ (\al_1, \al_2, \dots, \al_n)$ de $K/k$, et $e_1, e_2,\dots, e_s\in \{1,p\}$ tels que $({\al_1}^{e_1}, {\al_2}^{e_2},\dots, {\al_s}^{e_s})$ soit une r-base canoniquement
ordonn\'{e}e (respectivement, et modulaire) de $L/k$. De plus, pour tout $j\in \{1,\dots,  s\}$, on a $o_j (K/k) = o_j (L/k)$, auquel cas $e_j = 1$, ou $o_j (K/k) = o_j (L/k) + 1$, auquel cas
$e_j = p$.
\end{pro}
\pre Cf. {\cite{Che-Fli2}}, p. {146},
 proposition {8.4}. \cqfd \sk

Le th\'{e}or\`{e}me suivant de  Waterhouse joue un r\^{o}le important
dans l'\'{e}tude des extensions  modulaires (cf. \cite{Wat} Th\'{e}or\`{e}me 1.1).
\begin{thm} Soient $(K_{j})_{j\in I}$ une famille
quelconque de sous-corps d'un corps $\Omega$, et $K$ un autre sous-corps de $\Omega $. Si pour tout $j\in I$, $K$ et $K_j$ sont  $K\cap K_j$-lin\'{e}airement
disjoints,
alors $K$ et  $\di\bigcap_{j} K_{j}$ sont  $K\cap (\di\bigcap_{j} K_{j})$-lin\'{e}airement disjoint.
\end{thm}

Comme cons\'equence, la modularit\'e est stable par une intersection quelconque portant soit au dessus ou en dessous d'un corps commutatif. Plus pr\'ecis\'ement, on a :
\begin{cor} {\label{apr4}} Sous les m\^emes hypoth\`eses du th\'eor\`eme ci-dessus, on a :
\sk

\begin{itemize}{\it
\item[\rm{(i)}] Si pour tout $j\in I$, $K_{j}/k$ est modulaire, il en est de m\^eme de $\di\bigcap_{j} K_{j}/k$.
\item[\rm{(ii)}] Si pour tout $j\in I$, $K/K_j$ est modulaire, il en est de m\^eme de $K/\di\bigcap_{j} K_{j}$.}
\end{itemize}
\end{cor}

D'apr\`{e}s le th\'{e}or\`{e}me de Waterhouse, il existe une plus petite sous-extension $m/k$ de $K/k$ (respectivement une plus petite extension $M/K$) telle que $K/m$ (respectivement $M/k$) est modulaire. D\'esormais, on note $m=mod(K/k)$ et $M=clm(K/k)$. Toutefois, l'extension $clm(K/k)$ sera appel\'ee cl\^oture modulaire de $K/k$.

Comme application imm\'{e}diate de la proposition {\ref{pr7}}, on a
\begin{pro} {\label{pr24}} Etant donn\'{e}es une $r$-base modulaire $B$ d'une extension modulaire $K/k$ et une famille $(e_a)_{a\in B}$ d'entiers tels que $0\leq e_a\leq o(a,k)$.  Soit $L=k(({a}^{p^{e_a}})_{a\in B})$, alors $L/k$ et $K/L$ sont modulaires, et $(B\setminus L)$, $(({a}^{p^{e_a}})_{a\in B}\setminus L)$ sont deux $r$-bases modulaires  respectivement de $K/L$ et $L/k$. En outre, pour tout $a\in B$, $o(a,L)=e_a$.
\end{pro}
\pre On se ram\`ene au cas  fini auquel le r\'esultat d\'ecoule de la transitivit\'e de la lin\'earit\'e disjointe.
En outre, pour toute partie $\{a_1,\dots, a_n\}$ d'\'el\'ement de $B$, $[L(a_1,\dots,a_n):L]=\di\prod_{i=1}^{n}p^{e_{a_i}}$. \cqfd \sk

Dans la suite, pour tout $a\in B$, on pose $n_a=o(a,k)$. Consid\'erons maintenant les sous-ensembles $B_1$ et $B_2$ de $B$ d\'efinis par $B_1=\{a\in B\,|\, n_a>j\}$, $B_2=B\setminus B_1=\{a\in B\,|\, n_a\leq j\}$, ($j$ \'etant un entier ne d\'epassant pas $o(K/k)$).
\sk

Comme Application de la proposition pr\'ec\'edente, on a :
\begin{thm} {\label{thm11}} Sous les conditions pr\'ecis\'ees ci-dessus, pour tout entier $j< o(K/k)$, on a $k_j=k((a^{n_a-j})_{a\in B_1}, B_2)$.
\end{thm}
\pre   Comme $K/k$ est r\'eunion inductive d'extentions modulaires  engendr\'ees par des parties finies de $B$, et compte tenu de la distributivit\'e de l'intersection par rapport \`a la r\'eunion,   on peut supposer sans perdre de g\'en\'eralit\'e que $K/k$ est finie d'exposant not\'e $e$. Soient $\{\al_1,\cdots,\al_n\}$ une $r$-base modulaire  et canoniquement ordonn\'{e}e de $K/k$, et $m_j$ le j-i\`{e}me exposant de $K/k$.  D\'{e}signons par   $s$ le plus grand entier tel que $m_s>j$, et $L=k(\al_{1}^{p^{m_1-j}},\dots,\al_{s}^{p^{m_s-j}}, \al_{s+1},\dots,$ $ \al_n)$. On v\'{e}rifie aussit\^{o}t que :
\sk

\begin{itemize}{\it
\item[\rm{(i)}] $L\subseteq k_j$,
 \item[\rm{(ii)}] $K\simeq k(\al_1)\otimes_k \dots\otimes_k k(\al_n)\simeq L(\al_1)\otimes_L \dots \otimes_L L(\al_s)$.}
\end{itemize}
\sk

Ainsi,  pour tout $x\in K$, il existe des constantes uniques  $C_{\ep}\in L$ telles que
$x=\di\sum_{\ep \in \Lambda}{C_{\ep}}{(\al_1,
\dots, \al_{s})}^{\ep}$, o\`u $\Lambda=\{(i_1,\dots, i_{s})$ tel que
$0\leq i_1<p^{m_1-j},\dots, 0\leq i_{s}<p^{m_{s}-j}\}$, et donc $x^{p^j}=\di\sum_{\ep \in \Lambda}{C_{\ep}}^{p^j}{({\al_1}^{p^j},
\dots, {\al_{s}}^{p^j})}^{\ep}$. Compte tenu de la proposition {\ref{pr15}}, $x^{p^j}\in k$ (c'est-\`{a}-dire $x\in k_j$) si et seulement si $x^{p^j}={C_{(0,\dots, 0)}}^{p^j}$, ou encore $x=C_{(0,\dots,0)}$. Par suite $x\in k_j$ si et seulement si $x\in L$,  autrement dit  $k_j=L$. \cqfd
\sk

Comme cons\'equence imm\'ediate, dans le cas de modulaire le r\'esultat suivant exprime une propri\'et\'e de stabilit\'e de la taille d'un certains corps interm\'ediaires. Plus pr\'ecis\'ement,
\begin{cor}{\label{acor1}} Pour toute extension modulaire $K/k$, pour tout $n\in \N$, on a $di(k_n/k)=di(k_1/k)$. En particulier, $di(K/k)=di(k_1/k)$.
\end{cor}

Le r\'{e}sultat suivant est bien connu (cf. \cite{Kim}).
\begin{pro} {\label{pr23}} Soit $K/k$ une extension purement
ins\'{e}parable et modulaire~; soit pour tout $n\in\mathbf{N}$,
$K_n=k(K^{p^n})$. Alors $k_n/k$, $K/k_n$,  $K_n/k$ et $K/K_n$ sont modulaires.
\end{pro}
\begin{pro} {\label{apr1}} Soient $K_1$ et $K_2$ deux sous-extensions de $K/k$ telles que $K\simeq K_1\otimes K_2$. Si pour tout $i\in \{1,2\}$,  $K_i/k$ est modulaire,  il en est de m\^eme de $K/k$.
\end{pro}
\pre Cf. {\cite{Che-Fli4}}, p. {55},
 lemme {3.4}. \cqfd \sk

Le r\'{e}sultat suivant \'{e}tend trivialement les hypoth\`{e}ses de la proposition {3.3}, \cite{Mor-Vin},  p. {94}, ainsi que le th\'{e}or\`{e}me {3.2},  \cite{Dev1}, p. {289}.
Il utilise plus particuli\`{e}rement les propri\'{e}t\'{e}s du syst\`{e}me canoniquement g\'{e}n\'{e}rateur (pour plus d'information cf. \cite{Mor-Vin}, d\'efinition {1.32}, p. 29).
\begin{pro}  {\label{apr3}}Soient $K_1$ et $K_2$ deux sous-extensions de $K/k$ telles que $K\simeq K_1\otimes K_2$. Si $K/K_1$ est modulaire, et $K_2/k$ est d'exposant born\'e,  il existe une partie $B$ de $K$ telle que $K\simeq K_1\otimes_k (\otimes_k (k(\al)_{\al\in B})$.
\end{pro}
\pre D'abord, comme $K\simeq K_1\otimes_k K_2$, alors pour tout $i\in \N$, pour toute $r$-base $C$ de $k({K_2}^{p^i})/k$,  $C$ est aussi une $r$-base de $K_1({K_2}^{p^i})/K_1$.
Choisissons ensuite une $r$-base $B$ de $K_2/k$, comme $K_2/k$ est d'exposant fini,  alors $B$ est un $r$-g\'{e}n\'{e}rateur minimal de $K_2/k$. Soit $B_1,\dots, B_n$ une partition de $B$ v\'{e}rifiant $B_1=\{x\in B| \, o(x,k)=o_1(K_2/k)=e_1\}$ et, pour tout $1<i\leq n$,  $B_i=\{x\in B| \, o(x,k(B_1,\dots, B_{i-1}))=o_1(K_2/k(B_1,\dots,B_{i-1}))=e_i\}$. Il est clair que $e_1>\dots>e_n$, et  en vertu de la lin\'earit\'e disjointe,   pour tout $i\in \{1\dots,n\}$,  pour tout $x\in B_i$,  on a \'egalement $o(x,K_1(B_1,\dots, B_{i-1}))=o_1(K/K_1(B_1,\dots,B_{i-1}))\}=e_i$.  En particulier, pour tout $i\in \{2,\dots, n\}$, $(\di\prod_{\al}{(G)}^{{\al}p^{e_i}})_G$, o\`{u} $G$ est une partie finie d'\'{e}l\'{e}ments de $ B_1\cup\dots\cup B_{i-1}$ et les $\al$ sont convenablement choisis,  est une base respectivement de $k({K_2}^{p^{e_i}})$ sur $k$  et $K_1({K_2}^{p^{e_i}})=K_1({K}^{p^{e_i}})$ sur $K_1$. Notons $M_i$ cette base, et soit $x\in B_i$, il existe des $c_{\al}\in k$ uniques tels que $x=\di\sum_{\al}c_{\al}y_{\al}$, ($y_{\al}\in M_i$), en outre les $c_{\al}$ sont aussi uniques dans $K_1$.
 D'autre part, en vertu de la modularit\'e, pour tout $i\in \{1,\dots, n\}$,  $K^{p^{e_i}}$ et $K_1$ sont $K_1\cap K^{p^{e_i}}$-lin\'eairement disjointes. Comme $K_1({K_2}^{p^{e_i}})=K_1({K}^{p^{e_i}})$ et $M_i\subseteq K^{p^{e_i}}$, alors $M_i$ est aussi une base de $K^{p^{e_i}}$ sur $K_1\cap K^{p^{e_i}}$. En tenant compte de l'unicit\'e de l'\'ecriture de $x$ dans la base $M_i$,  on en d\'eduit par identification que les $c_{\al}\in k\cap K^{p^{e_i}}$, et donc ${B_i}^{p^{e_i}}\subseteq  k\cap K^{p^{e_i}}({K_1}^{p^{e_i}}({B_1}^{p^{e_i}},\dots, {B_{i-1}}^{p^{e_i}}))$ pour tout $i\in \{1\dots, n\}$. Par application du (\cite{Mor-Vin}, proposition {3.3}, p. {94}), il existe une sous-extension modulaire $J/k$ d'exposant finie de $K/k$ telle que $K\simeq K_1\otimes_k J$. Ainsi, le r\'esultat d\'ecoule imm\'ediatement du th\'eor\`eme de Swedleer.  \cqfd \sk

Dans le cas fini, le r\'esultat suivant g\'en\'eralise la proposion ci-dessus.
\begin{pro} {\label{ajpr1}} Soient $K_1$ et $K_2$ deux corps interm\'{e}diaires ; $k$-lin\'{e}airement disjoints d'une extension purement ins\'eparable finie $L/k$ avec $di(L/K_1)=di(K_2$ $/k)=n$.
Soit $s$ le plus petit entier tel que $o_s(K_2/k)=o_n(K_2/k)$. Si $L/K_1$ est modulaire, il existe une $r$-base $\{\al_1,\dots, \al_n\}$  canoniquement ordonn\'{e}e de
$K_1(K_2)/K_1$ v\'{e}rifiant $K_1(K_2)\simeq K_1\otimes k(\al_1,\dots, \al_s)\otimes_k k(\al_{s+1})\otimes_k \dots \otimes_k k(\al_n)$.
\end{pro}
\pre Pour simplifier l'\'ecriture, pour tout $j\in\{1,\dots, n\}$, on note  $o_j (K_2/k)$ $ = e_j$, et $K=K_1(K_2)$ . Soit $\{\al_1,\dots,\al_n\}$ une $r$-base canoniquement ordonn\'{e}e de $K_2/k$.
Compte tenu de la proposition {\ref{pr16}}, $\{\al_1,\dots,\al_n\}$ est aussi une $r$-base canoniquement ordonn\'{e}e de $K/K_1$, et pour tout $j\in\{1,\dots, n\}$,   $o_j (K/K_1) = e_j$.
D'apr\`{e}s la proposition {\ref{pr15}}, pour tout $i\in \{s,\dots, n\}$, il existe des constantes uniques $C_{\ep}^{i}\in k$ telles que  ${\al_i}^{p^{e_n}}=\di\sum_{\ep\in \Lambda_{s-1}} C_{\ep}^{i}{(\al_1\dots\al_{s-1})}^{p^{\ep}}$ ($*$). En vertu de la proposition {\ref{pr16}}, pour tout $i\in \{s\dots,n\}$,
l'\'{e}quation de d\'{e}finition de $\al_i$ par rapport \`{a} $K_1(\al_1,\dots,\al_{s-1})$ est aussi d\'efinie par la relation ($*$) ci-dessus.
Comme $L/K_1$ est modulaire, en se servant du crit\`{e}re de modularit\'{e}, pour tout $(i,\ep)\in \{s,\dots, n\}\times \Lambda_{s-1}$,
on aura ${(C_{\ep}^{i})}^{p^{-e_n}}\in L$.    Posons ensuite, $F=k({(C_{\ep}^{i})}^{p^{-e_n}})$ o\`u $(i,\ep)$ parcourt l'ensemble $\{s,\dots, n\}\times \Lambda_{s-1}$,
et $H=K_1(F)(\al_1,\dots,\al_{s-1})$. Il est clair que $o_1(F/k)\leq e_n$, et $K\subseteq H\subseteq L$. De plus, d'apr\`{e}s le th\'eor\`eme {\ref{thm4}} et la proposition {\ref{pr14}},
$n=di(K/K_1)\leq di(H/K_1)\leq di(L/K_1)=n$, et pour tout $i\in \{s, \dots, n\}$, $e_n=o_i(K/K_1)\leq o_i(H/K_1) \leq e_n$. Il en r\'{e}sulte que  $di(H/K_1)=n$,
et pour tout $i\in \{s, \dots, n\}$, $e_n=o_i(H/K_1)$. Comme $e_{s-1}>e_s=e_n$, d'apr\`es l'algorithme de la compl\'{e}tion des $r$-bases, il existe des \'el\'ements $b_s,\dots, b_n \in  F$
tels que $\{\al_1\dots,\al_{s-1},b_s,\dots,b_n\}$ soit une $r$-base canoniquement ordonn\'{e}e de $H/K_1$. En particulier, on aura :
\sk

\begin{itemize}{\it
\item[$\bullet$] Pour tout $i\in \{1 ,\dots, s-1\}$, $e_i=o_i(H/K_1)=o_i(K_1(\al_1,\dots,\al_{s-1})/K_1)= o_i(k(\al_1,\dots,\al_{s-1})/k)$.
\item[$\bullet$] Pour tout $j\in \{s,\dots,n\}$,  $e_n=o_j(H/K_1)=o(b_j, K_1(\al_1\dots,
\al_{s-1}, b_s,\dots, $ $b_{j-1})) \leq o(b_j, k(b_s,\dots,b_{j-1})/k)\leq o_1(F/k)\leq e_n$,
et donc $e_n=o_j(H/K_1)$ $=o_j(k(b_s,\dots,b_n)/k)$.}
\end{itemize}
\sk

D'o\`u, $H=K\simeq K_1\otimes k(\al_1,\dots, \al_{s-1})\otimes_k k(b_s)\otimes_k\dots\otimes_k k(b_n)$. \cqfd
\section{Extensions \'equiexponentielles }
\begin{pro} {\label{pr26}} Soit $K/k$ une extension purement ins\'{e}parable d'exposant $e$. Les  assertions suivantes sont \'equivalentes :
\sk

\begin{itemize}{\it
\item[\rm{(1)}] Il existe une $r$-base $G$ de $K/k$ v\'{e}rifiant $K\simeq \otimes_k (k(a))_{a\in G}$, et pour tout $a\in G$, $o(a,k)=e$.
\item[\rm{(2)}] Toute  $r$-base $G$ de $K/k$ satisfait $K\simeq \otimes_k (k(a))_{a\in G}$, et $o_1(K/k)=e$.
\item[\rm{(3)}] Il existe une $r$-base $G$ de $K/k$ telle que pour tout $a\in G$, $o(a,k(G\setminus \{a\}))=o(a,k)=e$.
\item[\rm{(4)}] Pour toute $r$-base $G$ de $K/k$, pour tout $a\in G$, $o(a,k(G\setminus \{a\}))=o(a,k)=e$.}
\end{itemize}
\end{pro}
\pre D'apr\`es le th\'eor\`eme de la $r$-base incompl\`ete, on se ram\`{e}ne au cas o\`{u} $K/k$ est finie auquel cas $[K:k]=p^{en}$, o\`u $e=o_1(K/k)$ et $n=di(K/k)$, et en vertu de la proposition {\ref{pr16}}, le r\'esultat est imm\'ediat. \cqfd
\begin{df} Une extension qui v\'{e}rifie l'une des conditions de la proposition ci-dessus est dite \'equiexponentielle d'exposant $e$.
\end{df}

Il est clair  que toute extension \'equiexponentielle est modulaire. De plus, on v\'erifie aussit\^ot qu'il est \'{e}quivalent de dire que :
\sk

\begin{itemize}{\it
\item[\rm{(1)}] $K/k$ est \'equiexponentielle d'exposant $e$.
\item[\rm{(2)}] Il existe une $r$-base $G$ de $K/k$, pour toute partie finie $G_1$ de $G$, on a $k(G_1)/k$ est \'equiexponentielle d'exposant $e$.
\item[\rm{(3)}] Pour toute $r$-base $G$ de $K/k$, pour toute partie finie $G_1$ de $G$, on a $k(G_1)/k$ est \'equiexponentielle d'exposant $e$.}
\end{itemize}
\begin{pro} {\label{thm13}}
Pour toute extension relativement parfaite et modulaire $K/k$, pour tout entier $n$, $k_n /k$ est \'equiexponentielle d'exposant $n$.
\end{pro}
\pre D'apr\`es le th\'eor\`eme {\ref{thm11}}, il suffit de montrer que $k({k_n}^p)=k_{n-1}$. Compte tenu de la modularit\'e de $K/k$,  $K^{p^n}$ et $k$ sont $k\cap K^{p^n}$-lin\'{e}airement disjointes pour tout $n\geq 1$, et en vertu de la transitivit\'{e} de la lin\'{e}arit\'{e} disjointe, $k^{p^{n-1}}(K^{p^n})$ et $k$ sont $k^{p^{n-1}}(k\cap K^{p^{n}})$-lin\'{e}airement disjointes. Or $K/k$ est relativement parfaite, donc $k^{p^{n-1}}(K^{p^n})=K^{p^{n-1}}$, et par suite $k\cap K^{p^{n-1}}=k^{p^{n-1}}(k\cap K^{p^n})$, ou encore $k({k_n}^p)=k_{n-1}$. \cqfd \sk

Le r\'{e}sultat suivant, rapporte plus de pr\'{e}cision \`a la proposition {\ref{thm13}} dans le cas des extensions $q$-finies, notamment aux extensions finies.
\begin{pro} {\label{pr27}}
Soit $K/k$ une extension purement ins\'{e}parable de degr\'e d'irra\-t\-ionalit\'e $t$, relativement parfaite et modulaire (respectivement finie
et \'{e}quiexpo\-n\-entielle). Soient $n$ et $m$ deux entiers naturels tel que $n< m\hbox{ }(\hbox{respectivement, } n$ $< o_1(K/k))$. Les
propri\'{e}t\'{e}s suivantes sont v\'{e}rifi\'{e}es:
\sk

\begin{itemize}{\it
\item[{\rm (1)}] $di(k_{m}/k_{n})=t$.
\item[{\rm (2)}] $k_{m}/k_{n}$ est \'{e}quiexponentielle d'exposant
$m-n$;
\item[{\rm (3)}] $k_{n}^{p^{-(m-n)}}\cap K=k_{m}\;$ et
$\;k(k_{m}^{p^{m-n}})=k_{n}$.}
\end{itemize}
\sk

En particulier, pour tout $n\in \N$,  on a $[k_{n},k]=p^{nt}$.
\end{pro}
\pre cf. \cite{Che-Fli2}, p. {147}, proposition {9.4}. \cqfd \sk

Comme cons\'{e}quence imm\'{e}diate, on a :
\begin{cor} {\label{cor10}}
Si $K/k$ est une extension \'equiexponentielle d'exposant $e$, alors:
\sk

\begin{itemize}{\it
\item[\rm{(i)}] Pour tout $i\in \{1,\dots, e\}$, $k_i/k$ et $K/k_i$ sont \'equiexponentielles  d'exposant respectivement $i$ et $e-i$.
\item[\rm{(ii)}] Pour tout $i\in\{ 1,\dots, e\}$, $k(K^{p^{i}})/k$ et $K/k(K^{p^i})$ sont \'equiexponentielles  d'exposant respectivement $e-i$ et $i$.}
\end{itemize}
\end{cor}
\pre Imm\'ediat. \cqfd \sk

Le th\'eor\`eme ci-dessus reproduit dans un cadre plus \'etendu le corollaire {4.5} qui se trouve dans \cite{Dev1}, p. {292}, et pour plus d'information au sujet d'extraction des $r$-bases modulaires, on se r\'{e}f\`{e}re aux \cite{Dev1} et \cite{Dev-Mor2}.
\begin{thm} {\label{thm12}} Soient $k\subseteq L\subseteq K$ des extensions purement ins\'{e}parables telles que $K/k$ est \'equiexponentielle d'exposant $e$. Si $K/L$ est modulaire, il existe une $r$-base $G$ de $K/k$ telle que $\{a^{p^{o(a,L)}}| \, {a\in G}$ et $o(a,L)<e\}$ est une $r$-base modulaire de $L/k$.
\end{thm}
\pre Comme $K/L$ est modulaire d'exposant fini,  il existe une $r$-base $B_1$ de $K/L$ telle que $K\simeq \otimes_L (\otimes_L L(a))_{a\in B_1})$, (*). Pour des raisons d'\'{e}criture, pour tout $a\in B_1$, on pose $e_a=o(a,L)$ et $C=(a^{p^{e_a}})_{a\in B_1}$. Soit $B_2$ une partie de $L$ telle que $B_2$ est une $r$-base de $L(K^p)/k(K^p)$. Compte tenu de la transitivit\'e de $r$-ind\'ependance, $B_1\cup B_2$ est aussi une $r$-base de $K/k$. Dans la suite, notons $M=k(C,B_2)$. Il est clair que $M\subseteq L$, de plus, comme $K/k$ est \'equiexponentielle, on aura $K\simeq \otimes_k (\otimes_k k(a))_{a\in B_1\cup B_2}$. En vertu de la transitivit\'{e} de la lin\'{e}arit\'{e} disjointe, $K\simeq \otimes_M (\otimes_M M(a))_{a\in B_1}$, (**). En particulier,  d'apr\`{e}s les relations (*) et (**),  pour toute famille finie $\{a_1,\cdots, a_n\}$ d'\'el\'ements de $B_1$,  $L(a_1,\dots, a_n)\simeq L(a_1)\otimes_L\dots \otimes_L L(a_n)$ et $M(a_1,\dots, a_n)\simeq M(a_1)\otimes_M\dots\otimes_M M(a_n)$. Par application de la proposition {\ref{pr16}}, on a successivement $ [L(a_1,\dots, a_n) :L]=\di\prod_{i=1}^{n}p^{e_{a_i}}$ et $[M(a_1,\dots, a_n) :M]=\di\prod_{i=1}^{n}p^{e_{a_i}}$, ou encore $L$ et $K$ sont $M$-lin\'{e}airement disjointes. D'o\`{u} $L=L\cap K=M$. \cqfd
\section{$q$-finitude et modularit\'{e}}

Soit $K/k$ une extension $q$-finie d'exposant non born\'e. Dans tout ce qui suit, nous utilisons les notations suivantes : $k_j= k^{p^{-j}}\cap K$, $U_{s}^j(K/k)=j-o_s(k_j/k)$, et  $Ilqm(K/k)$ d\'{e}signe
le premier entier $i_0$ pour lequel la suite  $(U_{i_0}^{j}(K/k))_{j\in\mathbf{N}}$ est non born\'{e}e.
Le r\'{e}sultat ci-dessus est une application imm\'{e}diate de la proposition {\ref{pr14}}.
\begin{pro} {\label{pr25}} Etant donn\'ee une extension $q$-finie $K/k$.
 Pour tout entier $s$,  la suite $(U_{s}^j(K/k))_{j\in\mathbf{N}}$ est croissante.
\end{pro}
\pre
\smartqed Comme ${k_{n+1}}^p\subseteq k_n$, il est clair que $o_s(k_n/k)\leq o_s(k_{n+1}/k) \leq o_s(k_n/k)+1$,  et donc $n+1-
o_s(k_{n+1}/k)\geq n-o_s(K_n/k)$~; c'est-\`{a}-dire  la suite
$(U_{s}^j(K/k))_{j\in\mathbf{N}}$ est croissante. \cqfd \sk

En outre, on v\'erifie aussit\^ot que :
\sk

\begin{itemize}{\it
\item[\rm{(i)}] Pour tout $s\geq Ilqm(K/k)$, $\di\lim_{n \rightarrow
+\infty}(U_{s}^n (K/k))=+\infty$.
\item[\rm{(ii)}]  Pour tout $s<Ilqm(K/k)$,  la suite
$(U_{s}^j(K/k))_{j\in\mathbf{N}}$ est born\'{e}e ;}
\end{itemize}

\noindent et par suite, pour tout
$n\geq \di\sup_{j\in\mathbf{N}}(\di\sup (U_{s}^{j}(K/$ $k$ $)))_{s<Ilqm(K/k)}$, on
a $U_{s}^n(K/k)=U_{s}^{n+1}(K/k)$. Autrement dit,
$o_s(k_{n+1}/k)=o_s(k_n/k)+1$.
\sk

Dans toute la suite,  on pose $e(K/k)=\di\sup_{j\in\mathbf{N}}(\di\sup (U_{s}^{j}(K/k)))_{s<Ilqm(K/k)}$, et pour tout $(s,j)\in {\N}^*\times {\N}^{*}$, $e_{s}^{j}=o_s(k_j/k)$
\begin{thm} {\label{thm10}}  Soit $K/k$ une extension $q$-finie, avec $t=di(rp(K/k)/k)$. Les affirmations suivantes sont \'equivalentes ;
\sk

\begin{itemize}{\it
\item[\rm{(1)}] $K/k$ est modulaire sur une extension finie de $k$.
\item[\rm{(2)}]  Pour tout $s\in\{1,2,\dots,t\}$,  la suite  $(U_{s}^{j}(K/k))_{j\in\mathbf{N}}$  est born\'{e}e.
\item[\rm{(3)}] $Ilqm(K/k)=t+1$.}
\end{itemize}
\end{thm}
\pre
Il est clair que $(2) \Leftrightarrow (3)$. Par ailleurs, compte tenu de la proposition {\ref{pr10}}, il existe un entier $j_0$ tel que $K/k_{j_0}$ est relativement parfaite et $k_{j_0}(rp(K/k))=K$, et d'apr\`es la proposition {\ref{pr20}}, on aura $di(K/k_{j_0})= di(rp(K/$ $k)/k)=t$. Supposons ensuite que la condition $(1)$ est v\'{e}rifi\'{e}e. On distingue deux cas :
\sk

 Si $K/k$ est modulaire, en vertu de la proposition {\ref{pr27}}, pour tout $j\geq j_0$, on a $k_j/k_{j_0}$ est  \'{e}quiexponentielle d'exposant $j-j_0$ et $di(k_j/k_{j_0})=t$. D'o\`{u} pour tout $s\in \{1,\dots,t\}$, on a $U_{s}^{j}(K/k)=U_{s}^{j+1}(K/k)$.
 \sk

Si $K$ est modulaire sur une extension finie $L$ de $k$,  compte tenu de la finitude de  $L/k$, il existe un entier naturel $n$ tel que $L\subseteq k_n$. Par suite, $L^{p^{-j}}\cap K\subseteq k_{n+j}$, et donc $U_{s}^{n+j}(K/k)\leq n+U_{s}^{j}(K/L)$. D'o\`{u}, la suite $(U_{s}^{j}(K/k))_j$ est stationnaire pour tout $s\in\{1,\dots, t\}$.
\sk

 Inversement, si la  condition $(2)$ est v\'{e}rifi\'{e}e, il existe $m_0\geq \sup(e(K/k),j_0)$,  pour tout $j\geq m_0$, pour tout $s\in \{1,\dots,t\}$,  on a $o_s(k_{j+1}/k)=o_s(k_j/k)+1$ (et $di(k_j/k_{m_0})=t$). Par suite, $k_j/k_{j_0}$ est \'{e}quiexponentielle, donc modulaire. D'o\`{u} $K=\di\bigcup_{j>m_0} k_j$ est modulaire sur $k_{j_0}$. \cqfd
 \begin{thm} {\label{thm27}}
La plus petite sous-extension $M/k$ d'une extension $q$-finie $K/k$ telle que
$K/M$ est modulaire n'est pas triviale ($M\not= K$). Plus
pr\'{e}cis\'{e}ment, si $K/k$ est d'exposant non born\'e, il en est de m\^eme de $K/M$.
\end{thm}
\pre
Le cas o\`{u} $K/k$ n'est pas relativement parfaite (en particulier le cas fini)  est trivialement \'{e}vident, puisque $K/k(K^p)$ est modulaire. Ainsi, on est amen\'{e} \`{a} consid\'erer  que $K/k$ est relativement parfaite d'exposant non born\'{e}. On  emploiera ensuite un raisonnement par r\'{e}currence sur $di(K/k)=t$. Si $t = 1$, ou encore si $K/k$ est $q$-simple, il est imm\'{e}diat que  $K/k$ est modulaire.
Supposons maintenant que $t > 1$,  si  $Ilqm(K/k)=t+1$, en vertu du th\'eor\`eme {\ref{thm10}}, $M/k$ est finie, et donc $K/M$ est d'exposant non born\'e. Si  $Ilqm(K/k)\leq t$, pour tout $j>e(K/k)$, pour tout $s\in[1; i-1]$ o\`u $i=Ilqm(K/k)$, on a $e^{j+1}_{s}=e^{j}_{s}+ 1$. Comme $k^{p}_{j+1}\subseteq  k_j$, d'apr\`es la proposition {\ref{proa1}}, il
existe une r-base canoniquement ordonn\'ee $(\al_1,\dots, \al_n)$ de $k_{j+1}/k$, il existe $\ep_i,\dots,\ep_t\in \{1,p\}$  tels que $(\al_{1}^{p},\dots,\al_{i-1}^{p},\al_{i}^{\ep_i} \dots,\al_{t}^{\ep_t})$  est une r-base canoniquement ordonn\'{e}e de $k_j/k$. Dans la suite, pour tout $j>e(K/k)$, notons $K_j = k(k_{j}^{p^{e^{j}_{i}}})$. D'une part,  $K_j=k(\al_{1}^{p^{e^{j}_{i}+1}},\dots,\al_{i-1}^{p^{e^{j}_{i}+1}})$ et $K_{j+1}= k(\al_{1}^{p^{e^{j+1}_{i}}},\dots,\al_{i-1}^{p^{e^{j+1}_{i}}})$. D'autre part, on a $e^{j+1}_{i}=e^{j}_{i}+\ep $, avec $\ep=0$ ou $1$, cela conduit \`a $K_j\subseteq K_{j+1}$. Toutefois,  par d\'{e}finition de $Ilqm(K/k)$, on a $1 + e^{j}_{i }> e^{j+1}_{i}$ (c'est-\`a-dire $e^{j}_{i} = e^{j+1}_{i}$) pour une infinit\'{e} de valeurs de $j$. Pour ces valeurs, on a $di(K_{j+1}/k) = i -1$, sinon d'apr\`es le lemme {\ref{lem2}},  $e^{j+1}_{i} = e^{j+1}_{i-1} = 1 + e^{j}_{i-1} = e^{j}_{i}$ , et donc $e^{j}_{i} > e^{j}_{i-1}$, ce qui contredit la d\'efinition des exposants.
Comme $(di(K_j/k))_{j>e(K/k)}$ est une suite croissante d'entiers born\'{e}e par $di(K/k)$, donc elle stationne sur $Ilqm(K/k)-1$. De plus, $K_j \not= K_{j+1}$, en effet si $K_j =K_{j+1} =k( K^{p}_{j+1})$,  comme $K_{j+1}/k$ est d'exposant born\'e, on aura $K_{j+1} = k$, ce qui est absurde.
Posons ensuite $H =\di\bigcup_{j>e(K/k)} K_j$. On v\'erifie aussit\^ot que $H/k$ est d'exposant non born\'{e} et $di(H/k) = i -1$, de plus $H/k$ est relativement parfaite car $k(K^{p}_{j+1})= K_j$ pour une infinit\'{e} de $j$. Par ailleurs, d'apr\`{e}s les corollaires {\ref{cor8}} et {\ref{ccor1}}, $di(K/H)<t$ et $K/H$ est d'exposant non born\'e.
Compte tenu de l'hypoth\`{e}se de r\'{e}currence appliqu\'ee \`a $K/H$, on aura $K$ est modulaire sur une extension $M'$ de $H$ avec
$K/M'$ est d'exposant non born\'e; comme $M \subseteq M'$, alors $K/M$ est aussi d'exposant non born\'e. \cqfd \sk

Une version \'equivalente de ce  r\'{e}sultat  se trouve dans \cite{Che-Fli2}. Toutefois, le th\'eor\`eme ci-dessus peut tomber en d\'efaut lorsque l'hypoth\`ese de la $q$-finitude n'est pas v\'erifi\'ee comme le montre le contre-exemple ci-dessus
\begin{exe}
Soient $Q$ un corps parfait de caract\'{e}ristique $p>0$, et $(X,(Y_i)_{i\in\N^*}, $ $(Z_i)_{i\in \N^*}, (S_i)_{i\in\N^*})$ une famille alg\'{e}briquement ind\'{e}pendante sur $Q$. Soit $k=Q(X,(Y_i)_{i\in\N^*}, (Z_i)_{i\in \N^*}, (S_i)_{i\in\N^*})$ le corps des fractions rationnelles aux ind\'{e}\-t\-ermin\'{e}es $(X,(Y_i)_{i\in\N^*}, (Z_i)_{i\in \N^*}, (S_i)_{i\in\N^*})$.
Posons ensuite :
\sk

\begin{itemize}{\it
\item[] $K_1=\di\bigcup_{n\geq 1} k(\te_{1,n})$, avec $\te_{1,1}={X}^{p^{-1}}$ et $\te_{1,n}={\te_{1,n-1}}^{p^{-1}}$ pour tout entier $n>1$.
\item[] $K_2=\di\bigcup_{n\geq 1} K_1(\te_{2,n})$, o\`u $\te_{2,1}={Z_1}^{p^{-1}}\te_{1,2}+{Z_2}^{p^{-1}}$, et  pour tout $n>1$, $\te_{2,n}={Z_1}^{p^{-1}}\te_{1,2n}+{\te_{2,n-1}}^{p^{-1}}$.
\item[]  Par r\'{e}currence, on pose $K_j=\di\bigcup_{n\geq 1} K_{j-1}(\te_{j,n})$, o\`u $\te_{j,1}={Z_{j-1}}^{p^{-1}}\te_{j-1,2}$ $+{Z_j}^{p^{-1}}$, et pour tout $n>1$, $\te_{j,n}={Z_{j-1}}^{p^{-1}}\te_{j-1,2n}$ $+{\te_{j,n-1}}^{p^{-1}}$.}
\end{itemize}
\sk

Enfin, on note  $K=\di\bigcup_{j\in\N^*} K_j$, et  par conventient on pose $K_0=k$, et pout tout $i\in \N$, $\te_{i,0}=0$.
Comme pour tout $j\in \N^*$, $K_j\subseteq K_{j+1}$, alors $K$ est un corps commutatif.
\end{exe}
\begin{thm} {\label{thm14}} Sous les conditions ci-dessus, la plus petite sous-extension $m$ telle que $K/m$ est modulaire est triviale, c'est-\`a-dire $mod(K/k)=K$
\end{thm}

Pour la preuve de ce th\'{e}or\`{e}me, on se servira des r\'{e}sultats  suivants :
\begin{lem} {\label{lem3}} Sous les m\^{e}mes conditions ci-dessus, pour tout $(j,n)\in \N\times \N^*$, $K_j(\te_{j+1,n})=K_j({\te_{j+1,{n+1}}}^{p})$ et $\te_{j+1,1}\not\in K_j$. En particulier,  $o(\te_{j+1,n},K_j)=n$.
\end{lem}
\pre Il est trivialement \'evident que $K_j(\te_{j+1,n})=K_j({\te_{j+1,{n+1}}}^{p})$ pour tout $(j,n)\in \N\times \N^*$.  Pour le reste, il suffit de remarquer que $K_j\subseteq k({X}^{p^{-\infty}},{Z_1}^{p^{-\infty}},$ $\dots, {Z_j}^{p^{-\infty}})$ et $k(\te_{j+1,n},{X}^{p^{-\infty}},{Z_1}^{p^{-\infty}},$ $\dots, {Z_j}^{p^{-\infty}})= k({Z_{j+1}}^{p^{-n}}, {X}^{p^{-\infty}},{Z_1}^{p^{-\infty}},$ $\dots, {Z_j}^{p^{-\infty}})$,  et donc, pour tout $n\in \N^*$, $n=o(\te_{j+1,n},k({X}^{p^{-\infty}},{Z_1}^{p^{-\infty}},$ $\dots, {Z_j}^{p^{-\infty}}$ $))\leq o(\te_{j+1,n},K_j)\leq n $. \cqfd \sk

Comme cons\'equence imm\'ediate, pour tout $j\in \N^*$, $K_j/K_{j-1}$  est $q$-simple d'exposant non born\'e. En particulier, $di(K_j/k)=j$.
\begin{lem} {\label{lem4}} Pour tout $i\in \N^*$, la famille $(Z_i,(S_j)_{j\in \N^*})$ est $r$-libre sur $K^p$.
\end{lem}
\pre
Puisque pour tout $i\geq 1$, ${S_i}^{p^{-1}}\not\in k({X}^{p^{-\infty}}, ({Z_{j}}^{p^{-\infty}})_{j\geq 1})({S_1}^{p^{-1}},\dots, $ ${S_{i-1}}^{p^{-1}})=K({Z_1}^{p^{-\infty}})({S_1}^{p^{-1}},\dots, {S_{i-1}}^{p^{-1}})$, il suffit de montrer que ${Z_i}\not\in K^p$; ou encore ${Z_i}^{p^{-1}}\not\in K$. Par construction, pour tout $j\in \{1,\dots, n\}$, on a $\te_{j,1}={Z_{j-1}}^{p^{-1}}\te_{j-1,2}$ $+{Z_j}^{p^{-1}}$ avec $K_n$ contient $\te_{j,1}$ et $\te_{j-1,2}$, et donc s'il existe $n>i$ tel que ${Z_i}^{p^{-1}}\in K_n$, par it\'eration, on aura  ${Z_{i-1}}^{p^{-1}},\dots, {Z_1}^{p^{-1}}\in K_n$ et ${Z_{i+1}}^{p^{-1}},\dots, {Z_n}^{p^{-1}}\in K_n$. Par suite, d'apr\`es le th\'eor\`eme {\ref{thm4}}, $di(k(X^{p^{-1}}, {Z_1}^{p^{-1}},$ $\dots, {Z_n}^{p^{-1}})/k)\leq di(K_n/k)$, ou encore $n+1\leq n$, absurde.
D'o\`u pour tout $n\in \N^*$, ${Z_i}^{p^{-1}}\not\in K_n$, et comme $K$ est r\'eunion  de la famille croissante d'extensions $(K_n)_{n\in \N^*}$, alors  ${Z_i}^{p^{-1}}\not\in K$. \cqfd \sk

{\bf Preuve du th\'{e}or\`{e}me {\ref{thm14}}.} Posons $m=mod(K/k)$. En utilisant un raisonnement par r\'ecurrence on va montrer que  $K_i\subseteq  m$  pour tout $i\in {\N}$, et par suite obtenir $K=m$. Il est imm\'ediat que  $K_0=k\subseteq m$, donc le r\'esultat est v\'erifi\'e pour le rang $0$. Soit $i\in \N^*$, supposons par application de  l'hypoth\`ese de r\'ecurrrence que $K_i\subseteq m$. S'il existe un entier naturel $s$ tel que $\te_{{i+1},s}\not \in m$,  d\'esignons par $n$ le plus grant entier tel que $\te_{{i+1},n}\in m$. D'o\`{u} pour tout $t\in \{0, \dots, n\}$, $\te_{{i+1},t}\in m$ et $\te_{{i+1},n+1}\not\in m$, en outre $\te_{{i+1},2n}^{p^n}\in m$, et ${\te_{{i+1},2(n+1)}}^{p^{n+1}}\not\in m$. Il en r\'{e}sulte que le syst\`{e}me $({\te_{{i+1},2(n+1)}}^{p^{n+1}},1)$ est libre sur $m$, en particulier, il en est de m\^{e}me sur $m\cap K^{p^{n+1}}$. Compl\'{e}tons ce syst\`{e}me en une base $B$ de $K^{p^{n+1}}$ sur $m\cap K^{p^{n+1}}$. Comme $K^{p^{n+1}}$ et $m$ sont $m\cap K^{p^{n+1}}$-lin\'{e}airement disjointes ($K/m$ est modulaire), $B$ est aussi une base de $m(K^{p^{n+1}})$ sur $m$. Or, par construction,   $\te_{{i+2},n+1}={Z_{i+1}}^{p^{-1}}\te_{{i+1},2(n+1)}+{\te_{{i+2},n}}^{p^{-1}}$, ou encore ${\te_{{i+2},n+1}}^{p^{n+1}}={Z_{i+1}}^{p^n}{\te_{{i+1},2(n+1)}}^{p^{n+1}}+{\te_{{i+2},n}}^{p^n}$, avec ${\te_{{i+2},n}}^{p^n}={Z_{i+1}}^{p^{n-1}}{\te_{i+1,2n}}^{p^n}+\cdots + Z_{i+1}{\te_{i+1,2}}^{p}+Z_{i+2}\in m$.  Par identification, ${Z_{i+1}}^{p^n}\in m\cap K^{p^{n+1}}\subseteq K^{p^{n+1}}$, et donc ${Z_{i+1}}^{p^{-1}}\in K$, absurde. D'o\`{u} pour tout $n\in \N^*$, $\te_{{i+1},n}\in m$, ou encore $K_{i+1}\subseteq m$.  D'o\`{u} $m=K$. \cqfd
\section{Extensions $lq$-finies}

Dans ce qui suit nous proposons un cadre r\'enov\'e qui g\'en\'eralise naturellement le cas $q$-fini.
\begin{df} Soit  $K/k$ une extension purement ins\'{e}parable de caract\'{e}ristique $p>0$. On dit que $K/k$ est $lq$-finie, si pour tout $n\in \N^*$,
$k^{p^{-n}}\cap K/k$ est finie.
\end{df}

Il est trivialement \'evident que :
\sk

\begin{itemize}{\it
\item[$\bullet$] La $q$-finitude entraine la $lq$-finitude.
\item[$\bullet$] La  $lq$-finitude est stable par inclusion. Plus pr\'ecis\'ement, toute sous-extension d'une extension $lq$-finie est $lq$-finie.}
\end{itemize}
\sk

Voici \'egalement un exemple non trivial d'extension $lq$-finie de degr\'e d'irratio\-n\-alit\'e infini.
\subsection{Exemple non trivial d'extension $lq$-finie}
\begin{exe} {\label{exe1}}
Soient $Q$ un corps parfait de caract\'{e}ristique $p>0$, et $(X, (Y_i,$ $Z_i)_{i\in \N^*})$ une famille alg\'{e}briquement ind\'{e}pendante sur $Q$,
et $k=Q(X, (Y_i,Z_i)_{i\in \N^{*}})$ le corps des fractions rationnelles aux ind\'{e}termin\'{e}s   $(X,(Y_i,Z_i)_{i\N^*})$.
Pour tout entier $n\geq 2$, pour tout $j=1,\dots n-1$, on pose $\al_{j}^{n}={Z_j}^{p^{-n+j}}X^{p^{-n}}+{Y_j}^{p^{-n+j}}$ et   $K_n=k(X^{p^{-n}}, \al_{1}^{n},\dots, \al_{n-1}^{n})$. Puisque $K_n\subseteq K_{n+1}$ et $K_n/k$ est purement ins\'{e}parable, alors $K=\di\bigcup_{n\geq 2}K_n$ est un corps purement ins\'{e}parable sur $k$. Par convention, on identifie $K_0$ \`a $k$ et $K_1$ \`a $k(X^{p^{-\infty}})$.
\end{exe}

On v\'erifie imm\'ediatement que pour tous entier $n\geq 2$, pour tout $j=1,\dots n-1$,  ${({\al_{j}^{n+1}})}^p= {(\al_{j}^{n})}$. En outre, pour tout $n\geq 2$, $k(K_{n+1}^p)=K_n$, et par suite $K/k$ est relativement parfaite.
\begin{thm} {\label{thm15}} L'extension $K/k$ ci-dessus est $lq$-finie de degr\'{e} d'irrationalit\'{e} infini.
\end{thm}

Pour la d\'emonstration, nous aurons besoin des  r\'{e}sultats suivants.
De plus, tout le long de cette section, pour tout entier $t\leq s$,  $\zeta_{s}^{t}$ d\'esigne l'ensemble $\{1,\dots, s\}\setminus \{t\}$.
Il est clair que pour tout entier $s\geq 2$, pour tout $i\in \{1,\dots,s-1\}$, $\al_{i}^s\in k(X^{p^{-s}}, {Z_i}^{p^{-s}},{Y_i}^{p^{-s}})$ ; et par suite,
$K_s\subseteq k(X^{p^{-s}},({Y_j}^{p^{-s}}, {Z_j}^{p^{-s}})_{j\in \ze_{s-1}^{i}}, \al_{i}^s)$.
\begin{lem} {\label{lem5}} Pour tout $j\geq 1$, $Z_j$ est $p$-ind\'{e}pendant sur $K^p$, ou encore ${Z_j}^{p^{-1}}\not\in K$.
\end{lem}
\pre Supposons que ${Z_j} ^{p^{-1}}\in K=\di\bigcup_{n\geq 2}K_n$, donc il existe $s> j$ tel que ${Z_{j}}^{p^{-1}}\in K_{s}$.
 Comme ${Y_j}^{p^{-1}}={(\al_{j}^{s})}^{p^{s-j-1}}-{Z_j} ^{p^{-1}}X^{p^{-j-1}}$, alors ${Y_j} ^{p^{-1}}\in K_s\subseteq k(X^{p^{-s}},({Y_i}^{p^{-s}}, {Z_i}^{p^{-s}})_{i\in \ze_{s-1}^{j}}, \al_{j}^s)$. En outre, $k(X^{p^{-1}},{Y_1}^{p^{-1}},\dots,$ ${Y_{s-1}}^{p^{-1}}, $ $ {Z_1}^{p^{-1}},$ $\dots,{Z_{s-1}}^{p^{-1}})\subseteq k(X^{p^{-s}},({Y_i}^{p^{-s}}, {Z_i}^{p^{-s}})_{i\in \ze_{s-1}^{j}}, \al_{j}^s)$. Puisque $(X,(Y_i, Z_i)_{i\geq 1})$ est alg\'{e}briquement ind\'{e}pendante sur $Q$, on aura $2(s-1)+1=di(k(X^{p^{-1}}, {Y_1}^{p^{-1}},\dots,$ ${Y_{s-1}}^{p^{-1}},{Z_1}^{p^{-1}},\dots,{Z_{s-1}}^{p^{-1}})/k)\leq di(k(X^{p^{-s}},$ $({Y_i}^{p^{-s}},{Z_i}^{p^{-s}})_{i\in \ze_{s-1}^{j}}, \al_{j}^{s})/k)$ $\leq 2($ $s-2)+2$. Il en r\'esulte que $1\leq 0$, c'est une contradiction, et par suite   ${Z_j} ^{p^{-1}}\not \in K$. 
\begin{lem} {\label{lem6}} Pour tout $n\geq 2$, la famille $\{X^{p^{-n-1}},\al_{1}^{n+1},\dots, \al_{n}^{n+1}\}$ est une $r$-base de $K_{n+1}/K_n$
\end{lem}
\pre Il est aussit\^ot que $K_n(X^{p^{-n-1}}, \al_{1}^{n+1},\dots,\al_{n}^{n+1})=K_{n+1}$,  et $X^{p^{-n-1}}\not\in K_n$ puisque $n+1=o(X^{p^{-n-1}},k)> n=o_1(K_n/k)$. S'il existe $i\in\{1,\dots, n\}$ tel que
$\al_{i}^{n+1}\in K_n(X^{p^{-n-1}},\al_{1}^{n+1},\dots, \al_{i-1}^{n+1},$ $\al_{i+1}^{n+1},\dots, \al_{n}^{n+1})$, comme $k({K_{n+1}}^{p})=K_n$ et $K_{n+1}/k$ est d'exposant fini,  on en d\'{e}duit que $k(X^{p^{-n-1}}, \al_{1}^{n+1},$ $\dots,\al_{i-1}^{n+1},$ $\al_{i+1}^{n+1},$ $
\dots,\al_{n}^{n+1})=K_{n+1}$. Notamment,  $\al_{i}^{n+1}\in k(X^{p^{-n-1}}, \al_{1}^{n+1},\dots,\al_{i-1}^{n+1},$ $\al_{i+1}^{n+1},$ $\dots,\al_{n}^{n+1})\subseteq
k({X}^{p^{-n-1}}, ({Z_{j}}^{p^{-n-1}}, {Y_{j}}^{p^{-n-1}})_{j\in \ze_{n}^{i}})$, et donc ${Z_i}^{p^{-1}}$ appartient \`a
$k({X}^{p^{-n-1}}, $ $({Z_{j}}^{p^{-n-1}},$ $ {Y_{j}}^{p^{-n-1}}$ $)_{j\in \ze_{n}^{i}}, {Y_i}{^{p^{-1}}})$. Comme la famille $(X,Y_i,Z_i)_{i\N^*}$ est alg\'{e}briquement ind\'{e}pendante sur $Q$,  on obtiendra $2n+1=di(k(X^{p^{-1}}, $ ${Y_1}^{p^{-1}},\dots,{Y_n}^{p^{-1}},{Z_1}^{p^{-1}},$ $\dots,{Z_n}^{p^{-1}})/k)\leq di(k({X}^{p^{-n-1}}, ({Z_{j}}^{p^{-n-1}}, $ ${Y_{j}}^{p^{-1}}$ $)_{j\in \ze_{n}^{i}}, $ ${Y_i}{^{p^{-n-1}}}$ $)/k)=2(n-1)+2$, donc $1\leq 0$, c'est une contradiction. D'o\`u $\{X^{p^{-n-1}},\al_{1}^{n+1},\dots, \al_{n}^{n+1}\}$ est une $r$-base de $K_{n+1}/K_n$. \cqfd \sk

D'apr\`es la proposition ci-dessus, pour tout $j>2$,  $\{X^{p^{-j-1}}, \al_{1}^{j+1},\dots,\al_{j}^{j+1}\}$ et $\{X^{p^{-j}}, {(\al_{1}^{j+1})}^p,\dots,$ ${(\al_{j-1}^{j+1})}^p\}$ sont deux  $r$-bases respectivement de $K_{j+1}/K_j$ et $K_{j}/K_{j-1}$. Comme $K_j= k({K_{j+1}}^p)$ et $K_{j-1}= k({K_{j}}^p)$, d'apr\`es le lemme {\ref{lem2}}, pour tout $j>2$, pour tout $i\in \{1,\dots, j-1\}$, on a $o(\al_{i}^{j+1}, K_{j-1}(X^{p^{-j-1}},$ $\al_{1}^{j+1},\dots,\al_{i-1}^{j+1}, \al_{i}^{j+1},\dots, \al_{j-1}^{j+1}))=o_1(K_{j+1}/K_{j-1})=2$. En particulier, $K_{j-1}($ $X^{p^{-j-1}},$ $ \al_{1}^{j+1},\dots,\al_{j-1}^{j+1})\simeq K_{j-1}(X^{p^{-j-1}})\otimes_{K_{j-1}} K_{j-1}(\al_{1}^{j+1})\otimes_{K_{j-1}} \dots \otimes_{K_{j-1}} K_{j-1}(\al_{i-1}^{j+1})$.
\begin{lem} {\label{lem7}}
Pour tout $m\in\N$, pour tout $n\in \N$,  on a $k^{p^{-n}}\cap K_m=K_n$ si $m\geq n$, et  $k^{p^{-n}}\cap K_m=K_m$ si $m\leq n$.
\end{lem}
\pre  Le lemme est v\'{e}rifi\'e pour $n=0$. Soit $n\geq 1$,
supposons que le r\'{e}sultat est satisfait pour tout entier naturel $i\leq n-1$. Puisque $o_1(K_m/k)=m$,  il est imm\'{e}diat que
$k^{p^{-n}}\cap K_m=K_m$ si $m\leq n$. On est amen\'e  au cas o\`{u} $n<m$. Compte tenu de l'hypoth\`{e}se du r\'ecurrence, on a $k^{p^{-n+1}}\cap K_m=K_{n-1}$,
donc  $k^{p^{-n}}\cap K_m={K_{n-1}}^{p^{-1}}\cap K_m$. Comme  $o_1(K_n/K_{n-1})=1$ et $K_n\subseteq K_m$, ($n<m$), on a $K_n\subseteq {K_{n-1}}^{p^{-1}}
\cap K_m$. Si ${K_{n-1}}^{p^{-1}}\cap K_m\not=K_n$, ou encore s'il existe $\te\in {K_{n-1}}^{p^{-1}}\cap K_m $ tel que $\te \not\in K_n$, soit $j$
le plus grand entier tel que $\te \not\in K_j$, donc $\te\in K_{j+1}$. En particulier $n\leq j<m$, et $o(\te, K_{j-1})=1$, (car ${\te}^p\in K_{n-1}\subseteq K_{j-1}$). Pour all\'eger l'\'ecriture, on pose provisoirement $\al_1=X^{p^{-j-1}}$, et pour tout $i=2,\dots, j+1$, $\al_i= \al_{i-1}^{j+1}$. On va montrer ensuite que $\{\te, \al_{1},\dots,\al_{j}\}$
est une $r$-base de  $K_{j+1} /K_j$. Par hypoth\`ese,   $\te\not\in K_j$ et,  pour tout $i\in \{1,\dots,j\}$, on a
$\al_{i}\not\in K_j(\te, \al_{1},\dots,\al_{i-1})$ ;
sinon on aura $1<o(\al_{i} ,K_{j-1}(\al_{1},\dots,\al_{i-1}))\leq  o_1(K_j(\al_{1},
\dots,\al_{i-1},\te)/K_{j-1}(\al_{1},\dots,\al_{i-1}))\leq o(K_j(\te)/K_{j-1})=1 $, ($o(\te, K_{j-1})=1$ et $o_1(K_j/K_{j-1})=1$),  c'est
une contradiction.  D'o\`{u} $\{\te,\al_{1},\dots,\al_{j}\}$ est $r$-libre sur $K_j$. Comme $di(K_{j+1}/K_j)=j+1$ et $k({K_{j+1}}^p)=K_j$, alors  $\{\te,\al_{1},\dots,\al_{j}\}$ est une $r$-base de $K_{j+1}/K_j$. D'autre part, on a $o(\te,K_{j-1})=1$ et $K_{j-1}(\al_{1},\dots,\al_{j})/K_{j-1}$ est \'equiexponentielle d'exposant $2$, ce qui entraine $K_{j+1}\simeq K_{j-1}(\al_{1},\dots,\al_{j})\otimes_{K_{j-1}} K_{j-1}(\te)\simeq K_{j-1}(\al_{1})\otimes_{K_{j-1}} \dots \otimes_{K_{j-1}} K_{j-1}(\al_{j})\otimes_{K_{j-1}} K_{j-1}(\te)$. D'o\`u, $K_{j+1}/K_{j-1}$ est modulaire.
Toutefois, l'\'{e}quation de d\'{e}finition de $\al_{j+1}$ sur $K_{j-1}(\al_1, \al_{2},\dots,\al_{j})$ s'\'{e}crit : ${\al_{j+1}}^p=Z_{j}X^{p^{-j}}+Y_{j}=Z_{j}{\al_1}^p+Y_{j}$, d'apr\`{e}s le crit\`{e}re de modularit\'{e}, on obtient $Z_{j}, Y_{j}\in {K_{j+1}}^p\cap K_{j-1}\subseteq {K_{j+1}}^p$. On en d\'eduit que ${Z_{j}}^{p^{-1}}, {Y_{j}}^{p^{-1}}\in {K_{j+1}}$. C'est une contradiction avec le lemme {\ref{lem5}}. D'o\`{u} $k^{p^{-n}}\cap K_m=K_n$. \cqfd
\sk

{\bf Preuve du th\'{e}or\`{e}me {\ref{thm15}}.} Pour tout $n\geq 1$, on $K_n\subseteq K_{n+1}$ et $K=\di\bigcup_{m\geq1} K_m$, donc en vertu du lemme {\ref{lem7}}, $k^{p^{-n}}\cap K=\di\bigcup_{m\geq1}k^{p^{-n}}\cap K_m=K_n$. Comme pour tout $n\geq 1$, on a $K_n/k$ est finie, on en d\'{e}duit que $K /k$ est $lq$-finie non triviale (qui n'est pas $q$-finie).\cqfd
\begin{rem}
L'exemple ci-dessus est aussi  bon pour affirmer que la $lq$-finitude n'est pas respect\'{e} si on change le corps de base dans le sens ascendant. En effet, si on pose $L=k(X^{p^{-\infty}})$, on v\'erifie imm\'{e}diatement que :
\sk

\begin{itemize}{\it
\item[\rm{(i)}] $k(({\al_{s-1}}^s)_{s>1})\subseteq L^{p^{-1}}\cap K$.
\item[\rm{(ii)}] La famille $(\al_{s-1}^s)_{s>1}$ est $r$-libre sur $k$.}
\end{itemize}
\sk

D'o\`u, $L^{p^{-1}}\cap K /k$ n'est pas finie, et par suite $K/L$ n'est pas $lq$-finie.
En particulier, il est fort probable que $K/k$ et $L/k$ sont $lq$-finies, mais $L(K)/L$ ne l'est pas.
\end{rem}

D'une fa\c{c}on plus g\'en\'erale, et contrairement \`{a} la $q$-finitude, la $lq$-finitude n'est pas transitive comme le montre l'exemple ci-dessous.
Il est \`a noter que cet exemple  modifie l\'{e}g\`{e}rement les conditions de l'exemple {\ref{exe1}}.
\begin{exe} {\label{exe2}} Soient $Q$ un corps parfait de caract\'{e}ristique $p>0$, et $(X,(Z_i)_{i\in {\N}^{*}})$ une famille alg\'{e}briquement libre sur $Q$, et $k= Q(X,(Z_i)_{i\in\N^*})$ le corps des fractions rationnelles aux ind\'{e}termin\'{e}es $(X,(Z_i)_{i\in {\N}^{*}})$. Dans cette partie,  on se servira des notations suivantes :
\sk

\begin{itemize}{\it
\item[] $K_1 = k(\te_{1,1})$, avec $\te_{1,1}=X^{p^{-1}}$.
\item[] $K_2=k(\te_{2,1},\te_{2,2})$, avec $\te_{2,1}={(\te_{1,1})}^{p^{-1}}=X^{p^{-2}}$ et $\te_{2,2}={Z_1}^{p^{-1}} \te_{2,1}+{Z_2}^{p^{-1}}={Z_1}^{p^{-1}} X^{p^{-2}}+{Z_2}^{p^{-1}}$.
\item[] $K_3=k(\te_{3,1},\te_{3,2},\te_{3,3})$, avec $\te_{3,1}={(\te_{2,1})}^{p^{-1}}$, $ \te_{3,2}={(\te_{2,2})}^{p^{-1}}$, et $\te_{3,3}={Z_2}^{p^{-1}}\te_{3,2}+{Z_3}^{p^{-1}}$
\item[] Par r\'{e}currence, on pose $K_n=k(\te_{n,1},\dots,\te_{n,n})$, o\`u  $\te_{n,i}={(\te_{n-1,i})}^{p^{-1}}$ pour tout $i=1,\dots, n-1$, et $\te_{n,n}={Z_{n-1}}^{p^{-1}}\te_{n,n-1}+{Z_n}^{p^{-1}}$.}
\end{itemize}
\sk

Posons $K=\di\bigcup_{n\geq 1}K_n$ et $L=K({Z_1}^{p^{-\infty}})$.
\end{exe}

On v\'{e}rifie aussit\^{o}t que :
\begin{itemize}{\it
\item Pour tout $n\geq 1$, $k({K_{n+1}}^p)=K_n$. En particulier, $K/k$ est relativement parfaite.
\item $L\simeq k(X^{p^{-\infty}})\otimes_k (\otimes_k (k({Z_i}^{p^{-\infty}}))_{i\geq 1})$
\item Pour tout $n\geq 1$, pour tout $i=1,\dots, n$, $Z_i$ est $p$-ind\'{e}pendant sur ${K_n}^p$, ou encore ${Z_i}^{p^{-1}}\not\in K_n$. En outre, ${Z_{i}}^{p^{-1}}\not\in K$.
\item Pour tout $n\geq 1$, $(\te_{n,1},\dots,\te_{n,n})$ est une $r$-base de $K_n/k$ (on utilise le m\^eme raisonnement que l'exemple ci-dessus).}
\end{itemize}
\begin{thm} {\label{thm16}} $K /k$ et $L/K$ sont $lq$-finies, mais $L/k$ ne l'est pas.
\end{thm}

Comme dans l'exemple pr\'ec\'edent, on emploiera le lemme technique suivant:
\begin{lem} {\label{lem8}} Pour tous $n, m \in \N^*$, $k^{p^{-n}}\cap K_m=K_m$ si $n\geq m$, et $k^{p^{-n}}\cap K_m=K_n$ si $n\leq m$. En particulier, $K/k$ est $lq$-finie.
\end{lem}
\pre D\'{e}monstration analogue \`{a} celle du lemme {\ref{lem7}}. \cqfd \sk

{\bf Preuve du th\'{e}or\`{e}me {\ref{thm16}}} D'une part, on a $K/k$ et $L/K$ sont $lq$-finies. D'autre part, pour tout $n\geq 1$, on a $k^{p^{-n}}\cap L\simeq k(X^{p^{-n}})\otimes_k(\otimes_k (k({Z_i}^{p^{-n}}))_{i\geq 1})$,  et donc $L /k$ n'est pas $lq$-finie. \cqfd
\begin{rem}
Egalement, cet exemple peut servir pour montrer que le produit ne respecte pas la $lq$-finitude. Pour cela, on pose $L_1=k({Z_1}^{p^{-1}})$ et $L_2=k(({Z_i}^{p^{-1}})_{i\geq 2})$. Par construction, pour tous entier $n\geq 2$,  $\te_{n,n}={Z_{n-1}}^{p^{-1}}\te_{n,n-1}$ $+{Z_n}^{p^{-1}}$, et  pour tout entier non nul $i$,  on montre par r\'ecurence que ${Z_i}^{p^{-1}}\in L_1(K)$, ou encore $L_2\subseteq L_1(K)$. Comme $di(L_2/k)$ est infinie, alors $L_1(K)/k$ n'est pas $lq$-finie m\^eme si $K/k$ et $L_1/k$ sont $lq$-finies.
\end{rem}
\begin{pro} {\label{pr31}} Soient $K/k$ une extension purement ins\'{e}parable d'exposant non born\'{e}, et $H$ l'ensemble des sous-extensions d'exposant non born\'{e} de $K/k$. Si $K/k$ est $lq$-finie, alors $H$ est inductif
pour la relation d'ordre d\'efinie par $K_1\leq K_2$ si et seulement si $K_2\subseteq K_1$.
\end{pro}
\pre
$H\not=\emptyset$, puisque $K\in H$.  Soit $(K_n/k)_{n\in I}$ une famille totalement ordonn\'{e}e de sous-extensions d'exposant
non born\'{e} de $K/k$,  donc pour tout $j\in {\N}^*$, la famille ${(k^{p^{-j}}\cap K_n)}_{n\in I}$ est aussi totalement ordonn\'{e}e.
Comme $K/k$ est $lq$-finie, donc la famille $([ k^{p^{-j}}\cap K_n :k])_{n\in I}$ des entiers naturels  est totalement ordonn\'{e}e ;
et par suite elle admet un plus petit \'{e}l\'{e}ment que l'on note $[ k^{p^{-j}}\cap K_{n_j} :k]$. Soit $I_1=\{s\in I|\, K_s\subseteq  K_{n_j}\}$ et $I_2=\{s\in I|\, K_{n_j}\subseteq K_s\}$, donc pour tout $m\in I_1$,  $[k^{p^{-j}}\cap K_m :k]\leq [k^{p^{-j}}
\cap K_{n_j}:k]$. Compte tenu de la propri\'et\'e caract\'eristique  du plus petit \'el\'ement, on en d\'eduit que $
[k^{p^{-j}}\cap K_m :k]=[k^{p^{-j}}\cap K_{n_j}:k]$, et par suite pour tout $m\in I_1$, on a $k^{p^{-j}}\cap K_m=k^{p^{-j}}\cap K_{n_j}$. D'o\`u, $k^{p^{-j}}\cap  (\di\bigcap_{m\in I_1}
K_m)= \di\bigcap_{m\in I_1} k^{p^{-j}}\cap K_m=k^{p^{-j}}\cap K_{n_j}$. Si on pose $L=\di\bigcap_{i\in I} K_i$, il est clair que
$L=\di\bigcap_{m\in I_1} K_m$. Comme $K/k$ est d'exposant non born\'e et $k^{p^{-j}}
\cap K_{n_j}\subseteq K_m$ pour tout $m\in I_1$, donc $o_1(k^{p^{-j}}
\cap L/k)=o_1(k^{p^{-j}}\cap K_{n_j}/k)=j$. D'o\`u, $L/k$ est d'exposant
non born\'{e}, et pour tout $n\geq 1$, $K_n\leq L$. Il en r\'{e}sulte que $H$ est inductif. \cqfd \sk

Comme cons\'equence imm\'ediate, on a :
\begin{cor} {\label{cor12}}Toute extension $lq$-finie d'exposant non born\'e admet une sous-extension minimale $M/k$ d'exposant non born\'e. En outre, $M/k$ est relativement parfaite.
\end{cor}
\pre Imm\'ediat. \cqfd
\begin{rem} La condition de $lq$-finitude est n\'{e}cessaire pour que $H$ soit inductif, comme le montre l'exemple suivant :
\end{rem}
\begin{exe} {\label{exe3}}Soient $Q$ un corps parfait de caract\'{e}ristique $p>0$, et $k=Q((X_i)_{i\in {\N}^*})$ le corps des fractions rationnelles aux ind\'{e}termin\'{e}es  $(X_i)_{i\in {\N}^*}$. Pour tout $n\geq 1$, notons  $K_n=k(((X_i)^{p^{-i}})_{i\geq n})$.
\end{exe}

On v\'erifie aussit\^ot que la famille $(K_n)_{n\in \N^*}$ de sous-extensions d'exposant non born\'{e} de $K_1/k$ est  totalement ordonn\'{e}e. Comme $((X_i)^{p^{-i}})_{i\in {\N}^*}$ est une $r$-base modulaire de $K_1=k(((X_i)^{p^{-i}})_{i\in {\N}^*})$ sur $k$, alors $\di\bigcap_{n\geq 1} K_n$  se r\'{e}duit \`{a} $k$, et donc $H$ n'est pas inductif.
\sk

Une autre cons\'equence de la proposition {\ref{pr31}} est le r\'esulatt suivant :
\begin{pro} {\label{pr32}} La cl\^{o}ture relativement parfaite $K_r$ d'une extension $lq$-finie d'exposant non born\'e  $K/k$  est non trivial. Plus pr\'{e}cis\'{e}ment,  $K_r/k$ est d'exposant non born\'{e} si $K /k$ l'est.
\end{pro}
\pre  D'apr\`{e}s la proposition {\ref{pr31}},  l'ensemble des sous-extensions  d'exposant non born\'{e} de $K/k$ admet une sous-extension $m/k$ minimale. N\'{e}cessairement,  $m/k$ est relativement parfaite, sinon $k(m^p) /k$ serait meilleure que $m$, contradiction ; et donc $K_r/k$ est d'exposant non born\'{e} ($m\subseteq K_r$). \cqfd \sk

Lorsque $K/k$ est $q$-finie, on peut situer $rp(K/k)$ par rapport \`a $K$, en particulier, on a $K/rp(K/k)$ est finie.  Cependant, dans le cas de la $lq$-finitude l'emplacement de la cl\^oture relativement parfaite varie d'une extension \`a l'autre.

A cet \'egard,  chacun des  exemples {\ref{exe1}} et {\ref{exe2}} ci-dessus pr\'esente une extension $K/k$ $lq$-finie et relativement parfaite (et donc $K/rp(K/k)$ est d'exposant fini). Toutefois, voici un exemple o\`u $K/rp(K/k)$ est d'exposant infini.
\begin{exe} {\label{exe4}} Soient $Q$ un corps parfait de caract\'{e}ristique  $p>0$, et $k=Q(X, (Y_i, $ $Z_i)_{i\in {\N}^*})$ le corps des fractions rationnelles aux ind\'{e}termin\'{e}es $(X, (Y_i,$ $ Z_i)_{i\in {\N}^*})$. Pour tout $n\geq 1$, on pose $\te_n= {{Y_n}}^{p^{-1}}X^{p^{-n-1}}+{Z_n}^{p^{-1}}$. Egalement, on note $K=k(X^{p^{-\infty}}, (\te_i)_{i\in{\N}^*})$,  et pour tout $n\in \N$,  $k_n=k^{p^{-n}}\cap K $. Par convention,  $\te_0$ d\'{e}signe $X^{p^{-1}}$.
\end{exe}
\begin{pro} {\label{pr33}} Sous les conditions ci-dessus, on a $k_1=k(\te_0)=k(X^{p^{-1}})$, et pour tout $n\geq 2$, $k_n=k(X^{p^{-n}}, \te_1,\dots, \te_{n-1})$. En particulier, $K/k$ est $lq$-finie.
\end{pro}
\begin{rem} Pour regrouper les deux conditions, la proposition ci-dessus peut s'\'enonc\'ee comme suit :  pour tout $n\geq 1$, $k_n=k(X^{p^{-n}},\te_0,\dots, \te_{n-1})$.
\end{rem}
\pre La d\'emonstration se fait par r\'ecurrence, et utilise (plus particuli\`erement au rang $1$) les m\^{e}mes techniques de raisonnement en passant d'un niveau \`{a} l'autre. Pour cela, on suppose que l'on a  $k_n=k(X^{p^{-n}}, \te_0,\dots, \te_{n-1})$ avec $n\geq 1$. Il est clair que $k_{n+1}={k_n}^{p^{ -1}}\cap K$ et $k_n(X^{p^{-n-1}},\te_n)\subseteq k_{n+1}$.  Toutefois, s'il existe $\al\in k_{n+1}$ tel que $\al\not\in k_n(X^{p^{-n-1}},\te_n)$ ; alors  $\al\not \in k_n(X^{p^{-\infty}},\te_n)$. Sinon, comme
$k_n(\te_n)(X^{p^{-\infty}})/k_n(\te_n)$ est $q$-simple et $o(\al , k_n(\te_n))=o({X}^{p^{-n-1}},$ $ k_n(\te_n))=1$, alors n\'{e}cessairement $k_n(\te_n,\al)=k_n(\te_n,$ $X^{p^{-n-1}})$, c'est une contradiction. Soit $i$ le
plus grand entier tel que $\al\not\in k_n(X^{p^{-\infty}},$ $ \te_n, \dots, \te_i)$. Pour all\'eger l'\'{e}criture, on pose
$L=k(X^{p^{-\infty}}, \te_1,\dots,\te_{i+1})=k_n(X^{p^{-\infty}}, $ $\te_n,\dots,\te_{i+1})$, donc
$\al\in L$ et $o(\al, k_n(X^{p^{-\infty}}, \te_n,\dots,\te_{i}))=o(\te_{i+1}, k_n(X^{p^{-\infty}},$ $ \te_n,\dots,\te_{i}))=1$. Il en r\'esulte que
$L= k_n(X^{p^{-\infty}}, \te_n,\dots,\te_{i})\otimes_{k_n} k_{n}(\al)$. Or, par construction, pour tout $j\geq 1$, $\te_j={Y_j}^{p^{-1}}X^{p^{-j-1}}+{Z_j}^{p^{-1}}$,
et donc $k(\te_j)\subseteq k({Y_j}^{p^{-1}},X^{p^{-j-1}},{Z_j}^{p^{-1}})\subseteq k({Y_j}^{p^{-1}},X^{p^{-\infty}},$ ${Z_j}^{p^{-1}})$. Notons $L_1=L(({Y_{j}}^{p^{-1}}, {Z_j}^{p^{-1}}$ $)_{n+1\leq j\leq i})$. Si $\te_{i+1}\in L_1$, comme $\te_{i+1}={Y_{i+1}}^{p^{-1}}X^{p^{-i-2}}+{Z_{i+1}}^{p^{-1}}$, alors $L_1($ ${Y_{i+1}}^{p^{-1}}$ $)=L_1({Z_{i+1}}^{p^{-1}})$, et donc ${Z_{i+1}}^{p^{-1}}\in k(X^{p^{-\infty}}, ({Y_j}^{p^{-1}},$ ${Z_j}^{p^{-1}}$ $)_{1\leq j\leq i},$ ${Y_{i+1}}^{p^{-1}})$, c'est une contradiction avec le fait que $(X, (Y_i, Z_i)_{i\in {\N}^*})$ sont al\'ebriquement ind\'ependent sur $Q$ et $k=Q((X, (Y_i, Z_i)_{i\in {\N}^*}))$. D'o\`u,
\begin{eqnarray*}
L_1 &=&k_n(X^{p^{-\infty}}, ({Y_{j}}^{p^{-1}}, {Z_j}^{p^{-1}})_{n+1\leq j\leq i}, \te_{i+1}))\\
& =& k_n(X^{p^{-\infty}}, ({Y_{j}}^{p^{-1}}, {Z_j}^{p^{-1}})_{n+1\leq j\leq i}, \al) \\
&\simeq& k_n(X^{p^{-\infty}})\otimes_{k_n} (( k_n({Y_{j}}^{p^{-1}})\otimes_{k_n} k_n ( {Z_j}^{p^{-1}}))_{n+1\leq j\leq i})) \otimes_{k_n} k_n(\al).
\end{eqnarray*}
On en d\'{e}duit que $L_1 /k_n$ est modulaire. En outre, ${L_1}^p$ et $k_n$ sont $k_n\cap {L_1}^p$-lin\'{e}airement disjointes.
Comme $X^{p^{-i-1}}\not\in k_n$ ($n\leq i$),
alors $(1,X^{p^{-i-1}})$ est $k_n$-libre, en particulier $k_n\cap {L_1}^{p}$-libre. Compl\'{e}tons ce syst\`{e}me en une base $B$ de ${L_1}^p$ sur $k_n\cap {L_1}^p$,
en vertu de la lin\'{e}arit\'{e} disjointe $B$ est aussi une base de $k_n({L_1}^p)$ sur $k_n$. Or, l'\'equation de d\'efinition de $\te_{i+1}$ par rapport \`a $L_1$ s'\'ecrit : ${\te_{i+1}}^p=Y_{i+1}X^{p^ {-i-1}}+Z_{i+1}$,
par identification on aura $Y_{i+1}, Z_{i+1}\in k_n\cap {L_1}^p$. En particulier, ${Y_{i+1}}^{p^{-1}}, {Z_{i+1}}^{p^{-1}}\in L_1$, et donc
$k(X^{p^{-\infty}}, {({Y_j}^{p^{-1}},{Z_j}^{p^{-1}})}_{1\leq j\leq i+1})\subseteq L({({Y_j}^{p^{-1}},{Z_j}^{p^{-1}})}_{1\leq j\leq i}),\te_{i+1})= k(X^{p^{-\infty}}, {({Y_j}^{p^{-1}},{Z_j}^{p^{-1}})}_{1\leq j\leq i}),\te_{i+1})$. Il en r\'{e}\-s\-u\-lte que
$2(i+1)+1=di(k(X^{p^{-\infty}}, $ $({Y_j}^{p^{-1}},{Z_j}^{p^{-1}}$ $)_{1\leq j\leq i+1})/k)\leq di(k(X^{p^{-\infty}}, ($ ${Y_j}^{p^{-1}},$ ${Z_j}^{p^{-1}}$ $)_{1\leq j\leq i},\te_{i+1})/k)\leq 2(i+1)$ ;
d'o\`{u} $1\leq 0$, c'est une contradiction. Par suite, $k_{n+1}=k(X^{p^{-n-1}},\te_1,$ $\dots, \te_n)$. \cqfd
\begin{pro} {\label{pr34}}Sous les hypoth\`eses ci-dessus, $rp(K/k)=k(X^{p^{-\infty}})$.
\end{pro}
\pre Il suffit de remarquer que  $(\te_i)_{i\geq 1}$ est une $r$-base modulaire de $K/k($ $X^{p^{-\infty}}$ $)$, et donc $k(X^{p^{-\infty}})\subseteq rp(K/k)\subseteq \di\bigcup_{i\in \N} k(K^{p^i})=\di\bigcup_{i\in \N} k(X^{p^{-\infty}})(K^{p^i})=k(X^{p^{-\infty}})$. D'o\`u $rp(K/k)=k(X^{p^{-\infty}})$ et $K/k(X^{p^{-\infty}})$ est d'exposant non born\'e. 
\sk

Non seulement $K/rp(K/k)$ est d'exposant non born\'e, mais $K/rp(K/k)$ n'est pas $lq$-finie.
\begin{pro} {\label{pr35}} Toute extension $lq$-finie et modulaire est $q$-finie.
\end{pro}
\pre Il suffit de remarquer que si $K/k$ est modulaire, on a $di(k^{p^{-n}}\cap K/k)=di(k^{p^{-1}}\cap K/k)<+\infty$ pour tout $n\geq 1$. \cqfd
\section{Extensions absolument $lq$-finies}

Soit $K/k$ une extension purement ins\'eparable. Si $K/k$ est $q$-finie, alors pour toute sous-extension  $L/k$ de $K/k$, $K/L$ est aussi $q$-finie. Toutefois, d'apr\`es l'exemple {\ref{exe4}}; cette proposition est g\'en\'eralement fausse dans le cas de la $lq$-finitude.
Cela nous am\`ene \`a \'etdier de p\`es la stabilit\'e au sens absolue de la $lq$-finitude au sein des extensions $lq$-finie. L'\'etude de cette propri\'et\'e fait l'objet de cette section.
\begin{pro} {\label{pr36}}Soit $K/k$ une extension purement ins\'eparable. Les assertions suivantes sont \'equivalentes :
\sk

\begin{itemize}{\it
\item[\rm{(1)}] Pour toute sous-extension $L/k$ de $K/k$, $K/L$ est $lq$-finie.
\item[\rm{(2)}] Toute sous-extension $L/k$ de $K/k$ satisfait $L/k(L^p)$ est finie.
\item[\rm{(3)}] Pour toute sous-extension $L/k$ de $K/k$, on a $L/rp(L/k)$ est finie.}
\end{itemize}
\end{pro}
\pre $1\Rightarrow 2\Leftrightarrow 3 $ est imm\'ediate, il suffit de remarquer que $L\subseteq {(k(L^p))}^{p^{-1}}\cap K$ et $K/k(L^p)$ est $lq$-finie pour toute sous-extension $L/k$ de $K/k$, et d'appliquer la proposition {\ref{arp1}}.
Egalement, si l'item $2$ est v\'erifi\'e, donc  pour tout $n\in \N^*$, on aura $L_n/k({L_n}^p)$ est finie, o\`u $L_n=L^{p^{-n}}\cap K$. En particulier, $di(L_n/L({L_n}^p))$ est fini. Comme $L_n/L$ est d'exposant fini, alors $di(L_n/L)=di(L_n/L({L_n}^p))<+\infty$,  ou encore $L_n/L$ est finie pour tout $n\in \N^*$.\cqfd
\begin{df}Une extension qui v\'erifie l'une des conditions \'equivalentes de la proposition ci-dessus s'appelle extension absolument $lq$-finie.
\end{df}

On v\'{e}rifie aussit\^{o}t que :
\sk

\begin{itemize}{\it
\item Toute extension absolument $lq$-finie est $lq$-finie.
\item Toute extension $q$-finie est absolument $lq$-finie.
\item $K/k$ est absolument $lq$-finie si et seulement si toute sous-extension $L_2/L_1$ de $K/k$ est $lq$-finie. En particulier, toute sous-extension d'une extension absolument $lq$-finie est absolument $lq$-finie.}
\end{itemize}
\sk

Voici un exemple qui montre que $lq$-finitude est distincte de la $lq$-finitude absolut.
\begin{exe} {\label{exe5}} Reprenons l'exemple {\ref{exe4}}, rappelons que $Q$ d\'esigne toujours un corps parfait  de caract\'{e}ristique $p>0$ et $k=Q(X, (Y_i,Z_i)_{i\geq 1})$ le corps des fractions rationnelles aux ind\'{e}termin\'{e}es  $(X, (Y_i,Z_i)_{i\geq 1})$.  Rappelons aussi que $\te_i={Y_i}^{p^{-1}}X^{p^{-i-1}}+{Z_i}^{p^{-1}}$ et $K=k(X^{p^{-\infty}}, (\te_i)_{i\geq 1})$.
\end{exe}

D'apr\`{e}s l'exemple {\ref{exe4}} $K/k$ est $lq$-finie. Comme $(\te_i)_{i\geq 1}$ est une $r$-base de $K/k(K^p)$, on en d\'{e}duit que $K/k$ n'est pas absolument $lq$-finie.
\begin{thm} {\label{thm19}} Une extension purement ins\'{e}parable $K/k$ est absolument $lq$-finie si et seulement si toute suite d\'{e}croissante de sous-extensions de $K/k$ est stationnaire.
\end{thm}
\pre  Supposons que la condition suffisante est v\'{e}rifi\'ee, et consid\'{e}rons  une sous-extension  $L/k$ de $K/k$. Soit $G$ une $r$-base de $L/k$. Si $|G|$ n'est pas fini, il existe une suite $(\al_n)_{n\in{\N}^*}$ d'\'{e}l\'{e}ments deux \`{a} deux distincts de $G$. Pour tout $n\geq 1$, posons $K_n=k(G\setminus \{\al_1,\dots, \al_n\}))$. Il est clair que la suite $(K_n/k)_{n\geq 1}$ de sous-extensions de $K/k$ est   strictement d\'{e}croissante (pour l'inclusion), contradiction avec l'hypoth\`{e}se de la condition suffisante. D'o\`{u} $|G|$ est fini, ou encore $L/k(L^p)$ est finie ; et donc $K/k$ est absolument $lq$-finie.
Inversement, supposons que  $K/k$ est absolument $lq$-finie, et soit $(K_n/k)_{n\geq 1}$ une suite d\'{e}croissante de sous-extensions  de $K/k$.  Pour tout $n\geq 1$, posons $L_n=rp(K_n/k)$. Comme $K/k$ est absolument $lq$-finie, pour tout $n\geq 1$, on a $K_n/k({K_n}^p)$ est finie. Par suite, d'apr\`es la proposition {\ref{arp1}},  il existe $e_n\in\N$, tel que $L_n=k({K_n}^{p^{e_n}})=\di\bigcap_{i\geq 1}k({K_n}^{p^i})$. Puisque la suite $(K_n/k)_{n\geq 1}$ est d\'{e}croissante, il en est de m\^{e}me de $(L_n/k)_{n\geq 1}$. On distingue deux cas :
\sk

1-ier cas : S'il existe $n_0\in \N$ tel que $L_n=L_{n_0}$ pour tout $n\geq n_0$, alors  $K_n/L_{n_0}$ sera finie pour tout $n\geq n_0$. En vertu de la monotonie de $(K_n)_{n\geq 1}$, on en d\'{e}duit que la suite des entiers naturels $([K_n :L_{n_0}])_{n\geq n_0}$ est d\'{e}croissante, et donc stationnaire ; c'est-\`{a}-dire il existe $e\geq n_0$ tel que  $[K_n :L_{n_0}]=[K_e :L_{n_0}]$  pour tout $n\geq e$. Or, pour tout $n\geq 1$, $K_{n+1}\subseteq K_n$, donc $K_{n+1}=K_n$ pour tout $n\geq e$.
\sk

2-i\`eme cas : Si $(L_n)_{n\geq 1}$ n'est pas stationnaire, on est amen\'e \`{a} consid\'{e}rer que la suite $(L_n)_{n\geq 1}$ est strictement d\'{e}croissante, et par suite on va construire par r\'{e}currence une suite $(\al_i)_{i\in \N^*}$ d\'el\'ements de $K$  telle que $\al_i\in L_i$ et $\al_i\not\in k(\al_1,\dots,\al_{i-1})(L_{i+1})$ pour tout $i\geq 1$.
Le cas o\`u $i=1$ est trivial, il suffit de choisir  $\al_1$ dans $L_1\setminus L_2$. De plus, $k(\al_1)L_2\subseteq L_1$ (strictement). Consid\'erons maintenant un entier naturel $n$ distinct de $0$ et $1$, et supposons que cette propri\'{e}t\'{e} est satisfait pour tout $i=1,\dots, n$. Comme $L_i/k$ est relativement parfaite pour tout $i\geq 1$ et, $(L_i)_{i\geq 1}$ est strictement d\'ecroissante, on en d\'{e}duit que $(k(\al_1,\dots, \al_n)L_i)_{i\geq 1}$ est aussi strictement d\'{e}croissante ; et donc il existe $\al_{n+1}\in L_{n+1}$ tel que $\al_{n+1}\not\in k(\al_1,\dots, \al_n)L_{n+2}$. Dans la suite, on pose
$G=(\al_1,\dots,\al_n,\dots)$ et $G_n=(\al_1,\dots,\al_n)$, et on va montrer que $G$ est $r$-libre sur $k$. Comme $G=\di \bigcup_{n\geq 1} G_n$ et $G_i\subseteq G_{i+1}$, il suffit de montrer que $G_n$ est $r$-libre sur $k$ pour tout $n\geq 1$. S'il existe $i\in \{1,\dots, n\}$ tel que $\al_i\in k(G_n\setminus \{\al_i\})$, donc en particulier $\al_i\in k(\al_1,\dots,\al_{i-1})K_{i+1}$  (puisque pour tout $j\geq i+1$, $k(\al_1,\dots, \al_{i-1})K_j\subseteq k(\al_1,\dots, \al_{i-1})(K_{i+1})$), c'est une contradiction par construction. Par suite, $di(k(G)/k)=di(k(G)/k({(k(G))}^p))=+\infty$ ; et donc $K/k$ n'est pas absolument $lq$-finie, absurde. \cqfd \sk

Soit $K/k$ une extension purement ins\'eparable. Si $K/k$ est $q$-finie, il est aussit\^ot que toute sous-extension $L/k$ de $K/k$ est finie sur $k(L^p)$, donc il est fort probable que la r\'eciproque soit aussi vraie. Autrement dit la $lq$-finitude absolut est synonyme de la $q$-finitude. Toutefois, en utilisant  la propri\'et\'e caract\'eristique de la $lq$-finitude absolut,  voici un exemple type d'extension absolument $lq$-finie qui n'est pas $q$-finie. Par ailleurs, les deux notions se confondent dans le cas de la modularit\'e.
\subsection{Existence effective de $lq$-finitude absolut}
\begin{exe} {\label{exe6}}
Soient $Q$ un corps parfait de caract\'{e}ristique $p>0$, et $(X,(Z_{i}^{j}, Y_{i}^{j}$ $))_{(i,j)\in {\N^*}\times (\N^*\setminus \{1\})}$ une famille alg\'{e}briquement ind\'{e}pendante  sur $Q$. Soit $k=Q((X,(Z_{i}^{j},$ $ Y_{i}^{j}))_{(i,j)\in {\N^*}\times (\N^*\setminus \{1\})})$ le corps des fractions rationnelles aux ind\'{e}\-t\-e\-rmin\'{e}es $(X,(Z_{i}^{j},$ $ Y_{i}^{j})_{(i,j)\in \N^*\times (\N^*\setminus \{1\})})$. Soit $(a_n)_{n\in \N^*}$ une famille d'\'el\'ements de $Q(X^{p^{-\infty}})$ telle que $Q(a_n)=Q(X^{p^{-n}})$ pour tout $n\in \N^*$. En particulier, $k(a_n)=k(X^{p^{-n}})$, et par suite pour tout $n\in N^*$, $o(a_n,k)=o(X^{p^{-n}},k)=n$.
Dans la suite on se servira des notations suivantes :
\sk

\begin{itemize}{\it
\item[] $K_1=\di\bigcup_{i\geq 1} k(\te_{i}^{1})$, avec $\te_{i}^{1}=a_i$ pour tout $i\in N^*$.
\item[] $K_2=\di\bigcup_{i\geq 1} K_1(\te_{i}^2)$, avec $\te_{1}^{2}={(Z_{1}^{2})}^{p^{-1}}\te_{2}^{1}+{(Y_{1}^{2})}^{p^{-1}}$,
et pour tout $i\geq 2$, $\te_{i}^{2}={(Z_{i}^{2})}^{p^{-1}}\te_{2i}^{1}+{(\te_{i-1}^{2})}^{p^{-1}}+{(Y_{i}^{2})}^{p^{-1}}$.
\item[] Par r\'{e}currence, on pose : $K_n= \di\bigcup_{i\geq 1} K_{n-1}(\te_{i}^n)$, avec $\te_{1}^n= {(Z_{1}^n)}^{p^{-1}}\te_{2}^{n-1}+(Y_{1}^{n})^{p^{-1}}$, et pour tout $i\geq 2$, $\te_{i}^{n}=(Z_{i}^n)^{p^{-1}}\te_{2i}^{n-1}+{(\te_{i-1}^{n})}^{p^{-1}}+ (Y_{i}^{n})^{p^{-1}}$.}
\end{itemize}
\sk

On note \'egalement $K=\di\bigcup_{i\geq 1} K_n$, et pour des raisons de formulation,  on pose par convention $K_0=k$, et pour tout $i\in N$, $\te_{0}^{i}=0$.
\end{exe}

On v\'{e}rifie aussit\^{o}t que :
\sk

\begin{itemize}{\it
\item[$\bullet$] Pour tout $(i,n)\in {\N^*}\times (\N\setminus \{0,1\})$, $$\te_{i}^{n}=(Z_{i}^{n})^{p^{-1}}\te_{2i}^{n-1}+\dots+ (Z_{1}^{n})^{p^{-i}}(\te_{2}^{n-1})^{p^{-i+1}}+(Y_{i}^{n})^{p^{-1}}+\dots+ (Y_{1}^{n})^{p^{-i}}.$$
\item[$\bullet$] Pour tout $n\in \N^*$, pour tout $i\in \N$, on a $K_{n-1}(\te_{i}^{n})\subseteq K_{n-1}(\te_{i+1}^{n})$, et donc $K_n/k$ est un corps commutatif. En outre, $K/k$ est purement ins\'eparable, et  pour tout $n\in \N^*$, $K_n/K_{n-1}$ est $q$-simple.}
\end{itemize}

\begin{thm} {\label{thm22}}
$(K_i/k)_{i\geq 1}$ sont les seules sous-extensions d'exposant non born\'{e} de $K/k$ \`{a} une extension finie pr\`{e}s ; c'est-\`{a}-dire pour toute sous-extension propre $L/k$  d'exposant non born\'{e} de $K/k$, il existe $i\in \N^*$ tel que $K_i\subseteq L$ et $L/K_i$ est finie.
\end{thm}

Comme application type de ce th\'eor\`eme, on a :
\begin{pro} {\label{pr38}} Toute suite d\'{e}croissante de sous-extensions de $K /k$ est stationnaire.
\end{pro}
\pre Soit $(F_i/k)_{i\geq 1}$ une suite d\'{e}croissante de sous-extensions de $K/k$. D'apr\`{e}s le th\'eor\`eme pr\'ec\'edent, pour tout $j\in {\N}^*$, il existe $i_j\in \N$ tel que $K_{i_j}\subseteq F_j$ et $F_j/K_{i_j}$ est finie. Comme $F_{j+1}\subseteq F_j$, alors $\di\bigcap_{n\in {\N}} k({F_{j+1}}^{p^n})=K_{i_{j+1}}\subseteq \di\bigcap_{n\in {\N}} k({F_j}^{p^n})=K_{i_j}$. D'o\`u, il existe $s\in \N^*$,  pour tout $j\geq s$, $K_{i_j}=K_{i_s}$. En particulier, pour tout $j\geq s$, $F_j/K_{i_s}$ est finie. Ainsi, pour tout $j\geq s$, $[F_{j+1} :K_{i_s}]\leq [F_j : K_{i_s}]\leq [F_s :K_{i_s}]$ ; ou encore la suite d'entiers naturels $([F_j :K_{i_s}])_{j\geq s}$ est d\'{e}croissante, donc stationnaire \`{a} partir d'un entier $n_0$, et par suite pour tout $n\geq n_0$, $F_n=F_{n_0}$. \cqfd \sk

Comme cons\'equence imm\'ediate, on a :
\begin{thm} {\label{thm23}} Sous les m\^emes hypoth\`eses ci-dessus, on a $K/k$ est absolument $lq$-finie.
\end{thm}
\pre Due au th\'eor\`eme {\ref{thm19}}. \cqfd \sk

Pour la preuve du  th\'{e}or\`{e}me {\ref{thm22}} ci-dessus, nous aurons besoin des r\'{e}sultats suivants :
\sk

Tout d'abord, d\'esormais,  pour tout $(i,s)\in (\N\setminus \{0,1\})\times \N^*$, on pose $\La_{s}^{i}=\{(j,l)\in {\N}^2| \, 2\leq j\leq i$ et $1\leq l\leq 2^{i-j}s\}$.
\begin{lem} {\label{lem10}} Pour tout $i\geq 2$, pour tout $s\geq 1$, on a
$\te_{s}^{i}\in Q(\te_{{2}^{i-1}s}^{1},({(Z_{l}^{j})}^{p^{-2^{i-j}s}},$ ${(Y_{l}^j)}^{p^{-2^{i-j}s}})_{(j,l)\in \Lambda_{s}^{i}})$.
\end{lem}
\pre On va utiliser une d\'{e}monstration par r\'{e}currence.  Par construction, pour tout $s\in\N^*$, $\te_{s}^2={(Z_{s}^{2})}^{p^{-1}}\te_{2s}^{1}+\dots+ {(Z_{1}^{2})}^{p^{-s}}{(\te_{2}^{1})}^{p^{-s+1}}+{(Y_{s}^{2})}^{p^{-1}}+\dots+{(Y_{1}^{2})}^{p^{-s}}$, d'o\`u $\te_{s}^{2}\in Q(\te_{2s}^{1},\dots, {(\te_{2}^1)}^{p^{-s+1}}, {(Z_{s}^{2})}^{p^{-1}},\dots, $ ${(Z_{1}^{2})}^{p^{-s}}, $ ${(Y_{s}^{2})}^{p^{-1}},$ $\dots, {(Y_{1}^{2})}^{p^{-s}})\subseteq Q(\te_{2s}^{1}, {(Z_{s}^{2})}^{p^{-s}},\dots, {(Z_{1}^{2})}^{p^{-s}},{(Y_{s}^{2})}^{p^{-s}},\dots, {(Y_{1}^{2})}^{p^{-s}})$, et par suite le lemme est v\'erifi\'e pour le premier rang. Supposons que la propri\'{e}t\'{e} de r\'{e}currence s'applique jusqu'\`{a} l'ordre $i>1$. Egalement, pour tout $s\in \N^*$,  on a $\te_{s}^{i+1}={(Z_{s}^{i+1})}^{p^{-1}}\te_{2s}^{i}+\dots+{(Z_{1}^{i+1})}^{p^{-s}}{(\te_{2}^i)}^{p^{-s+1}}+$ ${(Y_{s}^{i+1})}^{p^{-1}}+\dots+{(Y_{1}^{i+1})}^{p^{-s}}$, et donc $\te_{s}^{i+1}\in Q(\te_{2s}^{i},\dots,$ ${(\te_{2}^{i})}^{p^{-s+1}}, $ ${(Z_{s}^{i+1})}^{p^{-1}},$ $\dots,{(Z_{1}^{i+1})}^{p^{-s}}, {(Y_{s}^{i+1})}^{p^{-1}},\dots,$ ${(Y_{1}^{i+1})}^{p^{-s}})$. Or, d'apr\`{e}s l'hypoth\`{e}se de r\'{e}currence, pour tous $r\in \{1,\dots, s\}$, $\te_{2r}^{i}\in Q(\te_{2^{i-1}.2r}^{1},({(Z_{l}^{j})}^{p^{-2^{i-j}2r}}, $ ${(Y_{l}^{j})}^{p^{-2^{i-j}2r}}$ $)_{(j,l)\in \La_{2r}^{i}})$. Comme $Q$ est parfait,  pour tout $n\in\N$, on a ${Q}^{p^{-n}}=Q$ ; et par suite ${(\te_{2r}^{i})}^{p^{-(s-r)}}\in Q({(\te_{2^{i}r}^{1})}^{p^{-(s-r)}},$ $
({(Z_{l}^{j})}^{p^{-2^{i+1-j}r-(s-r)}}, {(Y_{l}^{j})}^{p^{-2^{i+1-j}r-(s-r)}})_{(j,l)\in \La_{2r}^{i}})$ pour tous $r\in \{1,\dots, s\}$. Pui\-sque $Q(\te_{2^{i}r}^{1})=Q(X^{p^{{-2^i}r}})$ et $2^{i+1-j}s=2^{i+1-j}(r+(s-r))=2^{i+1-j}r+2^{i+1-j}(s-r)\geq 2^{i+1-j}r+(s-r)$, on en d\'{e}duit les relations ci-dessus:
\sk

\begin{itemize}{\it
\item[$\bullet$] $Q({(\te_{2^{i}r}^{1})}^{p^{-(s-r)}})=Q(X^{p^{{-2^i}r-(s-r)}})\subseteq Q(X^{p^{{-2^i}s}})=Q((\te_{2^{i}s}^{1}))$,
\item[$\bullet$]   ${(Z_{l}^{j})}^{p^{-2^{i+1-j}r-(s-r)}}\in Q({(Z_{l}^{j})}^{p^{-2^{i+1-j}s}})$, et ${(Y_{l}^{j})}^{p^{-2^{i+1-j}r-(s-r)}}$ \'el\'ement de  $Q({(Y_{l}^{j})}^{p^{-2^{i+1-j}s}})$.}
\end{itemize}
\sk

D'o\`{u},
$\te_{s}^{i+1}\in Q(\te_{2^{i}s}^{1}, ({(Z_{l}^{j})}^{p^{-2^{i+1-j}s}}, {(Y_{l}^{j})}^{p^{-2^{i+1-j}s}})_{(j,l)\in \La_{s}^{i+1}})$. \cqfd \sk

 Dans toute la suite,  pour tout $n\geq 2$, on pose $S_n=k(X^{p^{-\infty}}, ({(Z_{i}^{j})}^{p^{-\infty}}, $ $(Y_{i}^{j})^{p^{-\infty}}$ $)_{(i,j)\in \Ga_n})$, o\`{u} $\Ga_n={\N}^*\times \{2,\dots,n\}$, et $S_1=k(X^{p^{-\infty}})$ par convention.
 \sk

D'apr\`es le lemme ci-dessus, pour tout $n\geq 2$, on v\'{e}rifie aussit\^{o}t que  $K_n\subseteq S_n$.  De plus, pour tout $(i,n)\in \N^*\times (\N\setminus \{0,1\})$,  ${(Z_{1}^n)}^{p^{-1}}\in S_{n-1}({(Y_{1}^{n})}^{p^{-1}},\te_{i}^{n})$ par construction.

\begin{pro} {\label{pr39}} La famille $(Z_{i}^{j})_{(i,j)\in \N^* \times (\N\setminus \{0,1\})}$ est $p$-ind\'{e}pendente sur $K^p$.
\end{pro}

Tout le long de cette d\'emonstration, on va utiliser les notations suivantes :  pour tous $(i,j)\in \N^*\times (\N\setminus\{0,1\})$, on pose
$A_{i}^{j}=({(Z_{s_1}^{s_2})}^{p^{-1}})_{(s_1,s_2)<(i,j)}$ o\`u   $\leq$ est la relation d'ordre hexad\'{e}cimal, et pour tout $j\geq 2$,  $C_{i}^{j}=({(Z_{l}^{j})}^{p^{-i}}, {(Y_{l}^{j})}^{p^{-i}})_{l<i}$.

\pre D'apr\`{e}s le lemme d'\'echange, pour tout $j\geq 2$, pour tout $i\in {\N}^*$, il suffit de montrer que $Z_{i}^{j}\not\in K^p({(A_{i}^{j})}^p)$, ou encore ${(Z_{i}^{j})}^{p^{-1}}\not\in K(A_{i}^{j})$. Pour cela, on  suppose  l'existence d'un couple  $(i,j)\in {\N}^*\times (\N\setminus \{0,1\})$ tel que ${(Z_{i}^{j})}^{p^{-1}}\in K(A_{i}^{j})$. Comme $K=\di\bigcup_{n\geq 1} K_n$, il existe $n\in \N^*$ tel que ${(Z_{i}^{j})}^{p^{-1}}\in K_n(A_{i}^{j})$. Soit $n$ le plus petit entier qui v\'{e}rifie cette relation. Par construction, il est clair que  $n>1$. De plus, on distingue deux cas :
\sk

1-ier cas : si $n\leq j$, donc ${(Z_{i}^{j})}^{p^{-1}}\in K_n(A_{i}^{j})\subseteq  K_j(A_{i}^{j})\subseteq S_{j-1}(K_j)(C_{i}^{j})$. Or, on a, $\te_{i}^{j}={(Z_{i}^{j})}^{p^{-1}}\te_{2i}^{j-1}+\dots+{(Z_{1}^{j})}^{p^{-i}}{(\te_{2}^{j-1})}^{p^{-i+1}}+{(Y_{i}^{j})}^{p^{-1}}+\dots+
{(Y_{1}^{j})}^{p^{-i}}$, et comme pour tout $(r,l)\in \{1,\dots, i-1\}\times \N^*$, on a ${(Z_{r}^{j})}^{p^{-i}}, {(Y_{r}^{j})}^{p^{-i}}\in C_{i}^{j}$ et $\te_{l}^{j-1}\in S_{j-1}$, on en d\'{e}duit que ${(Y_{i}^{j})}^{p^{-1}}\in S_{j-1}(K_j)(C_{i}^{j})$. Cela conduit en vertu du th\'{e}or\`{e}me  {\ref{thm4}}, \`{a} $di(S_{j-1}(C_{i}^{j},{(Z_{i}^{j})}^{p^{-1}}, {(Y_{i}^{j})}^{p^{-1}}) /S_{j-1}(C_{i}^{j}))=2\leq di(S_{j-1}(K_j,C_{i}^{j})/S_{j-1}($ $C_{i}^{j}))=1$ ; et donc $2\leq 1$, c'est une contradiction.
\sk

2-i\`{e}me cas : $n\geq j+1$. On a ${(Z_{i}^{j})}^{p^{-1}}\not\in K_{n-1}(A_{i}^{j})$ et ${(Z_{i}^{j})}^{p^{-1}}\in K_n(A_{i}^{j})$, comme $K_n(A_{i}^{j})/$ $K_{n-1}(A_{i}^{j})$ est $q$-simple et $\te_{1}^{n}\not\in K_{n-1}(A_{i}^{j})$, (sinon on aura ${(Z_{1}^n)}^{p^{-1}}\in S_{n-1}({(Y_{1}^{n})}^{p^{-1}})$, absurde),  alors $K_{n-1}(A_{i}^{j})(\te_{1}^{n})=K_{n-1}(A_{i}^{j})({(Z_{i}^{j})}^{p^{-1}})$ $\subseteq S_{n-1}$. En particulier, $\te_{1}^{n}\in S_{n-1}$  et donc ${(Z_{1}^{n})}^{p^{-1}}\in S_{n-1}({(Y_{1}^{n})}^{p^{-1}})$, (car $\te_{1}^n={(Z_{1}^{n})}^{p^{-1}}\te_{2}^{n-1}+{(Y_{1}^{n})}^{p^{-1}}$), c'est une contradiction. \cqfd
\begin{lem} {\label{lem11}} Pour tout $(i,s) \in\N\times \N^*$, pour tout $j> i$, $o(\te_{s}^{j},K_i)=2^{j-i-1}s$.
\end{lem}
 \pre Il est imm\'ediat que $o(\te_{s}^1,K_0)=s$ pour tout $s\in \N^*$, donc on se ram\`ene au cas o\`u $j\geq 2$. Dans la suite, on va utiliser un raisonnement  par r\'{e}currence sur $j$. Par construction, on a $\te_{s}^{i+1}={(Z_{s}^{i+1})}^{p^{-1}}\te_{2s}^{i}+ \dots+{(Z_{1}^{i+1})}^{p^{-s}} {(\te_{2}^{i})}^{p^{-s+1}}+{(Y_{s}^{i+1})}^{p^{-1}}+\dots+{(Y_{1}^{s+1})}^{p^{-s}}$, et donc  ${(\te_{s}^{i+1})}^{p^s}\in K_i$. Si ${(\te_{s}^{i+1})}^{p^{s-1}}\in K_i$, on aura ${(Z_{1}^{i+1})}^{p^{-1}}\in K_i({(Y_{1}^{i+1})}^{p^{-1}})$, c'est une contradiction. D'o\`u,  ${(\te_{s}^{i+1})}^{p^{s-1}}\not\in K_i$, et donc $o(\te_{s}^{i+1},K_i)=s$.
 Soit maintenant $j$ un entier tel que $j\geq i+2$. Egalement, pour tout $s\geq 1$, on a $\te_{s}^{j}={(Z_{s}^{j})}^{p^{-1}}\te_{2s}^{j-1}+\dots+ {(Z_{1}^{j})}^{p^{-s}}{(\te_{2}^{j-1})}^{p^{-s+1}}+{(Y_{s}^{j})}^{p^{-1}}+\dots+{(Y_{1}^{j})}^{p^{-s}}$. Or, d'apr\`{e}s l'hypoth\`ese de r\'{e}currence, pour tout $s\in N^*$, pour tout $n\in \{0,\dots, s-1\}$,  on a $o({\te_{2(s-n)}^{j-1}},K_i)=2^{j-1-i-1}.2(s-n)=2^{j-i-1}(s-n)$, d'o\`u
 pour tout $n\in \{1,\dots, s-1\}$, $$o({(\te_{2(s-n)}^{j-1})}^{p^{-n}},K_i)=2^{j-i-1}(s-n)+n<2^{j-i-1}s=o({\te_{2(s)}^{j-1}},K_i).$$ Il en r\'esulte que
  $o(\te_{s}^j,K_i)=o(\te_{2s}^{j-1},K_i)=2^{j-i-1}s$. \cqfd
\begin{lem} {\label{lem12}}
Pour tout $i\in \N$, ${K_i}^{p^{-1}}\cap K=K_i(\te_{1}^{i+1})$.
\end{lem}
\pre Il est clair que ${K_i}(\te_{1}^{i+1})\subseteq K_{i}^{p^{-1}}\cap K$. S'il existe $\al\in K_{i}^{p^{-1}}\cap K$ tel que $\al\not\in K_i(\te_{1}^{i+1})$, comme $K_{i+1}/K_i$ est $q$-simple, on a $\al\not\in K_{i+1}$, sinon $K_i(\te_{1}^{i+1})=K_i(\al)$. Soit $j$ le plus grand entier tel que $\al\not\in K_j$, donc $\al\in K_{j+1}$ avec $i+1\leq j$. En outre, on aura $o(\al,K_{j-1})=1$, et par suite $K_j(\al)=K_j(\te_{1}^{j+1})\simeq K_j\otimes_{K_{j-1}} K_{j-1}(\al)$ ; d'o\`u $K_j(\te_{1}^{j+1})/K_{j-1}$ est modulaire. Or, $o(\te_{2}^j,K_{j-1})=2$, donc le syst\`{e}me $(1, {(\te_{2}^j)}^p)$ est $K_{j-1}$ libre, en particulier il est libre sur $K_{j-1}\cap {[K_j(\te_{1}^{j+1})]}^p$. Compl\'{e}tons ce syst\`{e}me en une base $B$ de ${[K_j(\te_{1}^{j+1})]}^p$ sur $K_{j-1}\cap {[K_j(\te_{1}^{j+1})]}^p$.
Comme $K_j(\te_{1}^{j+1})/K_{j-1}$ est modulaire, donc ${[K_j(\te_{1}^{j+1})]}^p$ et $K_{j-1}$ sont $K_{j-1}\cap {[K_j(\te_{1}^{j+1})]}^p$-lin\'{e}airement disjointes ; et par suite $B$ est aussi une base de $K_{j-1}({[K_j(\te_{1}^{j+1})]}^p)$ sur $K_{j-1}$. Or
 l'\'{e}quation de d\'{e}finition de $\te_{1}^{j+1}$ sur $K_{j-1}$ s'\'{e}crit : ${(\te_{1}^{j+1})}^p=Z_{1}^{j+1}{(\te_2^{j})}^p+Y_{1}^{j+1}$, par identification on aura $Z_{1}^{j+1}, Y_{1}^{j+1}\in K_{j-1}\cap {[K_j(\te_{1}^{j+1})]}^p\subseteq K^p$, c'est une contradiction avec la proposition {\ref{pr39}}; et donc ${K_{i}}^{p^{-1}}\cap K=K_i(\te_{1}^{i+1})$. \cqfd
\begin{cor} {\label{cor14}} Pour tous $(i,n)\in \N\times {\N}^*$, on a ${K_i}^{p^{-n}}\cap K\subseteq {K_{i+1}}^{p^{-n+1}}\cap K$. En particulier, ${K_i}^{p^{-n}}\cap K\subseteq K_{n+i}$ pour tout entier $n$.
\end{cor}
\pre Il suffit de remarquer que ${K_{i}}^{p^{-n}}\cap K={({K_i}^{p^{-1}}\cap K)}^{p^{-n+1}}\cap K\subseteq {K_{i+1}}^{p^{-n+1}}\cap K$. \cqfd \sk

Comme cons\'{e}quence, on a :
\begin{cor} {\label{cor15}} Pour tout $j\in \N$, $K/K_j$ est $lq$-finie.
\end{cor}
\pre Application imm\'{e}diate du r\'{e}sultat pr\'{e}c\'{e}dente, et du fait que $di(K_{n+j}$ $/K_j)=n$ pour tous $j,n\in \N$. \cqfd \sk

Soit  $n$ un entier naturel non nul. Consid\'{e}rons la suite $(n_s)_{s\in \N^*}$ des entiers d\'{e}finie par la relation de r\'ecurrence suivante :
$n_1=n$, et pour tout $s\in N^*\setminus \{1\}$, $n_{s+1}=\di[\frac{n_s}{2}]$, o\`{u} $\di[.]$ d\'{e}signe la partie enti\`{e}re. La suite  des entiers ainsi obtenue est d\'ecroissante, et donc stationnaire. Le premier entier $e$ tel que $n_e=1$ s'appelle longueur de parit\'e inf\'erieur de $n$ et se note $lpi(n)$. On v\'{e}rifie imm\'{e}diatement que :
\sk

\begin{itemize}{\it
\item Pour tout $s> lpi(n)$, $n_s=0$.
\item Si $n\leq m$, alors $lpi(n)\leq lpi(m)$.
\item Pour tout $i\in \N^*$, $2n_{i+1}+\ep_i=n_i$, avec $\ep_i=0$ ou $1$ selon la parit\'{e} de $n_i$.}
\end{itemize}
\sk

En particulier, $2n_{i+1}\leq n_{i}$
et $n_i+1\leq 2(n_{i+1}+1)$, et par suite pour tout $j\geq i$, on en d\'eduit que $2^{j-i}n_j\leq n_i$ et $n_i+1\leq 2^{j-i}(n_j+1)$. Notamment, si l'on pose $p_1=lpi(n)$, on aura $2^{p_1-1}\leq n< 2^{p_1}$. Autrement dit, $p_1$ est le plus petit entier tel que $n< 2^{p_1}$.
\begin{thm} {\label{thm24}} Pour tout $(i, n)\in \N\times \N^*$, ${K_i}^{p^{-n}}\cap K\subseteq {K_{i+1}}^{p^{- \di[\frac{n}{2}]}}$.
\end{thm}

Pour la preuve de ce th\'eor\`eme, on aura besoin des r\'esultats suivants :
\begin{pro} {\label{pr40}} Soit $n$ un entier naturel non nul.  Pour tout $s\in \N^*$, pour tout $n\in \N$, on a $o(\te_{n_s}^{s+r}, K_r(\te_{n_1}^{1+r},$ $\dots,$ $\te_{n_{s-1}}^{s-1+r}))=n_s=\di\sup_{i\in \N^*}(o(\te_{n_i}^{i+r}, K_r(\te_{n_1}^{1+r},$ $\dots,$ $\te_{n_{s-1}}^{s-1+r}))$. En paticulier, $\{\te_{n_1}^{1+r},\dots,\te_{n_{p_1}}^{p_1+r}\}$ est une $r$-base canoniquement ordonn\'ee de $K_r((\te_{n_i}^{i})_{i\in \N^*})/K_r$, o\`u $p_1=lpi(n)$.
\end{pro}
\pre Comme les $K_i/k$, ($1\leq i$) sont construites de la m\^eme fa\c{c}on, on se contente de pr\'esenter la d\'emonstration dans le cas o\`u $r=0$. Par convention, pour tout $s>p_1$, $o(\te_{n_s}^s,k)=o(\te_{0}^s,k)=0$, donc on se ram\`ene au cas o\`u $s\leq p_1$.
D'apr\`es le lemme {\ref{lem11}}, pour tout $i\in \{1,\dots, p_1\}$, $o(\te_{n_i}^{i},k)=o(\te_{n_i}^{i},K_0)=2^{i-1}n_i\leq n_1=o(\te_{n_1}^{1},k)$, et donc $o(\te_{n_1}^{1},k)=\di\sup_{i\in \N^*} (o(\te_{n_i}^i,k))$. Soit $s\in \{2,\dots, p_1-1\}$, supposons en tenant compte de l'hypoth\`ese de r\'ecurrence que pour tout $j\in \{1,\dots,s\}$, on a $o(\te_{n_j}^j, k(\te_{n_1}^1,\dots,\te_{n_{j-1}}^{j-1}))=n_j=\di\sup_{i\in \N^*}(o(\te_{n_i}^i, k(\te_{n_1}^1,$ $\dots,\te_{n_{j-1}}^{j-1}))$. Pour simplifier l'\'ecriture, on pose $L_1=k(\te_{n_1}^{1},\dots,\te_{n_s}^s)$, $L_2=k(\te_{n_1}^{1},\dots,\te_{n_{p_1}}^{p_1})$ et $L_3=L_2(K_s)= K_s(\te_{n_{s+1}}^{s+1},\dots,\te_{n_{p_1}}^{p_1})$. D'une part, $L_3/k$ \'etant $q$-finie, donc en vertu du lemme {\ref{lem11}}, pour tout $i\geq s+1$, $o(\te_{n_i}^{i},K_s)=2^{i-s-1}n_i\leq n_{s+1}=o(\te_{n_{s+1}}^{s+1}, K_s)$. D'autre part, d'apr\`es la proposition {\ref{pr18}} et {\ref{pr19}}, $n_{s+1}=o(\te_{n_{s+1}}^{s+1}, K_s$ $) $ $\leq o(\te_{n_{s+1}}^{s+1}, L_1)\leq o_1(L_2/L_1)=o_{s+1}(L_2/k)\leq o_{s+1}(L_3/k)=o(\te_{n_{s+1}}^{s+1}, K_s)=n_{s+1}$, et donc $o(\te_{n_{s+1}}^{s+1}, L_1)=o_{s+1}(L_2/k)=o_1(L_2/$ $L_1)=n_{s+1}$. Par suite,  pour tout $s\in N^*$, $o(\te_{n_s}^s, k(\te_{n_1}^1,\dots,\te_{n_{s-1}}^{s-1}))$ $=n_s=\di\sup_{i\in \N^*}(o(\te_{n_i}^i, k(\te_{n_1}^1,\dots,\te_{n_{s-1}}^{s-1}))$. En particulier, par application de l'algorithme de la compl\'etion des $r$-bases, on aura $\{\te_{n_1}^1,\dots,\te_{n_{p_1}}^{p_1}\}$ est une $r$-base canoniquement ordonn\'ee de $k((\te_{n_i}^{i})_{i\in \N^*})/k$. \cqfd \sk

Soit  $n$ un entier naturel non nul. Consid\'{e}rons la suite $(n'_s)_{s\in \N^*}$ des entiers d\'{e}finie par la relation de r\'ecurrence suivante :
\sk

\begin{itemize}{\it
\item $n'_1=n$ si $n$ est paire, et $n'_1=n+1$ dans le cas contraire.
\item Pour tout $s\in \N^*$, $n'_{s+1}=\di\frac{n'_s}{2}$ si $n'_s$ est paire, et $n'_{s+1}=\di\frac{n'_s+1}{2}$ si $n'_s$ est impaire.}
\end{itemize}
\sk

La suite des entiers ainsi obtenue est d\'{e}croissante, donc stationnaire. Le plus petit entier $r$ tel que  $n'_r=1$ s'appelle longueur de parit\'e sup\'erieur de $n$ et se note $lps(n)$. On v\'{e}rifie imm\'{e}diatement que :
\sk

\begin{itemize}{\it
\item[$\bullet$] Pour tout $s\geq lps(n)$, $n'_s=1$.
\item[$\bullet$] Pour tout $i\in \N^*$, $2n'_{i+1}=n'_i+\ep_i$, avec $\ep_i=0$ ou $1$ selon la parit\'{e} de $n'_i$.}
\end{itemize}
\sk

En particulier, $n'_i\leq 2n'_{i+1}$ et $2(n'_{i+1}-1)\leq n'_i-1$. Par suite, pour tout $j\geq i$,  $n'_i\leq 2^{j-i}n'_j$ et $2^{j-i}(n'_j-1)\leq n'_i-1$. En outre, si on pose $p_2=lps(n)$, on aura $2^{p_2-2}< n \leq 2^{p_2-1}$. Autrement dit, $p_2-1$ est le plus petit entier  tel que $n \leq 2^{p_2-1}$.
De m\^eme, il est ais\'e d'obtenir l'expression $n_i\leq n'_i$,  pour tout $i\in \N^*$.  Toutefois, on montre par r\'{e}currence que $n'_i\leq n_i+1$. En effet, pour $i=1$ le r\'{e}sultat est v\'{e}rifi\'e, et si $n'_i\leq n_i+1$, ou encore $n'_i+1\leq 2(n_{i+1}+1)+1$, on aura $n'_{i+1}\leq \di\frac{n'_i+1}{2}\leq (n_{i+1}+1)+\di\frac{1}{2}$ ; et par suite $n'_{i+1}\leq [n_{i+1}+1+\di\frac{1}{2}]=n_{i+1}+1$.
En outre, $n'_i=n_i+\al_i$, avec $\al_i=0$ ou $1$ ; et par cons\'{e}quent $lps(n)=lpi(n)+\al$, avec $\al=0$ ou $1$. Egalement, on a $2^{p_2-2}\leq 2^{p_1-1}\leq n\leq 2^{p_2-1}\leq 2^{p_1}$, o\`u $p_1=lpi(n)$ et $p_2=lps(n)$, et comme cons\'equences, $p_2=p_1$ si et seulement si $n$ est puissance de $2$.
\sk

On pose ensuite,  $\Delta=\{(f,g)\in \N^2|\, \, f\geq 2$ et $1\leq g\leq n'_f\}$, et $Q_{0}^{n}=Q(\te_{n'_1}^{1}, ({(Z_{g}^{f})}^{p^{-n'_2}},$ $ {(Y_{g}^{f})}^{p^{-n'_2}}$ $)_{(f,g)\in \De})$.
\sk

Comme pour tout $i\in \N^*\setminus \{1\}$, et pour tout $j\in \{2,\dots, i\}$, $2^{i-j}n_i\leq n_j\leq n'_j\leq n'_2$, on en d\'{e}duit que $\{{(Z_{l}^{j})}^{p^{-2^{i-j}n_i}},{(Y_{l}^{j})}^{p^{-2^{i-j}n_i}}\}\subseteq Q_{0}^{n}$. D'apr\`es le lemme {\ref{lem10}}, pour tout $i\in \N^*$, on a $\te_{n_i}^{i}\in Q_{0}^{n}$ ; et donc $Q((\te_{n_i}^{i})_{i\in \N^*})\subseteq Q_{0}^{n}$.  Pour tout $i\in \N$, posons de m\^{e}me $\f_{i}^{n}=Q_{0}^{n}(K_i)$ et $\f^{n}=Q_{0}^{n}(K)$. On v\'{e}rifie aussit\^{o}t que:
\sk

\begin{itemize}{\it
\item[\rm{(1)}]  Pour tout $i\in\N$, $\f_{i}^{n}\subseteq \f_{i+1}^{n}$, et $\f^{n}=\di\bigcup_{i\in\N} \f_{i}^{n}$.\\
\item[\rm{(2)}] Pour tout $i\in\N$, $\f_{i+1}^{n}/\f_{i}^{n}$ est $q$-simple.\\
\item[\rm{(3)}] Pour tout $1\leq i$, ${({\te_{n'_i-1}^{i}})}^{p^{-1}}\in Q_{0}^{n}$.}
\end{itemize}
\sk

En effet,  $\te_{n'_i-1}^{i}\in Q(\te^{1}_{2^{i-1}(n'_i-1)}, ({(Z_{l}^{j})}^{p^{-2^{i-j}(n'_i-1)}},$ ${(Y_{l}^{j} )}^{p^{-2^{i-j}(n'_i-1)}}$ $)_{(j,l)\in \La_{n'_i-1}^{i}})$, en vertu du lemme {\ref{lem10}}.
Comme $\te_{n'_i-1}^{1}\in  {(Q_{0}^{n})}^p$ ($Q$ \'etant donn\'e parfait), et pour tout $j\in \{2,\dots, i\}$, $2^{i-j}(n'_i-1)\leq n'_{j}-1\leq n'_2-1$, (et donc en particulier $2^{i-1}(n'_i-1)\leq n'_1-1$),
alors $\te_{n'_i-1}^{i}\in Q(\te^{1}_{n'_1-1},({(Z_{l}^{j})}^{p^{-(n'_2-1)}},$ ${(Y_{l}^{j})}^{p^{-(n'_2-1)}}$ $)_{(j,l)\in
\La_{n'_i-1}^{i}})\subseteq {(Q_{0}^{n})}^p$ ;
ou encore ${(\te_{n'_i-1}^{i})}^{p^{-1}}\in {Q_{0}^{n}}$. En outre ${(\te_{n'_i-1}^{i})}^{p^{-1}}\in {\f_{0}^{n}}$, et
pour tout $i\geq 1$, $\te_{n'_i}^{i}\in \f_{i-1}^{n}$, (car $\te_{n'_i}^{i}= {(Z_{n'_i}^{i})}^{p^{-1}}\te_{2n'_i}^{i-1}+{(\te_{n'_i-1}^{i})}^{p^{-1}}+{(Y_{n'_i}^i)}^{p^{-1}}$, avec ${(Z_{n'_i}^{i})}^{p^{-1}}, {(Y_{n'_i}^i)}^{p^{-1}}\in \f_{0}^{n}$ et $\te_{2n'_i}^{i-1}\in K_{i-1}\subseteq \f_{i-1}^{n}$.
\begin{lem} {\label{lem13}} Pour tout $(s,j)\in \N^*\times \N^*$, $o(\te_{n'_s+j}^{s}, \f_{s-1}^{n})=j$.
\end{lem}

Pour simplifier l'\'ecriture, on pose
$L_2= k({(Q_{0}^{n})}^{p^{-\infty}}, {(Y_{n'_2+1}^{2})}^{p^{-\infty}})=k({X}^{p^{-\infty}},$ $({(Z_{l}^{s})}^{p^{-\infty}}, {(Y_{l}^{s})}^{p^{-\infty}})_{(s,l)\in \Delta},$ ${(Y_{n'_2+1}^{2})}^{p^{-\infty}} )$,
et pour tout $t> 2$, $L_t=S_{t-1}({(Q_{0}^{n})}^{p^{-\infty}},$ $ {(Y_{n'_t+1}^{2})}^{p^{-\infty}})=k({X}^{p^{-\infty}}, ({(Z_{i}^h)}^{p^{-\infty}}, {(Y_{i}^{h})}^{p^{-\infty}})_{(i,h)\in \Ga_{t-1}},$ $ ({(Z_{l}^{j})}^{p^{-\infty}}, {(Y_{l}^{j})}^{p^{-\infty}})_{(j,l)\in \Delta}, $ ${(Y_{n'_t+1}^{t})}^{p^{-\infty}})$, o\`{u} $\Ga_{t-1}={\N}^*\times \{2,\dots,t-1\}$. Il est clair que $\f_{s}^{n}\subseteq L_s$, par ailleurs, comme $Q$ est parfait, et  $\te_{i}^{s}=
{(Z_{i}^{s})}^{p^{-1}}\te_{2i}^{s-1}+\dots+ {(Z_{1}^{s})}^{p^{-i}}{(\te_{2}^{s-1})}^{p^{-i+1}}+{(Y_{i}^{s})}^{p^{-1}}+\dots +{(Y_{1}^{s})}^{p^{-i}}$, $s> 1$ et $i\in \N^*$, alors pour tout $j\in \N$,  ${(\te_{i}^{s})}^{p^{-j}}\in k({X}^{p^{-\infty}}, ({(Z_{r}^h)}^{p^{-\infty}}, {(Y_{r}^{h})}^{p^{-\infty}})_{(r,h)\in \Ga_{s-1}},$ $ {({(Z_{l}^{s})}^{p^{-\infty}}, {(Y_{l}^{s})}^{p^{-\infty}})}_{1\leq l\leq i})$. En particulier, ${(\te_{n'_s}^{s})}^{p^{-j}}\in L_s$.

\pre Pour $s=1$, le r\'esultat est trivial, donc on est amen\'e au cas  $s\geq 2$. Par construction, $\te_{n'_s+1}^{s}={(Z_{n'_s+1}^s)}^{p^{-1}}\te_{2(n'_s+1)}^{s-1}+{(\te_{n'_s}^{s})}^{p^{-1}}+{(Y_{n'_s+1}^s)}^{p^{-1}}=
{(Z_{n'_s+1}^{s})}^{p^{-1}}\te_{2(n'_s+1)}^{s-1}$ $+\dots+ {(Z_{1}^{s})}^{p^{-n'_s-1}}{(\te_{2}^{s-1})}^{p^{-n'_s}}+{(Y_{n'_s+1}^{s})}^{p^{-1}}+\dots +{(Y_{1}^{s})}^{p^{-n'_s-1}}$, donc ${(\te_{n'_s+1}^{s})}^p=Z_{n'_s+1}^{s}{(\te_{2(n'_s+1)}^{s-1})}^p+\te_{n'_s}^{s}+Y_{n'_s+1}^{s}$, et par suite
${(\te_{n'_s+1}^{s})}^p\in \f_{s-1}^{n}$. Toutefois, si $\te_{n'_s+1}^{s}\in \f_{s-1}^{n}$, on aura ${(Z_{n'_s+1}^{s})}^{p^{-1}}\in L_s$
Remarquons que pour tout $h\leq s-1$ et pour tout $i\in {\N}^*$, $Z_{i}^h\not= Z_{n'_s+1}^s$, \'{e}galement pour tout $1\leq l\leq n'_j$. $Z_{n'_s+1}^{s}\not= Z_{s}^j$ ; et comme la famille $(Z_{i}^{j}, Y_{i}^{j})_{(i,j)\in {{\N}^{*}}^2}$ est $p$-ind\'{e}pendante sur $k^p$, alors ${(Z_{n'_s+1}^{s})}^{p^{-1}}\not \in L_s$, c'est une contradiction.  D'o\`{u} $o(\te_{n'_s+1}^{s}, \f_{s-1}^{n})=1$. Soit $j$ un entier naturel diff\'erent de $0$ et $1$, supposons par application de la propri\'et\'e de  r\'ecurrence que $o(\te_{n'_s+i}^s,\f_{s-1}^{n})=i$ pour tout $i\in \{1,\dots,j\}$. Comme  $\te_{n'_s+j+1}^{s}={(Z_{n'_s+j+1}^{s})}^{p^{-1}}\te_{2(n'_s+j+1)}^{s-1}+{(\te_{n'_s+j}^{s})}^{p^{-1}}+{(Y_{n'_s+j+1}^{s})}^{p^{-1}}$,  on aura ${(\te_{n'_s+j}^{s})}^{p^{-1}}\in \f_{s-1}^{n}(\te_{n'_s+j+1}^s, {(Z_{n'_s+j+1}^{s})}^{p^{-1}}, $ ${(Y_{n'_s+j+1}^{s})}^{p^{-1}})$. D'o\`u, $j+1=o({(\te_{n'_s+j}^s)}^{p^{-1}},\f_{s-1}^{n})\leq o_1(\f_{s-1}^{n}(\te_{n'_s+j+1}^s, {(Z_{n'_s+j+1}^{s})}^{p^{-1}}, {(Y_{n'_s+j+1}^{s})}^{p^{-1}})$ $/$ $\f_{s-1}^{n})\leq o(\te_{n'_s+j+1}^{s},\f_{s-1}^{n})\leq j+1$,  et par suite $o(\te_{n'_s+j+1}^s,\f_{s-1}^{n})=j+1$. \cqfd \sk

Comme cons\'equence, on a :
\begin{cor} {\label{cor16}}
Pour tout $s\in \N^*$, ${(\te_{2(n'_{s+1}+1)}^s)}^p\not\in \f_{s-1}^{n}(\te_{2n'_{s+1}}^{s})$. En particulier, ${(\te_{n'_{s+1}+1}^{s+1})}^p\not\in \f_{s-1}^{n}($ $\te_{2n'_{s+1}}^{s})$.
\end{cor}
\pre Si ${(\te_{2(n'_{s+1}+1)}^s)}^p\in \f_{s-1}^{n}(\te_{2n'_{s+1}}^{s})$, alors $o($ ${(\te_{2(n'_{s+1}+1)}^s)}^p,$ $ \f_{s-1}^{n})=o({(\te_{n'_{s}+\ep_s+2}^s)}^p,$ $\f_{s-1}^{n})$ $=\ep_s+2-1\leq o(\te_{2n'_{s+1}}^s, \f_{s-1}^{n})=o(\te_{n'_{s}+\ep_s}^s, \f_{s-1}^{n})=\ep_s$, ce qui implique que $1<0$, contradiction. D'o\`{u}, ${(\te_{2(n'_{s+1}+1)}^s)}^p\not\in \f_{s-1}^{n}(\te_{2n'_{s+1}}^{s})$. D'autre part, comme ${(\te_{n'_{s+1}+1}^{s+1})}^{p}=Z_{n'_{s+1}+1}^{s+1}{(\te_{2(n'_{s+1}+1)}^{s})}^p+{(Z_{n'_{s+1}}^{s+1})}^{p^{-1}}\te_{2n'_{s+1}}^{s}+
{(\te_{n'_{s+1}-1}^{s+1})}^{p^{-1}}+Y_{n'_{s+1}+1}^{s+1}+ {(Y_{n'_{s+1}}^{s+1})}^{p^{-1}}$, et ${(\te_{n'_{s+1}-1}^{s+1})}^{p^{-1}}, {(Z_{n'_{s+1}}^{s+1})}^{p^{-1}}, {(Y_{n'_{s+1}}^{s+1})}^{p^{-1}}$ sont des \'el\'ements de $\f_{0}^{n}$, on en d\'eduit que ${(\te_{n'_{s+1}+1}^{s+1})}^p\not\in \f_{s-1}^{n}(\te_{2n'_{s+1}}^{s})$.   \cqfd
\begin{lem}{\label{lema14}} Pour tout $s\geq 1$, ${(Z_{n'_{s+1}+1}^{s+1})}^{p^{-1}} \not\in \f^n$.
\end{lem}
\pre
En posant $m_0=Q(({(Z_{g}^{f})}^{p^{-\infty}},$ $ {(Y_{g}^{f})}^{p^{-\infty}}$ $)_{(f,g)\in \De})$, $m=k(m_0)$ et $M_i=m(K_{i})$, $M=m(K)=\di\bigcup_{i\geq 1} M_i$, et pour tout $i\in \N^*$, ${Z'}_{i}^{j}=Z_{n'_j+i}^{j}$,  ${Y'}_{i}^{j}=Y_{n'_j+i}^{j}$ et $\al_{i}^j=\te_{n'_j+i}^j$ ;  on se retrouve dans des conditions analogues \`a celles de la proposition {\ref{pr39}}, o\`u $X$, ${Z'}_{i}^{j}$, ${Y'}_{i}^{j}$ et $\al_{i}^j$ jouent  des  r\^oles similaires que $X$, ${Z}_{i}^{j}$, ${Y}_{i}^{j}$ et $\te_{i}^j$ auquel cas le r\'esultat d\'ecoule imm\'ediatement.  \cqfd \sk

Toutefois, on v\'erifie aussit\^ot  que pour tout $i\in \N^*$,  $\f_{i}^{n}=\f_{i-1}^{n}((\al_{j}^{i})_{j\geq 1})$, et pour tout $j\geq 1$,  $o(\al_{j}^{i},\f_{i-1}^{n})=j$.
\begin{lem} {\label{lem14}}
${(\f_{0}^{n})}^{p^{-1}}\cap {\f}^n=\f_{0}^{n}(\al_{1}^{1})$.
\end{lem}
 \pre Il est clair que $\f_{0}^{n}(\al_{1}^{1})\subseteq {(\f_{0}^{n})}^{p^{-1}}\cap \f^n$. Supposons que ${(\f_{0}^{n})}^{p^{-1}}\cap \f^n\not= \f_{0}^{n}(\al_{1}^{1})=\f_{0}^{n}(\te_{n'_1+1}^{1})$, donc il existe $\te \in {(\f_{0}^{n})}^{p^{-1}}\cap \f^n$ tel que $\te \not\in \f_{0}^{n}(\te_{n'_1+1}^{1})$. Soit $s$ le plus grand entier tel que $\te \not\in \f_{s}^{n}$, donc $\te\in \f_{s+1}^{n}$ et $1\leq s$. Comme $\f_{s+1}^{n}/\f_{s}^{n}$ est $q$-simple et $o(\te, \f_{s}^{n})=o(\te, \f_{s-1}^{n})=1$, et compte tenu de la transitivit\'e de la lin\'earit\'e disjointe, on aura  $\f_{s}^{n}(\te)=\f_{s}^{n}(\al_{1}^{s+1})\simeq \f_{s}^{n}\otimes_{\f_{s-1}^{n}} \f_{s-1}^{n}(\te)\simeq \f_{s}^{n}\otimes_{\f_{s-1}^{n}(\te_{2n'_{s+1}}^{s})} \f_{s-1}^{n}(\te_{2n'_{s+1}}^s)(\te)$. D'o\`u $\f_{s}^{n}(\te)/\f_{s-1}^{n}(\te_{2n'_{s+1}}^s)$ est modulaire, en outre ${(\f_{s}^{n}(\te))}^p$ et $\f_{s-1}^{n}(\te_{2n'_{s+1}}^s)$ sont ${(\f_{s}^{n}(\te))}^p\cap \f_{s-1}^{n}(\te_{2n'_{s+1}}^s)$-lin\'{e}airement disjointes. Comme ${(\te_{2(n'_{s+1}+1)}^{s})}^p\not\in \f_{s-1}^{n}(\te_{2n'_{s+1}}^{s})$ d'apr\`es le corollaire  {\ref{cor16}},  alors $(1,{(\te_{2(n'_{s+1}+1)}^{s})}^p)$ est libre sur $\f_{s-1}^{n}(\te_{2n'_{s+1}}^s)$, d'o\`u $(1,{(\te_{2(n'_{s+1}+1)}^{s})}^p)$ est libre sur ${(\f_{s}^{n}(\te))}^p\cap \f_{s-1}^{n}(\te_{2n'_{s+1}}^s)$. Compl\'{e}tons ce syst\`{e}me en une base $B$ de ${(\f_{s}^{n}(\te))}^p$ sur ${(\f_{s}^{n}(\te))}^p\cap \f_{s-1}^{n}(\te_{2n'_{s+1}}^s)$. En vertu de la lin\'{e}arit\'{e} disjointe, $B$ est aussi une base de $( \f_{s-1}^{n}(\te_{2n'_{s+1}}^s))({(\f_{s}^{n}(\te))}^p)$ sur $\f_{s-1}^{n}(\te_{2n'_{s+1}}^s)$. Or, par construction, on a ${(\al_{1}^{s+1})}^p={(\te_{n'_{s+1}+1}^{s+1})}^p=Z_{n'_{s+1}+1}^{s+1}{(\te_{2(n'_{s+1}+1)}^{s})}^p+{(Z_{n'_{s+1}}^{s+1})}^
 {p^{-1}}\te_{2n'_{s+1}}^{s}+{(\te_{n'_{s+1}-1}^{s+1})}^{p^{-1}}+ Y_{n'_{s+1}+1}^{s+1}+{(Y_{n'_{s+1}}^{s+1})}^{p^{-1}}$, et ${(Z_{n'_{s+1}}^{s+1})}^
 {p^{-1}}\te_{2n'_{s+1}}^{s}+{(\te_{n'_{s+1}-1}^{s+1})}^{p^{-1}}+ Y_{n'_{s+1}+1}^{s+1}+{(Y_{n'_{s+1}}^{s+1})}^{p^{-1}}\in \f_{s-1}^{n}(\te_{2n'_{s+1}}^s)$, (car ${(Z_{n'_{s+1}}^{s+1})}^{p^{-1}}, {(Y_{n'_{s+1}}^{s+1})}^{p^{-1}},$ $ {(\te_{n'_{s+1}-1 }^{s+1})}^{p^{-1}}\in \f_{0}^{n}$). Par identification, on obtient ${(Z_{n'_{s+1}+1}^{s+1})}\in {(\f_{s}^{n}(\te))}^p\cap \f_{s-1}^{n}(\te_{2n'_{s+1}}^s)$, et donc
 ${(Z_{n'_{s+1}+1}^{s+1})}^{p^{-1}}\in {\f_{s}^{n}(\te)} \subseteq \f^n$,  c'est une contradiction d'apr\`es le lemme {\ref{lema14}} ci-dessus ; et
 par suite ${(\f_{0}^{n})}^{p^{-1}}\cap {\f}^n=\f_{0}^{n}(\al_{1}^{1})$. \cqfd
\begin{lem} {\label{lem15}} Pour tout $n\in \N^*\setminus \{1\}$, $o_2({(k(\te_{n_1}^{1},\dots,\te_{n_{p_1}}^{p_1}))}^{p^{-1}}\cap K/k)=\di[\frac{n_1+1}{2}]$, o\`u $p_1=lpi(n)$ et $n_1=n$, et  pour tout $s\in \N^*$, $n_{s+1}=\di[\frac{n_s}{2}]$.
\end{lem}
\pre Notons $L=k(\te_{n_1}^1,\dots,\te_{n_{p_1}}^{p_1})$, $L_1=L^{p^{-1}}\cap K$, et $L_2=k(\te_{n_1+1}^{1},\te_{n_2+1}^{2})$. On distingue deux cas.
\sk

1-ier cas : si $n$ est paire, alors $n'_1=n$ et $n'_2=n_2=\di[\frac{n}{2}]$.
 En vertu du lemme {\ref{lem14}}, $L_1\subseteq {(\f_{0}^{n})}^{p^{-1}}\cap K\subseteq {(\f_{0}^{n})}^{p^{-1}}\cap {\f}^n=\f_{0}^{n}(\al_{1}^{1})$. Or, $\f_{0}^{n}(\al_{1}^{1})\subseteq K_1(({(Y_{i}^{j})}^{p^{-n_2}}, {(Z_{i}^{j})}^{p^{-n_2}})_{(i,j)\in {{\N}^*}\times (\N\setminus \{0,1\})})$, donc pour tout $\be\in \f_{0}^{n}(\al_{1}^{1})$, on a $\be^{p^{n_2}}\in K_1$ ; et par suite $L^{p^{-1}}\cap K\subseteq \f_{0}^{n}(\al_{1}^{1})\cap K\subseteq  {K_1}^{p^{-n_2}}\cap K$. Compte tenu de la proposition {\ref{pr18}} et {\ref{pr19}}, on a $n_2=o_2(L/k)\leq o_2(L_1/k)\leq o_2({K_1}^{p^{-n_2}}\cap K/k)=o_1({K_1}^{p^{-n_2}}\cap K/K_1)=n_2$, et donc $o_2(L_1/k)=n_2=\di[\frac{n+1}{2}]$.
 \sk

 2-i\`eme cas : si $n$ est impaire, comme $k({L_1}^p)\subseteq L\subseteq L_1$, d'apr\`es la proposition {\ref{pr19}}, $o_i(L_1/k)=o_i(L/k)+\ep_i$, avec $\ep_i=0$ ou $1$. D'une part, par construction, on a ${(\te_{n_2+1}^{2})}^p=Z_{n_2+1}^{2}{(\te_{2(n_2+1)}^{1})}^p+{\te_{n_2}^{2}}+{Y_{n_2+1}^{2}}$. D'autre part, on a $2n_2+1=n=n_1$, et pour tout $i\in \N^*$, $Q(\te_{i}^1)=Q(X^{p^{-i}})$, donc en particulier $Q(\te_{2n_2+1}^{1})=Q({(\te_{2(n_2+1)}^{1})}^p)$ ; et par cons\'equent, ${(\te_{n_2+1}^{2})}^p\in L$. On en d\'eduit que $L_2\subseteq L_1$, et d'apr\`es la proposition {\ref{pr40}} et la proposition {\ref{pr18}}, on a   $o_2(L_2/k)=n_2+1\leq o_2(L_1/k)\leq n_2+1$. D'o\`u, $o_2(L_1/k)=n_2+1=\di[\frac{n_1+1}{2}]$. \cqfd
 \begin{lem} {\label{lem16}} Soient $n$ et $r$ deux entiers naturels avec $n$ non nul. Si ${K_r}^{p^{-n_1}}\cap K\subseteq {K_{r+1}}^{p^{-n_2}}\cap K\subseteq \dots\subseteq {K_{r+p_1-1}}^{p^{-n_{p_1}}}\cap K \subseteq K_{r+p_1}$, o\`u $p_1=lpi(n)$, $n_1=n$, et pour tout $s\in \N^*$, $n_{s+1}=\di[\frac{n_s}{2}]$, alors on aura  ${K_r}^{p^{-n}}\cap K=K_r(\te_{n_{r+1}}^1,\dots,\te_{n_{r+p_1}}^{p_1})$.
 \end{lem}
 \pre
 Comme les $K_i/k$, $1\leq i$ sont construites de fa\c{c}on analogue, on se ram\'ene au cas o\`u $r=0$. Pour tout $n\in \N^*$,  $k_n$ d\'esigne toujours $k^{p^{-n}}\cap K$. D'une part,  d'apr\`{e}s la proposition {\ref{pr40}}, on a $o_1(k(\te_{n_1}^{1},\dots, \te_{n_{p_1}}^{p_1})/k)=n_1$, et donc $k(\te_{n_1}^{1},\dots, \te_{n_{p_1}}^{p_1})\subseteq k_{n_1}$. D'autre part, comme $k_{n_1}\subseteq {K_{s-1}}^{p^{-n_s}}\cap K$, en vertu de la proposition {\ref{pr40}}, on aura $n_s=o_s(k(\te_{n_1}^{1},\dots, \te_{n_{p_1}}^{p_1})/k)\leq o_s(k_{n_1}/k)\leq o_s({K_{s-1}}^{p^{-n_s}}\cap K/k)=o_1({K_{s-1}}^{p^{-n_s}}\cap K/K_{s-1})=n_s$, (car ${K_{s-1}}^{p^{-n_s}}\cap K/k$ est $q$-finie et $K_{s-1}/k$ est relativement parfaite de degr\'e d'irrationalit\'{e} $s-1$) ; et donc $o_s(k_{n_1}/k)=n_s$. D'o\`{u}, $k_{n_1}= k(\te_{n_1}^1,\dots,\te_{n_{p_1}}^{p_1})$. \cqfd
 \sk

{\bf Preuve du th\'eor\`eme {\ref{thm24}}.}
D'abord, il est \`a signaler qu'on va continuer \`a  utiliser les notations du lemme ci-dessus. Soit $n$ un entier naturel non nul. D'apr\`{e}s le corollaire {\ref{cor14}}, pour tout $i\in \N$, ${K_i}^{p^{-1}}\cap K\subseteq {K_{i+1}}$ et ${K_i}^{p^{-2}}\cap K\subseteq {K_{i+1}}^{p^{-1}}\cap K$ ; donc le th\'eor\`eme est v\'{e}rifi\'e si $n\in\{1, 2\}$.  Supposons  que le th\'eor\`eme est satisfait pour tout $j\leq n$. Autrement dit, pour tout $j\leq n$, pour tout $i\in \N$, ${K_i}^{p^{-j}}\cap K\subseteq {K_{i+1}}^{p^{-\di[\frac{j}{2}]}}$, et montrons que ${K_i}^{p^{-n-1}}\cap K\subseteq {K_{i+1}}^{p^{-\di[\frac{n+1}{2}]}}$ pour tout $i\in \N$. Dans un premier temps, on se limite au cas  $i=0$.
D\'esormais et sauf mention du contraire, on notera $e=\di[\frac{n_1+1}{2}]$.

{\bf Cas particulier :}  Compte tenu de l'hypoth\`{e}se de r\'{e}currence, on a ${k}^{p^{-n_1}}\cap K\subseteq {K_{1}}^{p^{-n_2}}\cap K\subseteq \dots\subseteq {K_{p_1-1}}^{p^{-n_{p_1}}}\cap K\subseteq K_{p_1}$, o\`u $p_1=lpi(n)$. En vertu du lemme {\ref{lem16}}, $k_{n_1}= k(\te_{n_1}^1,\dots,\te_{n_{p_1}}^{p_1})$, et d'apr\`es le lemme {\ref{lem15}}, $o_2({k_{n_1}}^{p^{-1}}\cap K/k)=o_2(k_{n_1+1}/k)=e$.  Comme $o(\te_{n_1+1}^{1}, k)=n_1+1=o_1(k_{n_1+1}/k)$,  donc  $\te_{n_1+1}^{1}\in k_{n_1+1}$ ; et par application de l'algorithme de la compl\'etion des $r$-bases, on compl\`ete $\te_{n_1+1}^{1}$ en une $r$-base de $k_{n_1+1}/k$. D'o\`u, pour tout $\be\in k_{n_1+1}$,  on a $\be^{p^{e}}\in k(\te_{n_1+1}^1)\subseteq K_1$. Il en r\'esulte que $k_{n_1+1}\subseteq {K_1}^{p^{-e}}\cap K$.
\sk

{\bf Cas g\'{e}n\'{e}ral.}  Soit $s\geq 1$,  par application successive de l'hypoth\`ese de r\'ecurrence, on obtient ${K_s}^{p^{-n_1}}\cap K\subseteq {K_{s+1}}^{p^{-n_2}}\cap K\subseteq \dots \subseteq {K_{s+p_1-1}}^{p^{-n_{p_1}}}\cap K\subseteq K_{s+p_1}$. Dans la suite, on d\'esire prouver que ${K_s}^{p^{-n_1-1}}\cap K\subseteq {K_{s+1}}^{p^{-e}}$.  D'apr\`es le lemme {\ref{lem16}}, on a ${K_s}^{p^{-n}}\cap K=K_s(\te_{n_1}^{s+1},\te_{n_2}^{s+2},\dots,\te_{n_{p_1}}^{s+{p_1}})$. Soit $h_0=Q(X^{p^{-\infty}}, ({(Z_{i}^{j})}^{p^{-\infty}})_{(i,j)\in \Gamma_{s+1}\setminus \{(1,s+1)\}}, {(Y_{i}^{j})}^{p^{-\infty}})_{(i,j)\in \Gamma_{s+1}})$, o\`{u} $\Ga_{s+1}={\N}^* \times \{2,$ $\dots,$ $ s+1\}$. Il est clair que $h_0$ est parfait. Posons \'{e}galement, $h=k(h_0)$. D'apr\`es le lemme {\ref{lem10}}, pour tout $(i,j)\in \Ga_s$, on a $\te_{i}^{j}\in h_0$. En outre, $K_s\subseteq h$. De plus, comme pour tout $i\geq 1$, $\te_{i}^{s+1}={(Z_{i}^{s+1})}^{p^{-1}}\te_{2i}^{s}+\dots+ {(Z_{1}^{s+1})}^{p^{-i}}{(\te_{2}^{s})}^{p^{-i+1}}+{(Y_{i}^{s+1})}+\dots+ {(Y_{1}^{s+1})}^{p^{-i}}$, on en d\'eduit que $h_0(\te_{i}^{s+1})=h_0({(Z_{1}^{s+1})}^{p^{-i}})$. En modifiant l\'{e}g\`{e}rement les notations ci-dessus, on va se ramener au conditions de l'exemple ci-haut.  Pour se faire,  pour tout $i\geq 1$, on pose :
$H_i=h(K_{s+i})$, et par convention $H_0=h$.  Notons \'{e}galement : $H=h(K)=\di\bigcup_{i\geq 1} H_i$,  $X'=Z_{1}^{s+1}$,
et pour tout $j\geq 2$, ${Z'}_{i}^{j}=Z_{i}^{j+s}$ et ${Y'}_{i}^{j}=Y_{i}^{j+s}$.
On v\'{e}rifie aussit\^{o}t que
$H_j=\di\bigcup_{i\geq 1} H_{j-1} (\be_{i}^{j})$, o\`u $\be_{i}^{j}=\te_{i}^{j+s}={(Z_{i}^{j+s})}^{p^{-1}}\te_{2i}^{j+s-1}+{(\te_{i-1}^{j+s})}^{p^{-1}}+{(Y_{i}^{j+s})}^{p^{-1}}$ pour tout $j\geq 1$. En particulier, pour tout $j\geq 2$, $\be_{i}^{j}={({Z'}_{i}^{j})}^{p^{-1}}\be_{2i}^{j-1}+{(\be_{i-1}^{j})}^{p^{-1}}+{({{Y'}_{i}}^{j})}^{p^{-1}}$ ; et par suite on se retrouve de nouveau dans les m\^{e}mes conditions de l'exemple pr\'{e}c\'{e}dent.
\sk

Le reste de la d\'emonstration r\'esulte aussit\^ot du lemme suivant :
\begin{lem} {\label{lem17}}Pour tout $(i,j)\in \N^*\times \N$, $o(\be_{i}^{j+1},H_j)=o(\te_{i}^{s+j+1},h(K_{s+j}))=i=o(\te_{i}^{s+j+1},K_{j+s})$. En outre,  $K$ et $h$ sont $K_s$-lin\'eairement disjointes.
\end{lem}
\pre Comme $(X,({Z}_{i}^{j},{Y}_{i}^{j})_{(i,j)\in \N^*\times (\N\setminus\{0,1\})}$ est alg\'ebriquement libre sur $Q$ et $h_0=Q(X^{p^{-\infty}}, $ $({(Z_{i}^{j})}^{p^{-\infty}})_{(i,j)\in \Gamma_{s+1}\setminus \{(1,s+1)\}}, {(Y_{i}^{j})}^{p^{-\infty}})_{(i,j)\in \Gamma_{s+1}})$, alors
la famille $(X',({Z'}_{i}^{j},{Y'}_{i}^{j})_{(i,j)\in \N^*\times (\N\setminus\{0,1\})})$ est alg\'ebriquement ind\'ependante sur $h_0$. D'o\`u, de la m\^eme fa\c{c}on  qu'au lemme {\ref{lem11}},
 pour tout $(i,j)\in \N^*\times \N$, on a $o(\be_{i}^{j+1},H_j)=o(\te_{i}^{s+j+1},h(K_{s+j}))=i=(\te_{i}^{s+j+1},K_{j+s})$,
 ou encore pour tout $(i,j)\in \N^*\times \N$, $H_{j}$ et $K_{s+j}(\be_{i}^{j+1})$ sont $K_{s+j}$-lin\'eairement disjointes. Comme, $K_{j+s+1}$ est
 r\'eunion croissante d'extensions $(K_{j+s}(\be_{i}^{j+1}))_{i\geq 1}$, on en d\'eduit que $K_{j+s+1}$ et $H_j$ sont $K_{s+j}$-lin\'eairement disjointes pour tout $j\in \N$.
Par application successive de la transitivit\'e de la lin\'earit\'e disjointe, on aura $K_{j+s+1}$ et $H_0$  sont $K_{s}$-lin\'eairement disjointes. Or, $K$ est r\'eunion croissante des sous-extensions ${(K_n)}_{n\geq s}$, donc $K$ et $h$ sont $K_s$-lin\'eairement disjointes. \cqfd \sk

{\bf Suite de la preuve du th\'eor\`eme {\ref{thm24}}.}  Pour des raisons d'\'{e}criture, on pose $N= h(\be_{n_1}^{1},\dots,\be_{n_{p_1}}^{p_1})$. On montre \`{a}  la m\^{e}me fa\c{c}on du lemme {\ref{lem15}}  que $o_2({N}^{p^{-1}}\cap H/h)=e$. Or, ${K_s}^{p^{-n_1}}\cap K=K_s(\be_{n_1}^{1},\dots,\be_{n_{p_1}}^{p_1})\subseteq N$, donc ${K_s}^{p^{-n_1-1}}\cap K\subseteq {N}^{p^{-1}}\cap H$.
Comme $K$ et $h$ sont $K_s$-lin\'eairement disjointes, alors en particulier,   $h({K_s}^{p^{-n_1-1}}\cap K)\simeq h\otimes_{K_s} {K_s}^{p^{-n_1-1}}\cap K$ ; et par suite d'apr\`es les propositions {\ref{pr16}}, {\ref{pr19}}, $o_2(h({K_{s}}^{p^{-n_1-1}}\cap K)/h)=o_2( h\otimes_{K_s} {K_{s}}^{p^{-n_1-1}}\cap K/h)=o_2({K_{s}}^{p^{-n_1-1}}\cap K/K_s)\leq o_2({N}^{p^{-1}}\cap H/h)=e$. Comme $o(\te_{n_1+1}^{s+1},K_s)=n_1+1$, par application de l'algorithme de  la compl\'{e}tion des $r$-bases, on compl\`{e}te $\te_{n_1+1}^{s+1}$ en une $r$-base canoniquement ordonn\'{e}e de ${K_s}^{p^{-n_1-1}}\cap K/K_s$ ; et donc $o_2({K_s}^{p^{-n_1-1}}\cap K/K_s)=o_1({K_s}^{p^{-n_1-1}}\cap K/K_s(\te_{n_1+1}^{s+1}))\leq e$. Par cons\'{e}quent, pour tout $x\in {K_s}^{p^{-n_1-1}}\cap K$, $x^{p^{e}}\in K_s(\te_{n_1+1}^{s+1})\subseteq K_{s+1}$ ; il en r\'{e}sulte que  $x\in {K_{s+1}}^{p^{-e}}\cap K$. D'o\`{u} ${K_{s}}^{p^{-n_1-1}}\cap K\subseteq {K_{s+1}}^{p^{-e}}\cap K$. \cqfd
\begin{cor} {\label{cor17}} Pour tout  $(n,s)\in N^*\times \N$,  ${K_s}^{p^{-n}}\cap K=K_s(\te_{n_1}^{s+1},\te_{n_2}^{s+2},\dots,$ $\te_{n_{p_1}}^{s+p_1})$, o\`u $n_1=n$, $lpi(n)=p_1$, et pour tout $s\in \N^*$, $n_{s+1}=\di[\frac{n_s}{2}]$.
 \end{cor}
\pre D\'emonstration analogue \`a celle du lemme {\ref{lem16}}. \cqfd \sk

Voici une liste de  cons\'{e}quences du th\'{e}or\`{e}me {\ref{thm22}},
\begin{thm} {\label{thm25}} La plus petite sous-extension $m/k$ de $K/k$ telle que
$K/m$ est modulaire est triviale. Autrement dit $mod(K/k)=K$
\end{thm}
\pre Notons $m=mod(K/k)$. Si $m\not=K$, en vertu du th\'{e}or\`{e}me {\ref{thm22}}, il existe $i\in \N^*$ tel que $K_i\subseteq m\subseteq K_{i+1}$ et $m/K_i$ est fini. Comme $di(K/K_{i+1})$ est infini, il en est de m\^{e}me de $di(K/m)$.  En vertu du corollaire {\ref{acor1}}, $di(m^{p^{-1}}\cap K/m)=di(K/m)$, donc $di(m^{p^{-1}}\cap K/m) $ est infini. Par ailleurs, il existe  $e\in \N$, tel que $m\subseteq {K_i}^{p^{-e}}\cap K$, ($m/K_i$ est finie), et donc $m^{p^{-1}}\cap K\subseteq {K_i}^{p^{-e-1}}\cap K$. Compte tenu du th\'{e}or\`{e}me {\ref{thm4}}, on a $di(m^{p^{-1}}\cap K /m)\leq di({K_i}^{p^{-e-1}}\cap K/K_i) <+\infty$, c'est une contradiction avec le fait que $ di(m^{p^{-1}}\cap K /m)$ est infini. D'o\`{u} $m=K$.\cqfd
\begin{cor} {\label{cor18}} Pour tout $j\geq 1$, $mod(K/K_j)=K$.
\end{cor}
\pre Imm\'ediat. \cqfd
\begin{rem} Soit $\Om$ une cl\^{o}ture alg\'{e}brique de $K$. Dans $[k :\Om]$ on d\'{e}finit la relation $\sim$ de la fa\c{c}on suivante : $K_1\sim K_2$ si $K_1\subseteq K_2$ et $K_2/K_1$ est finie ou $K_2\subseteq K_1$ et $K_1/K_2$ est finie.
On v\'{e}rifie aussit\^{o}t que $\sim$ est r\'{e}flexive, sym\'{e}trique, cependant $\sim$ est g\'{e}n\'{e}ralement non transitive. Il est imm\'ediat que
\sk

\begin{enumerate}{\it
\item[\rm{(1)}] $\{K_i/k, \, i\geq 1\}$ est l'ensemble des repr\'{e}sentants des sous-extensions d'expo\-s\-ant non born\'{e} de $K/k$ pour la relation $\sim$.
\item[\rm{(2)}] $K/k$ et $(K_i/k)_{i\in \N}$ sont les seules sous-extensions relativement parfaites de $K/k$.
\item[\rm{(3)}] Pour tout $i,j \in \N$ tels que $i<j$, $mod(K_j/K_i)=K_{j-1}$.}
\end{enumerate}
\end{rem}
\subsection{Caract\'erisation des extensions $q$-finies}
\begin{pro} {\label{pr37}} Toute extension absolument $lq$-finie $K/k$ est $q$-finie sur $mod($ $K/k)$. En particulier, toute extension qui est \`a la fois modulaire et absolument $lq$-finie est $q$-finie
\end{pro}
\pre Notons $m=mod(K/k)$.  Le cas $m=K$ est trivial, pour cela on suppose que $m\not =K$. D'apr\`{e}s la proposition {\ref{pr23}}, $di(m^{p^{-1}}\cap K/m)=di(K/m)$. Comme $K/k$ est absolument $lq$-finie, il en est de m\^{e}me de $K/m$, et donc $di(K/m)$ est fini. \cqfd
\begin{thm} {\label{thma1}} Une extension  absolument $lq$-finie $K/k$ est $q$-finie si et seulement si pour tout corps interm\'{e}diaire $L$ de $K/k$, la plus petite sous-extension $m$ de $L/k$ telle que $L/m$ est modulaire est non triviale ($m\not= L$).
\end{thm}
\pre
La condition n\'{e}cessaire r\'{e}sulte du th\'{e}or\`{e}me {\ref{thm27}}. Inversement,  consid\'{e}rons la suite d\'{e}croissante $(m_i)$  de sous-extensions de $K/k$ d\'{e}finie par : $m_0=K$, et pour tout $i\geq 1$, $m_{i}=mod(m_{i-1}/k)$. Par hypoth\`{e}se, si $m_{i-1}\not=k$, alors $m_i\not=m_{i-1}$. Toutefois, comme toute suite d\'{e}croissante de sous-extensions de $K/k$ est stationnaire, on en d\'{e}duit l'existence d'un entier $j$ tel que $m_j=k$. Par suite, on aura $K=m_0 \longleftarrow m_1 \longleftarrow\dots \longleftarrow m_{j-1} \longleftarrow m_j=k$. En particulier, pour tout $i\in \{1,\dots,j\}$,  $m_{i-1}/m_i$ est absolument $lq$-finie, et par suite d'apr\`{e}s  la proposition {\ref{pr37}},   $m_{i-1}/m_i$  est $q$-finie. D'o\`u,  $K/k$ est $q$-finie. \cqfd \sk

Comme cons\'{e}quence imm\'{e}diate  des th\'{e}or\`{e}mes {\ref{thma1}} et {\ref{thm19}}, on a
\begin{cor}
Pour qu'une extension purement ins\'{e}parable $K/k$ soit $q$-finie il faut et il suffit qu'elle  satisfait les deux conditions suivantes :
\sk

\begin{itemize}{\it
\item[\rm{(i)}] Toutes  suite d\'{e}croissante de sous-extensions de $K/k$ est stationnaire.
\item[\rm{(ii)}] Pour tout corps interm\'{e}diaire $L$ de $K/k$, on a $mod(L/k)\not= L$.}
\end{itemize}
\end{cor}

Par ailleurs, voici une propri\'et\'e caract\'eristique qui permet d'identifier les extensions moduliares qui sont $q$-finies.
\begin{thm} {\label{thm20}} Soit $K/k$ une extension purement ins\'{e}parable et modulaire. Les assertions suivantes sont \'equivalentes :
\sk

\begin{itemize}{\it
\item[\rm{(1)}] Toute suite d\'{e}croissante de sous-extensions de $K/k$ est stationnaire.
\item[\rm{(2)}] Toute ensemble de sous-extensions  de $K/k$ admet un \'element minimal.
\item[\rm{(3)}] $K/k$ est $q$-finie.}
\end{itemize}
\end{thm}
\pre D\'ecoule aussit\^ot du th\'eor\`eme {\ref{thm19}} et la proposition  {\ref{pr37}}.\cqfd \sk

Comme application imm\'ediate, on a :
\begin{cor}Soit $k$ un corps commutatif de caracterstique $p>0$, et $\Omega$ une cl\^oture alg\'ebrique de $k$. Il est \'equivalent de dire que :
\sk

\begin{itemize}{\it
\item[\rm{(1)}] $di(k)$ est fini.
\item[\rm{(2)}] Toute ensemble de sous-extensions purement ins\'eparables de $\Omega/k$ admet un \'element minimal.
\item[\rm{(3)}] Toute suite d\'ecroissante de sous-extensions purement ins\'eparables de $\Omega/k$ est stationnaire.}
\end{itemize}
\end{cor}
\pre Il suffit de remarquer que la cl\^oture purement ins\'eparable $k^{p^{-\infty}}$ de $\Omega/k$ est modulaire sur $k$, et $di(k)=di(k^{p^{-\infty}}/k)$. Par suite, le r\'esulat d\'ecoule aussit\^ot du th\'eor\`eme {\ref{thm20}}. \cqfd \sk

Soit $K/k$ une extension modulaire d'exposant non born\'e. Pour  que $K/k$   soit $q$-finie, il suffit de modifier l\'eg\'erement la condition suffisante du th\'eor\`eme ci-dessus de la fa\c{c}on suivante :
\begin{thm} {\label{thm21}} Pour qu'une extension modulaire d'exposant non born\'e $K/k$ soit $q$-finie il faut et il suffit que toute suite d\'ecroissante de sous-extensions d'exposant non born\'e de $K/k$ soit stationnaire.
\end{thm}

Pour la preuve nous aurons besoin du r\'esultat suivant :
\begin{lem} {\label{lem9}} Soit $K/k$ une extension purement ins\'{e}parable d'exposant non born\'{e} et de degr\'e d'irratinalit\'e infini. Si de plus, $K/k$ est relativement parfaite et modulaire, alors $K/k$ contient une sous-extension propre $L/k$  d'exposant non born\'{e} et modulaire.
\end{lem}
\pre On va construire par r\'{e}currence une suite strictement croissante $(K_n/k)_{n\geq 1}$ de sous-extensions modulaires d'exposant $n$ de $K/k$. Comme $K/k$ est relativement parfaite, d'apr\`{e}s la proposition {\ref{thm13}} et le corollaire {\ref{acor1}}, pour tout $n\geq 1$, $ di(k^{p^{-n}}\cap K/k)=di(k^{p^{-1}}\cap K/k)=di(K/k)$ et $k^{p^{-n}}\cap K/k$ est \'equiexponentielle d'exposant $n$. Soit $G_1$ une $r$-base de $k^{p^{-1}}\cap K/k$, il en r\'{e}sulte que  $k^{p^{-1}}\cap K\simeq \otimes_k (k(a))_{a\in G_1}$. Choisissons un \'{e}l\'{e}ment $x$ de $G_1$, comme $G_1$ est infini,  il existe un sous-ensemble fini $G'_1$ de $G_1$ tel que $x \not \in k(G'_1)$. Posons $K_1=k(G'_1)$, il est clair que $K_1/k$ est modulaire. Supposons qu'on a construit une suite de sous-extensions finies  $k\subseteq K_1\subseteq K_2\subseteq \dots K_n$ de $K/k$ telle que
\sk

\begin{itemize}{\it
\item[\rm{(1)}] Pour tout $i\in \{1,\dots,n\}$, $K_i/k$ est modulaire.
\item[\rm{(2)}] Pour tout $i\in \{1,\dots,n\}$, $o_1(K_i/k)=i$.
\item[\rm{(3)}] $x\not\in K_n$.}
\end{itemize}
\sk

Soit $G_{n+1}$ une $r$-base de $k^{p^{-n-1}}\cap K/k$, d'apr\`{e}s la proposition {\ref{pr26}}, $k^{p^{-n-1}}\cap K\simeq \otimes_k (k(a))_{a\in G_{n+1}}$. Comme $o_1(K_{n}/k)=n$, on en d\'{e}duit que $K_n\subseteq k^{p^{-n-1}}\cap K$. Or $K_n/k$ est finie et $G_{n+1}$ est infinie, donc il existe une partie finie $G'_{n+1}$ de $G_{n+1}$ telle que $K_n\subseteq k(G'_{n+1})$. Deux cas peuvent se produire :

1-ier cas si $x\not\in k(G'_{n+1})$, alors $K_{n+1}=k(G'_{n+1})$ convient.
\sk

2-i\`{e}me cas si $x\in k(G'_{n+1})$, comme $k^{p^{-n-1}}\cap K\simeq \otimes_k (\otimes_k (k(a))_{a\in G'_{n+1}})\otimes_k (\otimes_k (k(a))_{a\in G_{n+1}\setminus G'_{n+1}})$, donc $x\not\in k(G_{n+1}\setminus G'_{n+1})$ ; sinon puisque $k(G'_{n+1})$ et $k(G_{n+1}\setminus G'_{n+1})$ sont $k$-lin\'{e}airement disjoints, alors $x \in k(G'_{n+1})\cap k(G_{n+1}\setminus G'_{n+1})=k$, absurde.  Soit $y$ un \'{e}l\'{e}ment  de $G_{n+1}\setminus G'_{n+1}$, ($y$ existe car $G_{n+1}$ est infinie et $G'_{n+1}$ est finie).  Notons $K_{n+1}=K_n(y)$, on v\'{e}rifie aussit\^ot que :
\sk

\begin{itemize}{\it
\item $K_{n+1}/k$ est finie, et $o_1(K_{n+1}/k)=o(y,k)=n+1$.
\item $K_{n+1}\simeq K_n\otimes_k k(y)$, (application de la transitivit\'{e} de la lin\'{e}arit\'{e} disjointe de $k(G'_{n+1})$ et $k(G_{n+1}\setminus G'_{n+1})$), et comme $K_n/k$ est modulaire, d'apr\`{e}s la proposition {\ref{apr1}}, $K_{n+1}/k$ est modulaire.
\item $x \not\in K_{n+1}$, sinon comme $k^{p^{{-n-1}}}\cap K\simeq k(G'_{n+1})\otimes_k k(G_{n+1}\setminus G'_{n+1})\simeq K_n(G'_{n+1})\otimes_{K_n} K_n(G_{n+1}\setminus G'_{n+1})$, alors $x\in k(G'_{n+1})\cap K_n(y)\subseteq K_n(G'_{n+1})$ $\cap K_n (G_{n+1}\setminus G'_{n+1})=K_n$, abbsurde.}
\end{itemize}
\sk

D'o\`u $K_{n+1}/k$ convient, et par suite $L=\di \bigcup_{i\geq 1}K_i$ satisfait les conditions du th\'{e}or\`{e}me ci-dessus.\cqfd
\sk

{\bf Preuve du th\'eor\`eme {\ref{thm21}}} Il suffit de montrer que $2\Rightarrow 1$. Pour cela, on va utiliser un raisonnement par contrapos\'e.
Supposons que $K/k$ n'est pas $q$-finie, et soit $G$ une $r$-base de $K/k(K^p)$. Si $|G|$ est infinie, il existe  une famille $(a_i)_{i\in \N^*}$ d'\'{e}l\'{e}ments de $G$ deux \`{a} deux distincts. Pour tout $n\in N^*$, notons $K_n=k(K^p)(G\setminus \{a_1,\dots,a_n\})$, et $K_0=K$. Il est clair que  la suite $(K_n)_{n\in N}$ de sous-extensions d'exposant non born\'e de $K/k$ est d\'{e}croissante, mais comme $G$ est une $r$-base de $K/k(K^p)$, alors pour tout $n\in\N$, $K_n\not= K_{n+1}$ ; et par suite la suite $(K_n/k)_{n\in \N}$ est non stationnaire, donc on est amen\'e au cas o\`u $|G|$ est fini. Posons $e=o_1(k(G)/k)$ et $L_1=k(K^{p^e})$. D'une part, d'apr\`es la proposition {\ref{pr23}}, $L_1/k$ est modulaire. D'autre part, Comme $K=k(K^p)(G)$, alors $k(K^{p^e})=k(K^{p^{e+1}})(G^{p^e})=k(K^{p^{e+1}})$ ; ou encore $L_1/k$ est relativement parfaite. On v\'{e}rifie aussit\^{o}t que :
\sk

\begin{itemize}{\it
\item[$\bullet$]  Si $L_1/k$ est $q$-finie, il en est de m\^{e}me de $K/k$ puisque $di(K/k)\leq di(K/L_1)$ $+di(L_1/k)$.
\item[$\bullet$] Si $L_1/k$ n'est pas $q$-finie, d'apr\`{e}s le lemme {\ref{lem9}}, $L_1/k$ contient une sous-extension propre modulaire d'exposant non born\'{e} $L_2/k$, et on distingue deux cas :}
\end{itemize}
\sk

1-ier cas : Si $L_2/k$ est $q$-finie, comme $L_1/k$ n'est pas $q$-finie et modulaire, d'apr\`{e}s le corollaire {\ref{acor1}}, $di(k^{p^{-1}}\cap L_1 /k)$ est infini. Notons $M=L_2(k^{p^{-1}}\cap L_1)$. Il est imm\'{e}diat que $di(M/k)\geq di(k^{p^{-1}}\cap L_1/k)$, et comme $L_2/k$ est $q$-finie, alors $di(M/L_2)$ est infini. Par suite si $G$ est une $r$-base de $M/L_2$, alors $|G|$ est infini. Comme dans le cas pr\'ec\'edent, on construit une suite strictement d\'ecroissante de sous-extensions de $K/k$.
\sk

2-i\`{e}me cas si $L_2/k$ n'est pas $q$-finie, on se ram\`{e}ne aux conditions de $K/k$, et en r\'{e}p\'{e}tant les m\^{e}me proc\'{e}dures, on construit une suite d\'{e}croissante non stationnaire  $(L_n/k)_{n\geq 1}$ de sous extensions de $K/k$.\cqfd
\subsection{L'absolument $lq$-finitude et la $w_0$-g\'en\'eratrice}

Rappelons qu'une extension $K/k$ est dite $w_0$-g\'en\'eratrice, s'elle n'admet aucune sous-extension propre d'exposant non born\'e. Pour plus de d\'etails sur ce sujet, se r\'ef\'erer aux articles \cite{Dev2} et \cite{Che-Fli2}.
\begin{thm} {\label{thm18}}Toute extension absolument $lq$-finie est compos\'{e}e d'\'{e}xtensions $w_0$-g\'{e}n\'{e}ratrices.
\end{thm}
\pre Si $K/k$ est finie ou $w_0$-g\'{e}n\'{e}ratrice, c'est termin\'{e}. Sinon, soit $H$ l'ensemble de sous-extensions d'exposant non born\'{e} de
$K/k$. D'apr\`{e}s la proposition {\ref{pr31}}, $H$ admet  une sous-extension minimale d'exposant non born\'{e} de $K/k$  que l'on note $K_1$.
N\'{e}cessairement, $K_1/k$ est $w_0$-g\'{e}n\'{e}ratrice.  Supposons qu'on a construit une suite de sous-extensions de $K/k$ telle que
$k=K_0\subseteq K_1\subseteq \dots\subseteq K_n$, et pour tout $i\in \{1,\dots, n\}$,  $K_{i}/K_{i-1}$ est $w_0$-g\'{e}n\'{e}ratrice.
Si $K=K_n$, on s'arr\^ete. Sinon, Comme $K/K_n$ est absolument $lq$-finie, on reprond avec $K/K_n$, et donc il existe une sous extension $K_{n+1}/K_n$ de $K/k_n$ telle que $K_{n+1}/K_n$  est $w_0$-g\'{e}n\'{e}ratrice.  D'o\`u, compte tenu de la propri\'et\'e de r\'ecurrence, il existe une suite de sous-extensions $(K_i/k)$ de $K/k$ telle que  $k=K_0\subseteq K_1\subseteq \dots\subseteq K_n\subseteq\dots\subseteq K$
et pour tout $i\geq 1$, $K_i/K_{i-1}$ est $w_0$-g\'{e}n\'{e}ratrice. \cqfd \sk

Comme application de la proposition {\ref{pr37}}, on a :
\begin{pro} {\label{cor13}} Soit $K/k$  une extension absolument $lq$-finie, si de plus $K/k$ est $w_0$-g\'{e}n\'{e}ratrice, alors $mod(K/k)\not=K$  si et seulement si $K/k$ est $q$-finie.
\end{pro}
\pre \smartqed Imm\'ediat. \cqfd
{\small
{\em Authors' addresses}:
{\em EL Hassane Fliouet}, Regional Center for the Professions of Education and Training, Agadir, Morocco
 e-mail: \texttt{fliouet@yahoo.fr}.
 }
\end{document}